\documentclass[3p]{elsarticle}
\usepackage{graphicx}
\usepackage{amssymb}
\usepackage{stmaryrd}
\usepackage{dsfont}
 \usepackage{amsthm}
  \usepackage{amsmath}
\usepackage{tikz-cd}

 \usepackage{float}
 \usepackage{newfloat}
\usepackage{comment}
\usepackage{ccicons}
\usepackage{epstopdf}

\usepackage{multirow}
\usepackage{color}
\usepackage{hhline}
\usepackage{booktabs}
\usepackage{rotating}
\usepackage{hyperref}
\usepackage{bm}

\definecolor{nverde}{RGB}{0,61,0} %nverde
\definecolor{cr1}{RGB}{200,0,0}
\definecolor{cr2}{RGB}{0,0,200}
\definecolor{cr12}{RGB}{100,0,100}

\newcommand{\halb}{\frac{1}{2}}
\newcommand{\n}{\mathbf{n}}
\renewcommand{\v}{\mathbf{v}}

\newcommand{\A}{\mathbf{A}}
\newcommand{\G}{\mathbf{G}}
\newcommand{\J}{\mathbf{J}}

\newcommand{\x}{\mathbf{x}}
\newcommand{\f}{\mathbf{f}}
\newcommand{\q}{\mathbf{q}}

\newcommand{\lnpc}{l^{pc} \n^{pc}} % lpc npc
\newcommand{\lnpb}{l^{pb} \n^{pb}} % lpc npc

\newcommand{\devG}{\mathring{G}}
\newcommand{\ISigma}[1]{\left(\Sigma^{n,c}_{#1}\right)^{-1}}

\journal{Journal of Computational Physics}

\allowdisplaybreaks

\begin{document}

\begin{frontmatter}

\title{A structure-preserving staggered semi-implicit four-split finite volume scheme for continuum mechanics on unstructured meshes}

\author[UniTN,SUSTech]{M. Dumbser \corref{cor1}}
\ead{michael.dumbser@unitn.it}
\cortext[cor1]{Corresponding author}

\author[UniSMB,UniFE]{W. Boscheri}
\ead{walter.boscheri@univ-smb.fr}

\author[UniVR]{M. Tavelli}
\ead{maurizio.tavelli@univr.it}

\author[Inria,UniTN]{A. Thomann}
\ead{andrea.thomann@inria.fr}

\address[UniTN]{Laboratory of Applied Mathematics, DICAM, University of Trento, via Mesiano 77, 38123 Trento, Italy}

\address[SUSTech]{Department of Mathematics, Southern University of Science and Technology, \\ Xueyuan Avenue 1088, 518055 Shenzhen, Guangdong, China}

\address[UniSMB]{Laboratoire de Math\'ematiques UMR 5127 CNRS
	Universit\'e Savoie Mont Blanc,
	73376 Le Bourget du Lac, France}

\address[UniFE]{Department of Mathematics and Computer Science, University of
	Ferrara, 44121 Ferrara, Italy}

\address[UniVR]{Department of Engineering for Innovation Medicine, University of Verona, Strada le Grazie 15, 37134, Verona, Italy}

\address[Inria]{Universit\'e de Strasbourg, CNRS, Inria, IRMA, Strasbourg, F-67000, France}

\begin{abstract}
We present a new semi-implicit structure-preserving (SP) finite volume discretization for the unified first-order hyperbolic model of continuum mechanics developed over the years by Godunov, Peshkov and Romenski, hereafter referred to as the GPR model. This framework provides a common mathematical description for fluids and solids and incorporates transport, viscous effects, heat conduction, and elastic deformations within a single system of hyperbolic partial differential equations. The coexistence of several physical mechanisms gives rise to multiple characteristic wave speeds, resulting in severe time-step restrictions for fully explicit discretizations. To overcome this difficulty, we develop a four-split scheme in which the governing equations are decomposed into convective, temperature, mechanical, and pressure subsystems. The convective  subsystem is the only one that is advanced explicitly in time, while the remaining subsystems are treated implicitly in a sequential manner. As a consequence, the time step restriction of the proposed method depends only on the material velocity and is independent of the acoustic, shear, and thermal wave speeds that characterize the model.

The scheme is designed to preserve key structural properties of the underlying equations on unstructured grids. In particular, a compatible vertex-staggered spatial discretization on triangles allows the curl-free involutions associated with the distortion field and thermal impulse to be respected whenever the relaxation source terms in the governing PDE are linear or absent. At the same time, the scheme is asymptotic preserving as it is consistent with the low Mach number limit and with the stiff relaxation limits of the GPR model that recover the classical Navier-Stokes-Fourier equations in fluid mechanics. A series of numerical experiments covering both fluids and solids demonstrates the accuracy, robustness, and multi-scale capabilities of the proposed approach.

\end{abstract}


\end{frontmatter}

%\linenumbers   % arXiv: no line numbers
%\tableofcontents

% % % % % % % % % % % % % % % % % % % % % % % % % % % % % %
%         Introduction
% % % % % % % % % % % % % % % % % % % % % % % % % % % % % %

\section{Introduction}\label{sec.intro}
Time-dependent nonlinear systems of partial differential equations (PDE) are widely used to model a broad range of physical phenomena, including fluid and solid mechanics, gas dynamics, multiphase flows, and plasma dynamics, to name just a few. These systems are often characterized by the \textit{presence of multiple time scales}, each associated with a different physical process. For example, diffusion typically evolves on different time scales than advection, and in low Mach number flows pressure waves propagate significantly faster than material interfaces or contact discontinuities. This implies that the governing equations exhibit stiffness, which is driven by the ratio of the fastest to the slowest physical processes, such as for example the ratio of the acoustic wave speed to the material velocity in low Mach number flows. Consequently, numerical schemes capable of accurately and efficiently handling multiple time scales simultaneously are essential for high-fidelity simulations of real-world applications. Another distinctive feature of many hyperbolic PDE models is the \textit{presence of stationary differential constraints}, so-called involutions, which remain automatically valid throughout the entire time evolution of the system if they are satisfied by the initial data. Classical examples include divergence-free or curl-free involutions, hence requiring the satisfaction, at the discrete level, of the fundamental vector calculus identities $\nabla \cdot (\nabla \times \mathbf{a})=0$  and $\nabla \times \nabla \psi = 0$, with $\mathbf{a} \in \mathds{R}^3$ a generic vector field and $\psi \in \mathds{R}$ a scalar field. Numerical methods that are designed and proven to preserve these constraints at the discrete level are commonly referred to as structure-preserving methods. An additional characteristic of many systems of PDE is the \textit{presence of relaxation mechanisms}, which are mathematically represented by source terms in the governing equations and are used to model non-equilibrium dynamics. These source terms are typically characterized by a stiffness parameter $\tau$, commonly referred to as the relaxation time. In the stiff asymptotic limit $\tau \to 0$, the original model approaches an equilibrium or a reduced model. For example, in the Baer-Nunziato model \cite{BaerNunziato1986,SaurelAbgrall}, pressure and velocity relaxation drive the system toward mechanical equilibrium in the stiff relaxation limit. Numerical schemes designed for the full relaxation system that are able to recover a consistent discretization of the corresponding limit model, while employing a time-marching strategy independent of the relaxation parameter $\tau$, are called asymptotic-preserving schemes.

All the aforementioned properties rarely coexist within the same PDE system, since each mathematical model typically describes only a subset of physical processes. For example, the compressible Navier-Stokes equations describe viscous heat-conducting flows, plasma dynamics in the continuum limit are described by the magnetohydrodynamics (MHD) equations, while hyperelasticity models are used in nonlinear solid mechanics. As a consequence, the associated numerical schemes have also been developed with the aim of simultaneously satisfying only a subset of the corresponding structural properties, depending on the specific model under consideration. However, in \cite{PeshRom2014}, a unified first-order hyperbolic model of continuum mechanics was proposed, building upon the class of hyperbolic thermodynamically compatible (HTC) systems \cite{God1961,Rom1998}. This model, commonly referred to as the GPR model, is able to describe both fluid and solid mechanics through relaxation terms in the evolution equations of the distortion field, which govern the mechanical response of the medium. In the stiff relaxation limit, either inviscid or viscous fluid behavior is recovered, thus reducing the model to the well-known Euler or compressible Navier-Stokes equations. Moreover, heat conduction is modeled through a first-order hyperbolic evolution equation for the thermal impulse supplemented by relaxation source terms, which recover the Fourier law in the stiff limit. The homogeneous GPR system furthermore possesses two curl-type involutions associated with the distortion field and the thermal impulse. The GPR model is inherently multiscale, as it embeds transport, diffusion, shear deformation, and heat-conduction effects within one single mathematical model, in addition to featuring involutions and stiff relaxation terms. Therefore, the numerical discretization of this model poses the challenge of simultaneously preserving a large number of additional structural properties which concern the stiff relaxation source terms, the involutions as well as the low Mach number limits.

The numerical treatment of multiscale PDE systems remains a major challenge due to the coexistence of physical processes evolving on very different time scales. A widely adopted strategy consists in identifying the mechanisms responsible for the stiffness of the governing equations and separating them from the non-stiff dynamics. This decomposition permits the use of different discretization techniques for the resulting subsystems: the fast processes are typically integrated implicitly in order to avoid severe stability restrictions on the time step, whereas the slow dynamics can be advanced explicitly to preserve computational efficiency and resolution. Such an approach naturally leads to the solution of coupled subsystems associated with the different temporal scales present in the original model. Among the available methodologies, semi-implicit schemes constitute one of the most established classes of multiscale integrators \cite{MunzPark,Degond2,KleinMach,Casulli1984,Casulli1990,BosFil2016,DumbserCasulli2016,TavelliDumbser2017}. Their key idea is to treat the stiff contributions through suitable linearizations, yielding linearly implicit discretizations that retain the stability advantages of fully implicit methods while significantly reducing the computational cost associated with the solution of nonlinear systems \cite{BP2021}. A more general framework for the time integration of multiscale systems is provided by partitioned methods, among which implicit-explicit (IMEX) Runge-Kutta schemes represent one of the most successful and widely adopted approaches \cite{AscRuuSpi,BP2017,BosRus,PR_IMEX,BPR2017,Hofer,Thomann2020}. By combining implicit and explicit discretizations within the same time-stepping procedure, IMEX methods achieve high-order accuracy while avoiding the severe stability restrictions associated with the fastest scales of the problem. Moreover, suitably designed stiffly accurate IMEX schemes have been shown to possess the asymptotic-preserving (AP) property, ensuring that the asymptotic limit of the governing equations is consistently recovered at the discrete level without requiring the time step to resolve the fast dynamics \cite{JINAP1999,JP2001,KLARAP1999}. Most semi-implicit and IMEX methods are designed to handle two distinct time scales and, consequently, two associated subsystems. For example, in the case of the Euler equations, the advection subsystem is treated explicitly, while the pressure subsystem is discretized implicitly according to the flux splitting introduced in \cite{MunzPark,ToroVazquez}. However, many PDE models exhibit more than two relevant time scales, requiring the original system to be decomposed into a larger number of subsystems. For instance, three subsystems naturally arise when the viscous terms of the compressible Navier-Stokes equations are separated from the advection and pressure contributions, leading to the three-split schemes proposed in \cite{TavelliDumbser2017,BosTav22}. Another recent example is provided by the three-split formulations for the magnetohydrodynamics equations introduced in \cite{Fambri_3splitMHD,BosTho_3splitMHD}, where advection, magnetic, and pressure effects are treated as separate subsystems.

\vspace{2mm} 

The natural objective of this work is to design a novel numerical scheme on unstructured simplex meshes that can deal simultaneously with the multiple time scales, respect involutions and satisfy the asymptotic limits of the GPR model for continuum mechanics at the discrete level. Four subsystems naturally arise for this model: i) convection; ii) temperature (heat conduction); iii) deformation gradient and velocity with mechanical and thermal stress; iv) pressure. A four-split scheme for the GPR model has been very recently proposed in \cite{FourSplit}, where the first subsystem associated with transport phenomena is solved explicitly, while the remaining three subsystems are discretized implicitly. Consequently, the time step restriction of the method is dictated by the classical CFL condition based only on the velocity field of the medium and not by the other fast waves present in the system. Here, we aim to retain the same time-marching algorithm while extending the spatial discretization to unstructured triangular meshes. The resulting scheme is capable of dealing with four different time scales of the original model, and thanks to the semi-implicit time discretization, it is consistent with the stiff relaxation limits of the original GPR model.

The adoption of unstructured grids poses challenges for the satisfaction of the curl involutions related to the distortion field and the thermal impulse. The construction of discrete differential operators that exactly reproduce fundamental vector-calculus identities has received considerable attention over the last decades, motivated by applications ranging from computational electromagnetics to fluid and plasma dynamics. A large variety of approaches has been proposed in the literature. In the context of Maxwell and magnetohydrodynamics equations, compatible discretizations can be traced back to the pioneering work of Yee \cite{Yee66} and have subsequently been extended in numerous finite difference and finite volume formulations \cite{DeVore,BalsaraSpicer1999,Balsara2004,GardinerStone,balsarahlle2d,ADERdivB}. Alternative compatible frameworks include mimetic finite difference methods \cite{HymanShashkov1997,JeltschTorrilhon2006,Torrilhon2004,Margolin2000,Lipnikov2014,Carney2013}, finite element discretizations based on discrete differential forms and finite element exterior calculus \cite{Nedelec1,Nedelec2,Cantarella,Hiptmair,Monk,Arnold,Alonso2015,CAMPOSPINTO2016,Zampa1,Zampa2}, as well as more recent discontinuous Galerkin formulations \cite{SPDG2023,Ern2023,perrier2024}. Related developments relying on discrete de Rham complexes and Compatible Discrete Operators (CDO) can be found in \cite{bonelle2015,DiPietro2023}. A recurring ingredient of many compatible discretizations is the use of staggered grids. In these methods, different physical quantities are associated with different geometric entities of the mesh, such as cells, faces, edges, or vertices. Classical edge-based staggering \cite{Yee66,Margolin2000} exploits the interplay between primal and dual grids to enforce the desired algebraic identities and has recently been employed in the context of the GPR model \cite{SIGPR}. More recently, vertex-based staggering strategies have been shown to provide an effective and flexible alternative, particularly on unstructured meshes \cite{Barsukow2024,Sidilkover2025,CompatibleDG1}. Thus, in this work we rely on vertex-staggered grids, which collocate the variables differently with respect to our previous work \cite{FourSplit}. In particular, the momentum, distortion field, and thermal impulse are discretized at the cell centers of the triangular elements, while the remaining quantities are located at the mesh vertices. Finally, the proposed structure-preserving discretization allows the recovery of the expected $O(M_a^2)$ scaling in the low Mach number limit of the GPR model.
To the best of our knowledge, the method proposed here is the first all Mach number flow solver for the GPR model of continuum mechanics on unstructured meshes which also respects the two basic vector calculus identities at the discrete level. It differs from the four-split scheme of \cite{FourSplit} not only by the use of unstructured simplex meshes instead of uniform Cartesian grids, but also by the collocation of the variables, since here the momentum, the distortion field and the thermal impulse are located at the cell centers of the primal mesh, while the density, the pressure and the total energy are located at the vertices, that is in the centers of the dual mesh, so that the compatible discrete operators, which on Cartesian grids follow directly from the tensor product structure of the mesh, have to be constructed explicitly. It differs from the structure-preserving and thermodynamically compatible cell-centered Lagrangian scheme of \cite{HTCLagrangeGPR}, which is likewise formulated on unstructured meshes, by its Eulerian and semi-implicit character, since that explicit Lagrangian discretization is subject to a time step restriction governed by the acoustic, the shear and the heat wave speeds, whereas the time step of the present scheme is dictated by the bulk flow velocity alone.

The rest of the paper is organized as follows. In Section \ref{sec.pde} we briefly recall the GPR model of continuum mechanics and we present the splitting strategy with the associated subsystems. In Section \ref{sec.method} we introduce the notation for the unstructured grid, we present the new four-split structure-preserving finite volume scheme and we analyse its asymptotic limits. Computational results for a wide range of mechanical regimes are shown in Section \ref{sec.results}.
Section \ref{sec.conclusions} concludes the paper providing some remarks and an outlook to future work.

% % % % % % % % % % % % % % % % % % % % % % % % % % % % % %
%         splitting
% % % % % % % % % % % % % % % % % % % % % % % % % % % % % %

\section{Governing PDE system and splitting}\label{sec.pde}

The governing equations are given by the unified first order hyperbolic model of continuum mechanics proposed in \cite{PeshRom2014} that belongs to the class of hyperbolic and thermodynamically compatible (HTC) systems \cite{God1961,Rom1998}. Throughout this paper we will use the Einstein summation convention over repeated indices, and also adopt bold symbols to label vectors and matrices. In three space dimensions, using the indices $1\leq i,k,m \leq 3$, the mathematical model writes
\begin{subequations}\label{eqn.GPR}
	\begin{align}
		& \frac{\partial \rho}{\partial t}+\frac{\partial (\rho v_k)}{\partial
			x_k}=0,\label{eqn.conti}\\[2mm]
		&\frac{\partial \rho v_i}{\partial t}+\frac{\partial \left(\rho v_i v_k + p \, \delta_{ik} +
			\sigma_{ik} + \omega_{ik} \right)}{\partial x_k}=0, \label{eqn.momentum}\\[2mm]
		&\frac{\partial A_{i k}}{\partial t}+\frac{\partial (A_{im} v_m)}{\partial x_k} +
		v_m \left(\frac{\partial A_{ik}}{\partial x_m}-\frac{\partial A_{im}}{\partial x_k}\right)
		=-\dfrac{ \alpha_{ik} }{{\theta}_1(\tau_1)},\label{eqn.deformation}\\[2mm]
		&\frac{\partial J_k}{\partial t}+\frac{\partial \left( J_m v_m + T \right)}{\partial x_k} +
		v_m \left(\frac{\partial J_{k}}{\partial x_m}-\frac{\partial J_{m}}{\partial x_k}\right)  =
		-\dfrac{\beta_k}{{\theta}_2(\tau_2)}, \label{eqn.heatflux}\\[2mm]
		& \frac{\partial \mathcal{E}}{\partial t}+\frac{\partial \left( v_k \mathcal{E} + v_i (p \, \delta_{ik}
			+ \sigma_{ik} + \omega_{ik} ) + q_k \right)}{\partial x_k}=0, \label{eqn.energy}
	\end{align}
\end{subequations}
which also satisfies the entropy inequality
\begin{equation}
	\frac{\partial \rho S}{\partial t}+\frac{\partial \left( \rho S v_k  + \beta_k \right)}{\partial x_k} = \frac{\alpha_{ik} \alpha_{ik}}{T \, {\theta}_1(\tau_1)} + \frac{\beta_k \beta_k}{T \, {\theta}_2(\tau_2)} \geq 0, \label{eqn.entropy}
\end{equation}
with $t \in \mathds{R}^+_0$ being the time and $\x=\{x_k\} \in \mathds{R}^3$ denoting the spatial position vector. The state vector $\mathbf{q}= (\rho, \rho v_i, \mathcal{E}, A_{ik}, J_k)^T$ accounts for the mass density $\rho$, the velocity field $v_i$, the distortion field $A_{ik}$, which for pure elasticity corresponds to the inverse deformation gradient, the specific thermal impulse $J_k$ and the total energy density $\mathcal{E} = \rho E = \mathcal{E}_1 + \mathcal{E}_2 + \mathcal{E}_3 + \mathcal{E}_4$, which is obtained as the sum of internal energy, kinetic energy, energy due to the elastic deformation of the medium and energy due to the specific thermal impulse:
\begin{equation}
	\mathcal{E}_1 = \frac{\rho^\gamma}{\gamma-1} e^{S/c_v}, \quad
	\mathcal{E}_2 = \halb \rho v_i v_i,
	\quad
	\mathcal{E}_3 = \frac{1}{4} \rho c_s^2 \mathring{G}_{ij} \mathring{G}_{ij},
	\quad
	\mathcal{E}_4 = \halb c_h^2 \rho J_i J_i.
\end{equation}
The metric tensor is denoted by $\mathbf{G}$ with $
{G}_{ik} = A_{ji} A_{jk} $ and its trace-free (deviatoric) part $\mathring{\mathbf{G}}$ is given by
$\mathring{G}_{ik} = {G}_{ik} - \frac{1}{3} \, G_{mm} \delta_{ik}
$. The polytropic index $\gamma=c_p/c_v$ is given as the ratio of specific heats at constant pressure and volume, namely $c_p$ and $c_v$, respectively. Furthermore, $T > 0$ is the temperature and $p>0$ is the pressure, which are assumed to be positive.
The elastic and thermal stress tensors read
\begin{equation}
	\sigma_{ik} = A_{ji} \partial_{A_{jk}} \mathcal{E} =
	A_{ji} \alpha_{jk} = \rho c_s^2 G_{ij} \mathring{G}_{jk},
	\qquad
	\omega_{ik} = J_i \partial_{J_{k}} \mathcal{E} = J_i \beta_{k} = \rho c_h^2 J_i J_k,
\end{equation}
and the heat flux is given by
\begin{equation}
	q_k = \partial_{\rho S} \mathcal{E} \, \partial_{J_k} \mathcal{E} = T \beta_k = \rho c_h^2 T J_k.
\end{equation}
They are expressed in terms of the thermodynamic dual variables $\alpha_{ik} = \partial_{A_{ik}} \mathcal{E}$ and $\beta_{k} = \partial_{J_{k}} \mathcal{E}$, the shear sound speed $c_s$, and the parameter $c_h$ related to the heat wave propagation speed $c_T$:
\begin{equation}
	c_T = c_h \sqrt{\frac{T}{c_v}}.
	\label{eqn.cT}
\end{equation}
The source terms in \eqref{eqn.deformation} and \eqref{eqn.heatflux} contain two functions $\theta_1(\tau_1)>0$ and $\theta_2(\tau_2)>0$ that depend on the two relaxation times $\tau_1>0$ and $\tau_2>0$ and are defined as
\begin{equation}\label{eqn:theta1}
	\theta_1 = \frac{1}{3}  \rho \tau_1 \, c_s^2 \, \left| \mathbf{A} \right|^{-\frac{5}{3}},
	\qquad
	\theta_2 = \rho c_h^2 \tau_2.
\end{equation}
The determinants of $\A$ and $\G$ are related to the density of the medium and the reference density $\rho_0$ by the relations
\begin{equation}
	|\A| = \frac{\rho}{\rho_0}, \qquad |\G| = \left( \frac{\rho}{\rho_0} \right)^2.
\end{equation}
The asymptotic limit of the model \eqref{eqn.GPR} has been analyzed in \cite{GPRmodel} at the continuous level and in \cite{LGPR} in the fully discrete setting, showing that for small relaxation times, i.e. when $\tau_1 \to 0$ and $\tau_2 \to 0$, the Navier-Stokes-Fourier limit is obtained. Indeed, the stress tensor $\sigma_{ik}$ and the heat flux $q_k$ tend to
\begin{equation}
	\sigma_{ik} = -\frac{1}{6} \rho_0 c_s^2 \tau_1 \left( \partial_k v_i + \partial_i v_k
	- \frac{2}{3} \left( \partial_m v_m\right) \delta_{ik} \right),
	\label{eqn.asymptoticlimit.stress}
\end{equation}
and
\begin{equation}
	 q_k = -\rho T c_h^2 \tau_2 \partial_k T.
	\label{eqn.asymptoticlimit.heatflux}
\end{equation}
In terms of the relaxation times, \eqref{eqn.asymptoticlimit.stress} and \eqref{eqn.asymptoticlimit.heatflux} identify the dynamic viscosity of the Navier-Stokes-Fourier limit as $\mu = \frac{1}{6} \rho_0 c_s^2 \tau_1$, while the corresponding heat conduction coefficient $\lambda = \rho c_h^2 \tau_2 T$ is not a constant, since it depends on the temperature.

Let us also introduce the governing equation for the metric tensor, that can be derived from the definition $\G=\A^\top \A$ and the evolution equation \eqref{eqn.deformation}:
\begin{equation}
\frac{\partial G_{i k}}{\partial t}
+ v_m \frac{\partial G_{ik} }{\partial x_m}
+ G_{im} \frac{\partial v_m }{\partial x_k}
+ G_{mk} \frac{\partial v_m }{\partial x_i}
=-\dfrac{ 2 \sigma_{ik} }{\theta_1(\tau_1)}.
\label{eqn.pde.metric}
\end{equation}

The time evolution for the trace-free part $\mathring{\G} = \G - \frac{1}{3} \text{tr}(\G) \mathbf{I}$ is instead given by 
\begin{equation}
	\frac{\partial \mathring{G}_{i k}}{\partial t}
	+ v_m \frac{\partial \mathring{G}_{ik} }{\partial x_m}
	+ G_{im} \frac{\partial v_m }{\partial x_k}
	+ G_{mk} \frac{\partial v_m }{\partial x_i}
	- \frac{2}{3} \delta_{ik} G_{nm} \frac{\partial v_m }{\partial x_n}
	=-\dfrac{ 2 \mathring{\sigma}_{ik} }{\theta_1(\tau_1)}.
	\label{eqn.pde.devmetric}
\end{equation}
and will become an essential part later in this paper.

Finally, the equation for the specific thermal impulse \eqref{eqn.heatflux} can also be rewritten as
\begin{equation}
	\frac{\partial J_{k}}{\partial t}
	+ v_m \frac{\partial J_{k} }{\partial x_m}
	+ J_{m} \frac{\partial v_m }{\partial x_k}
	+ \frac{\partial T }{\partial x_k}
	=-\dfrac{ \beta_{k} }{{\theta}_2(\tau_2)}.
	\label{eqn.pde.J}
\end{equation}

In addition to the material wave induced by convection, three other time scales characterize the model \eqref{eqn.GPR}, hence yielding four subsystems discussed in detail in the following subsections.

\subsection{Convective subsystem}
The convective subsystem only accounts for advection phenomena based on the material speed $\| \v \|$ of the medium, and it reads
\begin{subequations}\label{eqn.GPR.s1}
	\begin{align}
		& \frac{\partial \rho}{\partial t}+\frac{\partial (\rho v_k)}{\partial
			x_k}=0,\label{eqn.conti.1}\\[2mm]
		&\frac{\partial \rho v_i}{\partial t}+\frac{\partial \left(\rho v_i v_k  \right)}{\partial x_k}=0, \label{eqn.momentum.1}\\[2mm]
		&\frac{\partial A_{i k}}{\partial t}  +
		v_m  \frac{\partial A_{ik}}{\partial x_m}
		=0,\label{eqn.deformation.1}\\[2mm]
		&\frac{\partial J_k}{\partial t}
		 + v_m \frac{\partial J_{k}}{\partial x_m} = 0, \label{eqn.heatflux.1}\\[2mm]
		& \frac{\partial \mathcal{E}}{\partial t}+\frac{\partial \left( v_m \mathcal{E}_{2,3,4}  \right)}{\partial x_m}=0, \label{eqn.energy.1}
	\end{align}
\end{subequations}
with $\mathcal{E}_{2,3,4} = \mathcal{E}_{2}+\mathcal{E}_{3}+\mathcal{E}_{4}$.
The convective subsystem in direction $x_1$ has the eigenvalues $\lambda_1=0$ and $\lambda_i = v_1$, $i \in \left\{2, 3, \cdots, 17 \right\}$. The convective subsystem is only weakly hyperbolic, since a full set of linearly independent eigenvectors cannot be established, exactly as in the well-known Toro-V\'azquez splitting, see \cite{ToroVazquez}. 

\subsection{Temperature subsystem}
By introducing the heat Mach number
\begin{equation}
	M_h = \frac{\| \v \|}{c_T},
\end{equation}
with respect to the heat wave propagation speed $c_T$ given by \eqref{eqn.cT}, we can identify the temperature subsystem which couples the temperature gradient in \eqref{eqn.heatflux} and the heat flux of the total energy equation \eqref{eqn.energy}:
\begin{subequations}\label{eqn.GPR.s2}
	\begin{align}
		&\frac{\partial \rho}{\partial t} = 0, \\[2mm]
		&\frac{\partial J_k}{\partial t}
		+  \frac{\partial T}{\partial x_k} = -\frac{1}{\theta_2} \beta_k, \label{eqn.heatflux.2} \\[2mm]
		& \frac{\partial \mathcal{E}}{\partial t}+\frac{\partial q_k }{\partial x_k}=0. \label{eqn.energy.2}
	\end{align}
\end{subequations}
The eigenvalues of the temperature subsystem in the $x_1$ direction are
$\lambda_{1,5} = \mp c_h \sqrt{T/c_v} = \mp c_T$ and $\lambda_{2,3,4}=0$, and the temperature subsystem is hyperbolic for $\rho >0$ and $T>0$.
In terms of the primitive variables $J_k$ and $T$ and setting $\theta_2 = \tau_2 \rho c_h^2$ it takes the simple form
\begin{subequations}\label{eqn.GPR.s2b}
	\begin{align}
		&\frac{\partial \rho}{\partial t} = 0, \\[2mm]
		&\frac{\partial J_k}{\partial t}
		+  \frac{\partial T}{\partial x_k} = -\frac{1}{\tau_2} J_k, \label{eqn.heatflux.2b} \\[2mm]
		& \frac{c_v}{T c_h^2} \frac{\partial T}{\partial t}+\frac{\partial J_k }{\partial x_k}=0. \label{eqn.energy.1b}
	\end{align}
\end{subequations}

\subsection{$\G$-$\J$-$\v$ subsystem}
The third subsystem governs the mechanical response of the medium and is therefore called the mechanical subsystem. The momentum equation \eqref{eqn.momentum} is coupled via the elastic and the thermal stress tensors $\sigma_{ik}$ and $\omega_{ik}$ with the distortion field $\A$ and the specific thermal impulse $\J$, where $\sigma_{ik}$ becomes a viscous stress only in the stiff relaxation limit $\tau_1 \to 0$; although $\omega_{ik}$ is of thermal origin, it enters the momentum balance and therefore contributes to the mechanical response in the same way as $\sigma_{ik}$ does. Therefore, this subsystem is related not only to the heat Mach number $M_h$ introduced above, but also to the shear Mach number
\begin{equation}
	M_s = \frac{\| \v \|}{c_s},
\end{equation}
defined with respect to the shear sound speed. The $\G$-$\J$-$\v$ subsystem reads
\begin{subequations}\label{eqn.GPR.s3}
	\begin{align}
	&\frac{\partial \rho v_i}{\partial t}
	+\frac{\partial \left(\sigma_{ik} + \omega_{ik} \right)}{\partial x_k}=0, \label{eqn.momentum.3}\\[2mm]
	& \frac{\partial J_{k}}{\partial t}
	+ J_{m} \frac{\partial v_m }{\partial x_k}
	=-\dfrac{ J_{k} }{\tau_2},
	\label{eqn.heatflux.3} \\[2mm]
	&	\frac{\partial G_{i k}}{\partial t}
	+ G_{im} \frac{\partial v_m }{\partial x_k}
	+ G_{mk} \frac{\partial v_m }{\partial x_i}
	=-\dfrac{ 2 \sigma_{ik} }{\theta_1(\tau_1)}, \label{eqn.metric.3}  \\[2mm]
	& \frac{\partial \mathcal{E}}{\partial t}+\frac{\partial v_i (\sigma_{ik}+\omega_{ik}) }{\partial x_k}=0. \label{eqn.energy.3}
	\end{align}
\end{subequations}
The eigenvalues of \eqref{eqn.GPR.s3} cannot be determined analytically, but for $\G = \text{diag}(G_{kk})$, $J_3 = 0$, the eigenvalues and a full set of linearly independent eigenvectors can be calculated with a computer algebra system, although the resulting expressions are too lengthy to be reported here.

\subsection{Pressure subsystem}
The last subsystem is the classical pressure subsystem related to the acoustic Mach number
\begin{equation}
	M_a = \frac{\| \v \|}{c_0}, \qquad c_0 = \sqrt{\frac{\gamma p}{\rho}},
\end{equation}
with $c_0$ being the adiabatic sound speed. Originally proposed in \cite{MunzPark,ToroVazquez}, this subsystem couples the momentum with the pressure gradient in \eqref{eqn.momentum} and the specific enthalpy in the energy equation \eqref{eqn.energy}. It is given by
\begin{subequations}\label{eqn.GPR.s4}
	\begin{align}
		&\frac{\partial \rho}{\partial t} = 0, \\[2mm]
		&\frac{\partial \rho v_i}{\partial t}+\frac{\partial  p    }{\partial x_i}=0, \label{eqn.momentum.4}\\[2mm]
		& \frac{\partial \mathcal{E} }{\partial t}+\frac{\partial \left( v_k ( \mathcal{E}_1 +  p ) \right)}{\partial x_k}=0. \label{eqn.energy.4}
	\end{align}
\end{subequations}
The eigenvalues in the $x_1$ direction are
$\lambda_{1,5}=\halb \left( v_1 \mp a \right) $  and $\lambda_{2,3,4}=0$, with
$a^2 = v_1^2 + 4 c_0^2$. For $\rho >0$ and $p>0$, the pressure subsystem is hyperbolic since a full set of linearly independent eigenvectors does exist \cite{MunzPark,ToroVazquez}. The proposed splitting is additive in the left hand side of \eqref{eqn.GPR}, where the equation for the metric tensor takes the place of the one for the distortion field through $G_{ik} = A_{ji} A_{jk}$. The relaxation source of the specific thermal impulse, by contrast, deliberately appears both in the temperature subsystem and in the $\G$-$\J$-$\v$ subsystem, since each of them has to be individually consistent with its own stiff relaxation limit.

% % % % % % % % % % % % % % % % % % % % % % % % % % % % % %
%        Numerical method.
% % % % % % % % % % % % % % % % % % % % % % % % % % % % % %

\section{Numerical method}\label{sec.method}
We consider a two-dimensional domain $\Omega$ with boundary $\partial \Omega$ which is discretized by a set of non-overlapping triangular control volumes $\omega_c$ of volume $|\omega_c|$. In what follows, the indices $a$, $c$, $d$ refer to cells, while the indices $p$, $q$, $r$ refer to vertices. Grid indices are written as superscripts on the discrete fields, while their subscripts are reserved for tensor indices. The set of vertices belonging to cell $\omega_c$ is referred to with $\mathcal{P}_c$, while the set of cells sharing node $p$ is indicated with $\mathcal{C}_p$.
\begin{figure}[!htbp]
	\begin{center}
		\begin{tabular}{cc}
			\includegraphics[trim=0 0 0 0,clip,width=0.47\textwidth]{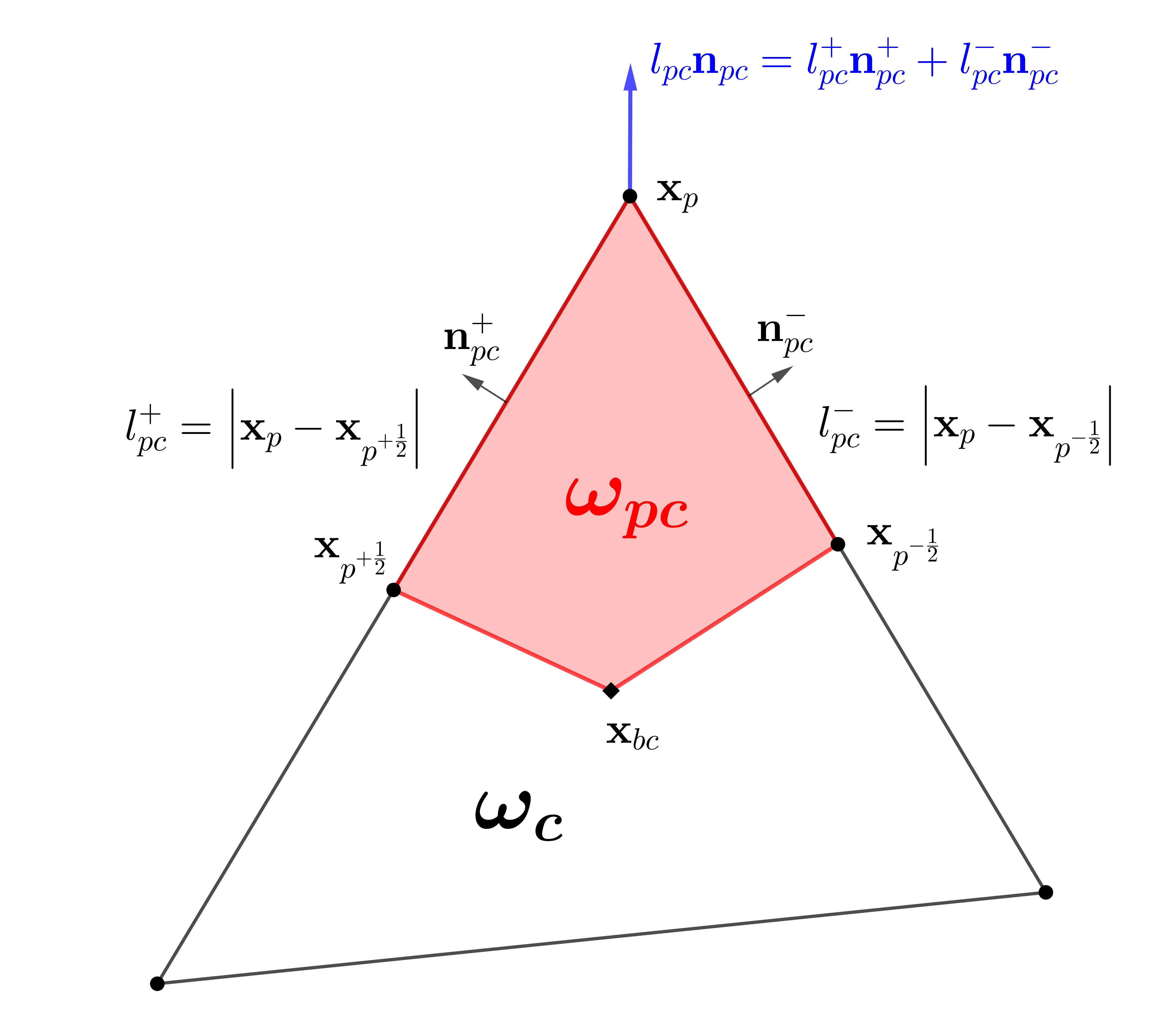} &
			\includegraphics[trim=100 100 100 100,clip,width=0.47\textwidth]{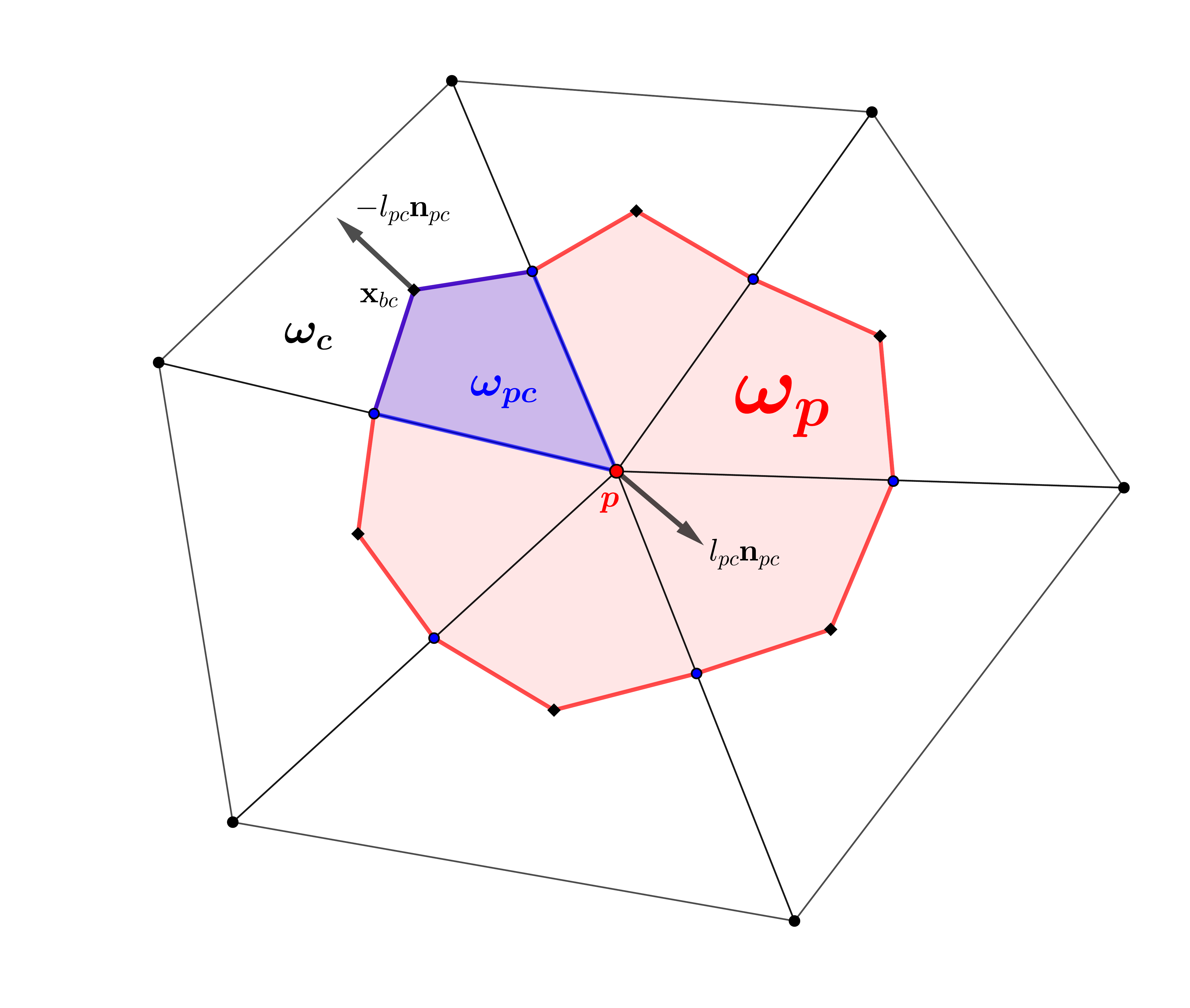} \\
		\end{tabular}
		\caption{Notation for the primary cell $\omega_c$ (left) and for the dual cell $\omega_p$ (right).}
		\label{fig.Grid}
	\end{center}
\end{figure}
For each cell $\omega_c$, depicted in the left panel of Figure \ref{fig.Grid}, we define its barycenter of coordinates $\x^{bc}=1/(d+1) \, \sum_{p \in \mathcal{P}_c}\x^p$. With respect to cell $\omega_c$, one can also introduce the left and right midpoints of the edges impinging on node $p$, namely $p^{-\halb}$ and $p^{+\halb}$, respectively. The half lengths of these edges are given by $l^{pc,-}=\left| \x^p - \x^{p-\halb} \right|$ and $l^{pc,+}=\left| \x^p - \x^{p+\halb} \right|$, and the corner normal is computed as
\begin{equation}
	\lnpc = l^{pc,+} \n^{pc,+} + l^{pc,-} \n^{pc,-}.
	\label{eqn.lpcnpc}
\end{equation}
The corner normal $\n^{pc}$ is oriented from the cell center $\x^{bc}$ towards the vertex $\x^p$, hence outward with respect to the primal cell $\omega_c$, following the convention of the Lagrangian finite volume literature \cite{Maire2007,Despres2009,Maire2011,Maire2020}, while the opposite corner normal $\n^{cp} = -\n^{pc}$ points from the vertex towards the cell center and is therefore inward with respect to $\omega_c$.
By construction the corner vectors satisfy the discrete Gauss theorem
\begin{equation}
	\sum \limits_{p \in \mathcal{P}_c} \lnpc = 0.
	\label{eqn.gauss}
\end{equation}
Let us now define the sub-cell $\omega_{pc}$ by connecting the vertex with position $\x^p$, the cell barycenter $\x^{bc}$, and the left and right midpoints of the edges impinging on node $p$, namely the points of coordinate $\x^{p-\halb}$ and $\x^{p+\halb}$. The union of all sub-cells around a node $p$ yields the dual cell $\omega_p$ given in the right panel of Figure \ref{fig.Grid}.

\begin{figure}[!htbp]
	\begin{center}
		\begin{tabular}{c}
			\includegraphics[width=0.65\textwidth]{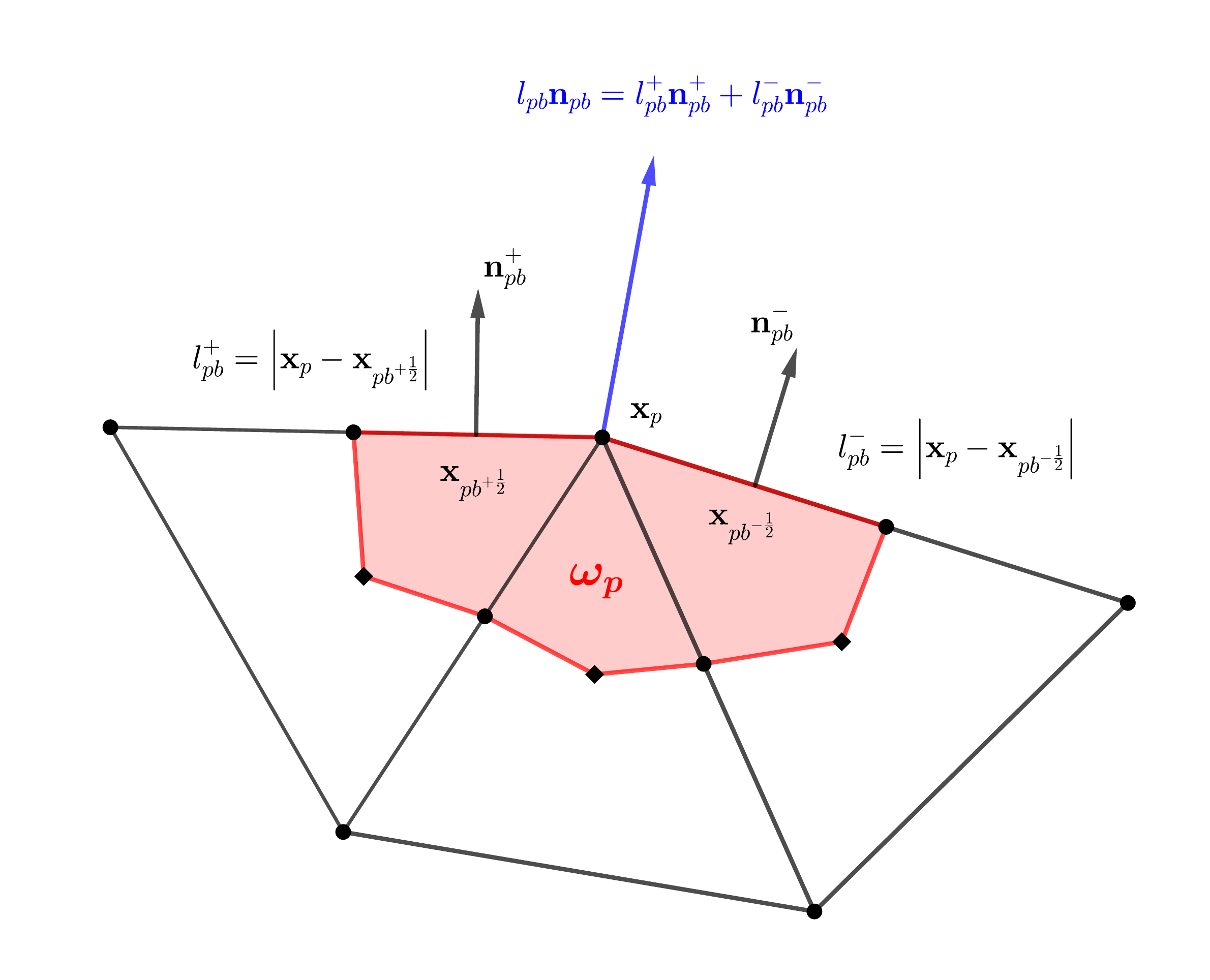}
		\end{tabular}
		\caption{Dual cell $\omega_p$ located on the boundary $\partial \Omega$ of the computational domain.}
		\label{fig.bnd}
	\end{center}
\end{figure}

To provide a complete geometry description, let us also address the domain boundary $\partial \Omega$, thus considering a dual cell $\omega_p$ lying on $\partial \Omega$ as depicted in Figure \ref{fig.bnd}.
The corner vector on the boundary is given by
\begin{equation}
	\lnpb = l^{pb,+} \n^{pb,+} + l^{pb,-} \n^{pb,-},
	\label{eqn.lpcnpcb}
\end{equation}
where the outward normals $l^{pb,+} \n^{pb,+}$ and $l^{pb,-} \n^{pb,-}$ are referred to the boundary sides impinging on node $p$. These normals are essential to satisfy the discrete Gauss theorem also for a boundary node, and they are indeed linked to the corner vectors with the following relation:
\begin{equation}
	\sum \limits_{c \in \mathcal{C}_p} \lnpc = \lnpb.
	\label{eqn.gaussb}
\end{equation}

We adopt a vertex-staggered discretization, thus considering scalar fields $\phi^p=\phi(\x^p)$ and vector fields $\mathbf{a}^p=\mathbf{a}(\x^p)$ defined at the vertices of the primal grid, and scalar and vector fields defined at the centers of the primal grid $\phi^c=\phi(\x^c)$ and $\mathbf{a}^c=\mathbf{a}(\x^c)$, respectively. We recall that geometry-related indices such as cell and vertex number are reported as superscripts, as well as the time index, while vector and tensor indices of field variables are written as subscripts. Based on this notation, the discrete versions of gradient ($\mathbf{u} = \nabla \phi$), divergence ($\psi = \nabla \cdot \mathbf{a}$) and curl ($\mathbf{b} = \nabla \times \mathbf{a}$) operators and their duals are naturally defined as follows (see also \cite{Maire2007,Despres2009,Maire2011,Maire2020,HTCLagrange,HTCLagrangeGPR}):
\begin{alignat}{2}
	& u^c_k=\partial^{cp}_k \phi^p   = \frac{1}{|\omega_c|} \sum_{p \in \mathcal{P}_c}  l^{pc} n^{pc}_k \,  \phi^p,  \qquad
	&&u^p_k=\partial^{pc}_k \phi^c   = \frac{1}{|\omega_p|} \sum_{c \in \mathcal{C}_p}  l^{pc} n^{cp}_k \,   \phi^c,
	\\
	& \psi^c = \partial^{cp}_k a^p_k   =  \frac{1}{|\omega_c|} \sum_{k} \sum_{p \in \mathcal{P}_c}  l^{pc} n^{pc}_k \,  a^p_k,  \qquad
	&&\psi^p = \partial^{pc}_k a^c_k   =  \frac{1}{|\omega_p|} \sum_{k} \sum_{c \in \mathcal{C}_p}  l^{pc} n^{cp}_k \,  a^c_k,
	\\
	&b^c_i= \epsilon_{ijk} \partial^{cp}_j a^p_k   = \frac{\epsilon_{ijk}}{|\omega_c|} \sum_{j} \sum_{k} \sum_{p \in \mathcal{P}_c}  l^{pc} n^{pc}_j   a^p_k,  \qquad
	&&b^p_i= \epsilon_{ijk} \partial^{pc}_j a^c_k   = \frac{\epsilon_{ijk}}{|\omega_p|} \sum_{j} \sum_{k} \sum_{c \in \mathcal{C}_p}  l^{pc} n^{cp}_j   a^c_k.
\end{alignat}
%On the regular Cartesian meshes considered in this paper we have
%$ |\Omega_c| = |\Omega_p| = \Delta x \Delta y$ and the corner normal vectors $n_k^{cp}$ are simply
%$ n_k^{pc} = \halb ( \pm \Delta y, \pm \Delta x)^T$, with $n_k^{cp} = -n_k^{pc}$. For the definition of the corner normals on general meshes see \cite{Maire2007,Despres2009,Maire2011,Maire2020,HTCLagrange,HTCLagrangeGPR}.
%
The use of the two opposite corner normals makes the two families of operators adjoint to each other with respect to the volume-weighted scalar products on the primal and on the dual grid, since
\begin{equation}
	\sum \limits_{p} |\omega_p| \, \phi^p \, \partial^{pc}_k a^c_k = - \sum \limits_{c} |\omega_c| \, a^c_k \, \partial^{cp}_k \phi^p
	\label{eqn.adjoint}
\end{equation}
follows directly from the definitions, both sides being one and the same sum over all pairs $(p,c)$ of a vertex and one of its adjacent cells. This identity is the discrete counterpart of the integration by parts $\int \phi \, \nabla \cdot \mathbf{a} \, d\mathbf{x} = - \int \mathbf{a} \cdot \nabla \phi \, d\mathbf{x}$, see also \cite{CompatibleDG1}.
In \cite{HTCLagrangeGPR,spfvEulerHeat} it is shown that the discrete operators introduced above satisfy the discrete vector calculus identities
\begin{equation}
	\epsilon_{ijk} \partial_j^{pc} \partial_k^{cq} \phi^q = 0, \qquad
	\epsilon_{ijk} \partial_i^{pc} \partial_j^{cq} a_k^q = 0,
\end{equation}
which are discrete analogues of $\nabla \times \nabla \phi = 0$ and $\nabla \cdot \nabla \times \mathbf{a}= 0$. This property is well-established by now and is therefore not proven here again. 

\medskip

We set the following vertex-staggered discretization:
\begin{itemize}
	\item primal triangular cells: momentum density $(\rho \v)^c$, distortion field $\A^c$ and specific thermal impulse $\J^c$;
	\item dual polygonal cells: mass density $\rho^p$, total energy density $\mathcal{E}^p$, pressure $p^p$ and temperature $T^p$.
\end{itemize}
To transfer quantities between the two grids, we use volume-weighted averages. For instance, for obtaining a cell-centered pressure and a vertex-centered velocity vector we compute
\begin{equation}
	p^c = \frac{1}{|\omega_c|} \sum_{p \in \mathcal{P}_c} |\omega_{pc}| \, p^p, \qquad  \v^p = \frac{1}{|\omega_p|} \sum_{c \in \mathcal{C}_p} |\omega_{pc}| \, \v^c,
	\label{eqn.averages}
\end{equation}
where $|\omega_{pc}|$ is the volume of the sub-cell $\omega_{pc}$ as shown in Figure \ref{fig.Grid}.
The velocity is used on both grids: in the cells it follows from the cell-centered momentum and the cell-averaged density obtained by an average analogous to \eqref{eqn.averages}, while at the vertices it is either obtained by the volume-weighted average described above, or computed directly as the unknown of the vector wave equation \eqref{eqn.velocity.wave}.
Throughout this section the superscripts $n$, $*$, $**$, $***$ and $n+1$ count the sequential steps of the algorithm that a quantity has already undergone, so that a variable which is not modified by a given step simply keeps its superscript and the number of stars differs from one variable to another; the density, for example, is already final after the convective step, whereas the specific thermal impulse is updated in the convective, in the temperature and in the mechanical step before it receives its compatible update at time $t^{n+1}$.

%In our scheme the density, the pressure, the temperature, the distortion field and the specific thermal impulse are defined in the cells $c$, while the momentum is defined in the vertices $p$. When needed, simple arithmetic averaging between the grid location is carried out to transfer data from one grid to the other one.
%
%
%In the following we present the fully discrete semi-implicit four-split scheme using the discrete operators given above.
%Therein, the steps of the numerical scheme are given by the order of the subsystems given above, i.e. we first solve the convective, then the temperature subsystem, followed by the G-J-v and pressure sub-systems.
%
\subsection{Discrete convective subsystem}
In compact notation \eqref{eqn.GPR.s1} can be written as
\begin{equation}
	\frac{\partial \q}{\partial t} + \frac{\partial \f_k}{\partial x_k} + B_k \frac{\partial \q}{\partial x_k} = 0,
\end{equation}
where $\f$ denotes the flux tensor and $B$ is the matrix containing the non-conservative contributions.
An explicit finite volume scheme is used to evolve the convective subsystem. Note that, besides the density and the total energy, also the pressure and the temperature are defined at the nodes $p$, while the momentum, the distortion field and the thermal impulse are defined in the cell centers $c$. In order to be compatible with the low Mach number limit of the equations we use the following special update formulas for the vertex-based mass density and total energy density: 
\begin{eqnarray}
	\rho^{*,p}        & =&  \rho^{n,p} - \Delta t \, \partial_k^{pc}  \, f^{\rho,c}_k, \label{eqn.fv.mass} \\ 
	\mathcal{E}^{*,p} & =&  \mathcal{E}^{n,p} - \Delta t \, \partial_k^{pc}  \, f^{\mathcal{E},c}_k, 
	\label{eqn.fv.rho}
\end{eqnarray}
with the multidimensional numerical flux tensors chosen of the Rusanov-type and which read 
\begin{eqnarray}
	 f^{\rho,c}_k   & = & v_k^{n,c}  \rho^{n,c} - \halb \ell^c s^c_{\max} \partial_k^{cp} \rho^{n,p}, \label{eqn.rusanov.mass} \\
	 f^{\mathcal{E},c}_k & = & v_k^{n,c} \, \mathcal{E}_{2,3,4}^{n,c} - \halb \ell^c s^c_{\max}  \partial_k^{cp} \mathcal{E}^{n,p},
	\label{eqn.rusanov.md}
\end{eqnarray}
where $\ell^c = 4 |\omega_c| / |\partial \omega_c|$ and $\ell^p = 4 |\omega_p| / |\partial \omega_p|$ are characteristic
length scales in the primal cell and in the dual cell, respectively, while
$s^c_{\max} = \max \limits_{p \in \mathcal{P}_c}(|\v^p|)$ and
$s^p_{\max} = \max \limits_{c \in \mathcal{C}_p}(|\v^c|)$ are the maximum convective speeds
in the primal cell and in the dual cell, respectively. Furthermore we set 
$\rho^{n+1,p}= \rho^{*,p}$ since the density is not modified in the remaining sub-systems.  

The nonlinear advection of the cell-centered momentum $m_k^{n,c} = \rho^{n,c} v_k^{n,c}$ is carried out with the aid of a classical unstructured finite volume scheme with edge-based numerical flux 
\begin{equation}
   m_i^{*,c} = m_i^{n,c} - \frac{\Delta t}{|\omega_c|} \sum \limits_{a \in \mathcal{N}_c} | \partial \omega_{ac} | f^{m,ac}_{i}  
	\label{eqn.mom.transport}
\end{equation}
with $a$ being the index of the neighbor cell, $| \partial \omega_{ac} |$ the length of the common edge between cells $\omega_a$ and $\omega_c$ and $\mathcal{N}_c$ the set of neighbors that share a common edge with cell $\omega_c$. The numerical flux in normal direction $n^{ac}_k$ to the edge is given by 
\begin{equation}
	f^{m,ac}_{i} = \halb \left( v_j^{n,a} + v_j^{n,c} \right) n^{ac}_j \, \halb \left( m_i^{n,a} + m_i^{n,c} \right) - \halb s_{\max}^{ac} \left( m_i^{n,a} - m_i^{n,c} \right). 
\end{equation}
where $\n^{ac}$ is the unit normal to the common edge $\partial \omega_{ac}$ pointing from the cell $\omega_c$ to its neighbour $\omega_a$, and where
\begin{equation}
	s_{\max}^{ac} = \max \left( \left| v_k^{n,a} n_k^{ac} \right|, \left| v_k^{n,c} n_k^{ac} \right| \right)
	\label{eqn.smax.ac}
\end{equation}
is the maximum wave speed on that edge.
Last but not least the transport step of the distortion field $\A$ and of the specific thermal impulse $\J$ reads 
\begin{equation}
	A_{ik}^{*,c} = A_{ik}^{n,c} - \frac{\Delta t}{|\omega_c|} \sum \limits_{a \in \mathcal{N}_c} | \partial \omega_{ac} | D^{A,ac}_{ik},
	\label{eqn.A.transport}  
\end{equation}
and 
\begin{equation}
	J_{k}^{*,c} = J_{k}^{n,c} - \frac{\Delta t}{|\omega_c|} \sum \limits_{a \in \mathcal{N}_c} | \partial \omega_{ac} | D^{J,ac}_{k},
	\label{eqn.J.transport}  
\end{equation}
with the numerical fluctuations in normal direction $n^{ac}_k$ to the edge given by 
\begin{equation}
	D^{A,ac}_{ik} = \halb \left( v_j^{n,a} + v_j^{n,c} \right) n^{ac}_j \, \halb \left( A_{ik}^{n,a} - A_{ik}^{n,c} \right) - \halb s_{\max}^{ac} \left( A_{ik}^{n,a} - A_{ik}^{n,c} \right), 
\end{equation}
and 
\begin{equation}
	D^{J,ac}_{k} = \halb \left( v_j^{n,a} + v_j^{n,c} \right) n^{ac}_j \, \halb \left( J_{k}^{n,a} - J_{k}^{n,c} \right) - \halb s_{\max}^{ac} \left( J_{k}^{n,a} - J_{k}^{n,c} \right).  
\end{equation}

This subsystem is the only one to be solved explicitly, thus implying a CFL-type time step restriction which merely depends on the material speed, that is
\begin{equation}
    \Delta t = \text{CFL} \min \limits_{c} \left( \frac{\ell^c}{s_{\max}^{c}}
    \right),
	\label{eqn.cfl}
\end{equation}
for a CFL coefficient $\text{CFL} \leq \halb$.
\subsection{Discrete temperature subsystem}
After the nonlinear convective terms, which are the only ones to be discretized explicitly, we solve the auxiliary temperature subsystem \eqref{eqn.GPR.s2b}, which couples $T$ and $J_k$ and which contains the temperature gradient of \eqref{eqn.pde.J}, implicitly. This avoids a CFL-type time step restriction based on the propagation speed of the heat waves.
Applying our compatible discrete derivative operators, the discrete temperature subsystem reads
\begin{subequations}
	\begin{align}
    J_k^{**,c} &= J_k^{*,c} - \Delta t \, \partial^{cp}_k T^{**,p} - \frac{\Delta t}{\tau_2} J_k^{**,c}, \label{eqn.heatflux.disc.2} \\
    \frac{c_v}{T^{n,p} c_h^2} T^{**,p} &= \frac{c_v}{T^{n,p} c_h^2}  T^{*,p} - \Delta t \,
    \partial^{pc}_k \, J_k^{**,c}. \label{eqn.energy.disc.2}
    \end{align}
\end{subequations}
Inserting \eqref{eqn.heatflux.disc.2} in equation \eqref{eqn.energy.disc.2} leads to a discrete scalar wave equation for the unknown temperature $T^{\ast\ast,p}$
\begin{equation}
	  \left( 1+\frac{\Delta t}{\tau_2} \right) \frac{c_v}{T^{n,p} c_h^2} \, T^{**,p}  - \Delta t^2
	\partial^{pc}_k \left( \partial^{cq}_k T^{**,q} \right)
	= \left( 1+\frac{\Delta t}{\tau_2} \right) \frac{c_v}{T^{n,p} c_h^2} \, T^{*,p} - \Delta t \,
		\partial^{pc}_k \left( J_k^{*,c} \right). \label{eqn.wave.T}
\end{equation}
The resulting algebraic system is symmetric and positive definite. We solve it via the classical matrix-free conjugate gradient method.
Once $T^{**,p}$ is known, the discrete thermal impulse and total energy are updated as follows:
\begin{subequations}
	\begin{align}\label{eqn.heat.disc.final}
		J^{**,c}_k &= \left(\frac{\tau_2}{\tau_2 + \Delta t }\right)(J^{*,c}_k - \Delta t \partial_k^{cp} T^{**,p}), \\
		\mathcal{E}^{**,p} &= \mathcal{E}^{*,p} - \Delta t \, \partial_k^{pc} q_k^{**,c}, \label{eqn.heat.disc.final.E}
	\end{align}
\end{subequations}
where $q_k^{**,c} =  \rho^{*,c} \, c_h^2 \, T^{**,c} \, J_k^{**,c}$ is the heat flux, with $T^{**,c}$ and $\rho^{*,c}$ obtained by averaging from the nodes to the cells.
\subsection{Discrete $\mathring{\mathbf{G}}$-$\J$-$\v$ subsystem}
We now discretize the $\mathring{\mathbf{G}}$-$\J$-$\v$ subsystem, see also \cite{FourSplit}. Rather than discretizing the evolution equation \eqref{eqn.metric.3} of the metric tensor directly, we evolve its trace-free part, so that \eqref{eqn.devG.1} below is the discrete counterpart of \eqref{eqn.pde.devmetric}; it is obtained by applying the deviator operator to the discrete form of \eqref{eqn.metric.3}, where the symmetry of $\G$ can be used in order to simplify the trace contribution; the intermediate steps of this computation can be found in \cite{FourSplit}. Among the different discrete formulations of the mechanical subsystem that have been investigated, the one based on $\mathring{\mathbf{G}}$ is the only one for which the asymptotic-preserving property with respect to the Navier-Stokes stress tensor could be proven, see \cite{FourSplit}, so that the discrete equation for $\mathring{\mathbf{G}}$ is an essential ingredient of the method proposed here; the corresponding analysis is carried out in Section \ref{sec.asymptotic}. The aim of the manipulations that follow is to express $\devG_{ik}^{**,c}$ and $J_k^{***,c}$ in terms of the nodal velocity, so that a single vector wave equation for $v_i^{**,p}$ remains to be solved. The discrete subsystem reads:
\begin{subequations}
	\begin{align}
		\rho^{*,p} v_i^{**,p}  &= \rho^{*,p} v_i^{*,p} - \Delta t \, \partial^{pc}_k
		\left(  \sigma_{ik}^{**,c} + \omega_{ik}^{***,c} \right), \label{eqn.mom.disc.3} \\
			\devG_{ik}^{**,c}
		&= \devG_{ik}^{*,c}
		- \Delta t\left(G_{im}^{n,c}\partial_k^{cp} v_m^{**,p}+G_{mk}^{n,c}\partial_i^{cp} v_m^{**,p}-\frac{2}{3}\delta_{ik}G_{nm}^{n,c}\partial_n^{cp}v_m^{**,p}\right)
		-\frac{2\Delta t}{\theta_1^{n,c}}\mathring{\sigma}_{ik}^{**,c},
		 \label{eqn.devG.1} \\
		J_{k}^{***,c} &= J_{k}^{**,c}
		- \Delta t \, J_{m}^{**,c} \, \partial_k^{cp} v_m^{**,p}
		-  \dfrac{  \Delta t  }{ \tau_2} \, J_{k}^{***,c},
		\label{eqn.J.disc.3}
	\end{align}
\end{subequations}
where $v_i^{**,p}$, $\devG_{ik}^{**,c}$ and $J_k^{***,c}$ are the unknowns, and where the evolution equation of the thermal impulse discretizes the remaining, velocity gradient contribution of \eqref{eqn.pde.J}. Note that $G_{ik}^{*,c} = 
A_{ji}^{*,c} A_{jk}^{*,c}$, with $A_{jk}^{*,c}$ obtained from the explicit transport stage. 
Since our goal is to obtain a linearly implicit method, we introduce the following semi-implicit approximations for the two stress tensors:
\begin{equation}
	\sigma_{ik}^{**,c} = \rho^{n,c} c_s^2 G_{il}^{n,c}\devG_{lk}^{**,c},
	\qquad
	\omega_{ik}^{***,c} = \rho^{n,c} c_h^2 J_i^{**,c} J_k^{***,c}.
	\label{eqn.si.sigma.omega}
\end{equation}
The semi-implicit discretization of $\mathring{\sigma}$, required in \eqref{eqn.devG.1}, reads
\begin{equation}
	\mathring{\sigma}_{ik}^{**,c} = \rho^{n,c} \, c_s^2 \left(G_{il}^{n,c}\devG_{lk}^{**,c}-\frac{1}{3}G_{mn}^{n,c}\devG_{nm}^{n,c}\delta_{ik}\right).
	\label{eqn.si.devsigma}
\end{equation}
Substituting \eqref{eqn.si.devsigma} into \eqref{eqn.devG.1} and rearranging terms leads to
\begin{eqnarray}
	\left(\delta_{il}+\frac{2\Delta t \rho^{n,c} c_s^2}{\theta_1^{n,c}}G_{il}^{n,c}\right)\devG_{lk}^{**,c} &=& \devG_{ik}^{*,c}
	- \Delta t\left(G_{im}^{n,c}\partial_k^{cp} v_m^{**,p}+G_{mk}^{n,c}\partial_i^{cp} v_m^{**,p}-\frac{2}{3}\delta_{ik}G_{nm}^{n,c}\partial_n^{cp}v_m^{**,p}\right)  \nonumber \\
	&&+\frac{2\Delta t \rho^{n,c} c_s^2}{3 \theta_1^{n,c}} G_{mn}^{n,c}\devG_{nm}^{n,c}\delta_{ik}.
	\label{eqn.devG.3}
\end{eqnarray}
We now define
\begin{equation}
	\Sigma_{il}^{n,c}:=\delta_{il}+\frac{2\Delta t \rho^{n,c} c_s^2}{\theta_1^{n,c}}G_{il}^{n,c},
\end{equation}
which is a symmetric positive definite $3\times 3$ matrix that can be directly inverted, see \cite{FourSplit}.
Multiplication of \eqref{eqn.devG.3} by $\ISigma{li}$ leads to an expression for $\mathring{\mathbf{G}}^{**,c}$ in terms of the sole remaining unknown $v_m^{**,p}$:
\begin{align}
	\devG_{lk}^{**,c} &= \ISigma{li}\devG_{ik}^{*,c} - \Delta t  \ISigma{li}\left(G_{im}^{n,c}\partial_k^{cp} v_m^{**,p}+G_{mk}^{n,c}\partial_i^{cp} v_m^{**,p}-\frac{2}{3}\delta_{ik}G_{nm}^{n,c}\partial_n^{cp}v_m^{**,p}\right) \nonumber \\*
	&\quad +\frac{2\Delta t \rho^{n,c} c_s^2}{3 \theta_1^{n,c}}  \ISigma{lk}G_{mn}^{n,c}\devG_{nm}^{n,c}.
	\label{eqn.devG.4}
\end{align}
For $J_k$ one obtains
\begin{equation}
	J_k^{***,c}=\frac{1}{1+\frac{\Delta t}{\tau_2}}J_k^{**,c}-\frac{\Delta t}{1+\frac{\Delta t}{\tau_2}} J_m^{**,c}\partial_k^{cp}v_m^{**,p}.
	\label{eqn.J.disc.5}
\end{equation}
Substituting \eqref{eqn.J.disc.5} into \eqref{eqn.si.sigma.omega} yields
\begin{equation}
	\omega_{ik}^{***,c}=\frac{\rho^{n,c} c_h^2}{1+\frac{\Delta t}{\tau_2}}J_i^{**,c}J_k^{**,c}-\frac{\Delta t \rho^{n,c} c_h^2}{1+\frac{\Delta t}{\tau_2}}J_i^{**,c}J_m^{**,c}\partial_k^{cp}v_m^{**,p}.
	\label{eqn.si.omega.2}
\end{equation}
Inserting Eqs. \eqref{eqn.si.omega.2}, \eqref{eqn.devG.4} and \eqref{eqn.si.sigma.omega} into the discrete momentum equation  \eqref{eqn.mom.disc.3} leads to
\begin{eqnarray}
	\rho^{*,p}v_i^{**,p}&=& \rho^{*,p}v_i^{*,p}-\Delta t \, \partial^{pc}_k \left(  \sigma_{ik}^{**,c} + \omega_{ik}^{***,c} \right) \nonumber \\
	&=& \rho^{*,p}v_i^{*,p}-\Delta t \partial_k^{pc} \rho^{n,c} c_s^2 G_{il}^{n,c}\ISigma{la}\devG_{ak}^{*,c}-\Delta t^2 \partial_k^{pc} \frac{2 \left(\rho^{n,c}\right)^2 c_s^4}{3 \theta_1^{n,c}} G_{il}^{n,c} \ISigma{lk}G_{mn}^{n,c}\devG_{nm}^{n,c} \nonumber \\
	&&+\Delta t^2 \partial_k^{pc} \rho^{n,c} c_s^2 G_{il}^{n,c}\ISigma{al}\left(G_{am}^{n,c}\partial_k^{cp} v_m^{**,p}+G_{mk}^{n,c}\partial_a^{cp} v_m^{**,p}-\frac{2}{3}\delta_{ak}G_{nm}^{n,c}\partial_n^{cp}v_m^{**,p}\right) \nonumber \\
	&&-\partial_k^{pc}\frac{\Delta t \rho^{n,c} c_h^2}{1+\frac{\Delta t}{\tau_2}}J_i^{**,c}J_k^{**,c}
	+\frac{\Delta t^2 \rho^{n,c} c_h^2}{1+\frac{\Delta t}{\tau_2}} \partial_k^{pc}J_i^{**,c}J_m^{**,c}\partial_k^{cp}v_m^{**,p}.
	\label{eq.linsys.devGJv}
\end{eqnarray}
The density in the cell centers $\rho^{n,c}$ is obtained by averaging from the density in the nodes $\rho^{n,p}$.
Defining a new rank 4 tensor
\begin{equation}
	H_{ik\,nm}^{n,c} = C_{ik\,nm}^{n,c} + \frac{\tau_2 \rho^{n,c} c_h^2}{\Delta t+\tau_2}\delta_{kn}J_i^{**,c}J_m^{**,c}\label{eqn.H.r4.1}
\end{equation}
with
\begin{equation}
	C_{ik\,nm}^{n,c} = \rho^{n,c} c_s^2 G_{il}^{n,c}\ISigma{al}\left(G_{am}^{n,c}\delta_{kn}+G_{mk}^{n,c}\delta_{an}-\frac{2}{3}\delta_{ak}G_{nm}^{n,c}\right), \label{eqn.C.r4.1}
\end{equation}
the final vector wave equation for the discrete velocity field becomes
\begin{equation}
	\rho^{*,p}v_i^{**,p} - \Delta t^2 \partial_k^{pc} H_{ik\,nm}^{n,c} \partial_n^{cp}v_m^{**,p} =  b_i^{n,p},
	\label{eqn.velocity.wave}
\end{equation}
with the right hand side
\begin{align}
	b_i^{n,p} &= \rho^{*,p}v_i^{*,p}-\Delta t \partial_k^{pc} \rho^{n,c} c_s^2 G_{il}^{n,c}\ISigma{la}\devG_{ak}^{*,c}-\Delta t^2 \partial_k^{pc} \frac{2 \left(\rho^{n,c}\right)^2 c_s^4}{3 \theta_1^{n,c}} G_{il}^{n,c} \ISigma{lk}G_{mn}^{n,c}\devG_{nm}^{n,c} \nonumber \\*
	&\quad -\partial_k^{pc}\frac{\Delta t \rho^{n,c} c_h^2}{1+\frac{\Delta t}{\tau_2}}J_i^{**,c}J_k^{**,c}.
	\label{eqn.b.rhs}
\end{align}
The coefficient matrix of this system is in general not symmetric, so that \eqref{eqn.velocity.wave} is solved with a matrix-free GMRES method, in contrast to the temperature and the pressure systems, for which the conjugate gradient method can be used.
Once the new auxiliary velocity field $v_i^{**,p}$ in the nodes is known, we compute the auxiliary nodal quantity 
\begin{equation}
	\Phi^p_{i} = \partial^{pc}_k \left(  \sigma_{ik}^{**,c} + \omega_{ik}^{***,c} \right)
\end{equation}
using the definition of the stress tensors given in \eqref{eqn.si.sigma.omega}. 
These expressions have to be evaluated, since the linear solver returns only the velocity field. 
Their computation in this form is essential for the asymptotic-preserving property of the scheme. 
We then can update the cell-centered momentum as 
\begin{equation}
	\label{eqn.GJv.mom}
	(\rho v_i)^{**,c} =	(\rho v_i)^{*,c} - \Delta t \Phi^c_{i}, 
\end{equation}
where $\Phi^c_{i}$ is the previous auxiliary nodal quantity averaged back to the cell centers. 

\subsection{Discrete pressure subsystem}
The pressure and the specific enthalpy that enter this subsystem are not evolved by any of the previous steps, but are recovered from the total energy density obtained after the temperature step. Subtracting the kinetic, the elastic and the thermal contributions from $\mathcal{E}^{**,p}$ gives the internal energy density, from which the pressure and the enthalpy follow via the ideal gas law,
\begin{equation}
	\mathcal{E}_1^{**,p} = \mathcal{E}^{**,p} - \mathcal{E}_2 - \mathcal{E}_3 - \mathcal{E}_4, \qquad
	p^{**,p} = (\gamma-1) \, \mathcal{E}_1^{**,p}, \qquad
	h^{**,p} = \frac{\mathcal{E}_1^{**,p} + p^{**,p}}{\rho^{*,p}},
	\label{eqn.p.enthalpy}
\end{equation}
where the three subtracted contributions are evaluated with the most recent values available at this stage of the algorithm, that is with $v_i^{**,p}$ and with $\devG_{ik}^{**,c}$ and $J_k^{***,c}$ averaged from the cells to the vertices, since no better approximation exists at this point of the time step.
The discrete pressure subsystem reads
\begin{subequations}
	\begin{align}
		(\rho v_i)^{n+1,c}  &= (\rho v_i)^{**,c} - \Delta t \, \partial^{cp}_i p^{n+1,p}  , \label{eqn.mom.disc.4} \\
		\frac{1}{\gamma - 1} p^{n+1,p} &= \frac{1}{\gamma - 1} p^{**,p} - \Delta t \,
		\partial^{pc}_i \left( h^{**,c} (\rho v_i)^{n+1,c} \right) , \label{eqn.energy.disc.4}
	\end{align}
\end{subequations}
where we have used the ideal gas law to express the internal energy density in terms of the pressure as $\mathcal{E}_1 = p/(\gamma - 1)$.
Substitution of \eqref{eqn.mom.disc.4} into \eqref{eqn.energy.disc.4} leads to the following discrete scalar wave equation for the pressure
\begin{equation}
	\frac{1}{\gamma - 1} p^{n+1,p}  - \Delta t^2
	\partial^{pc}_i \left( h^{**,c} \, \partial^{cq}_i p^{n+1,q} \right)
	= \frac{1}{\gamma - 1} p^{**,p} -   \Delta t \,
	\partial^{pc}_i \left( h^{**,c} \, (\rho v_i)^{**,c} \right). \label{eqn.wave.p}
\end{equation}
Again, the enthalpy in the cell center $h^{**,c}$ is obtained by averaging from the nodes.
The coefficient matrix of this system is symmetric and positive definite, and \eqref{eqn.wave.p} can therefore be solved easily via a classical matrix-free conjugate gradient method. Note that the positive definiteness is due to the term $\frac{1}{\gamma-1} p^{n+1,p}$, which contributes a strictly positive multiple of the identity matrix to the system.
The final momentum is then given by \eqref{eqn.mom.disc.4}.
All linear systems arising in the implicit stages are solved to a residual tolerance of $10^{-14}$, without preconditioning.

\subsection{Compatible update of $\A$ and $\J$ and final total energy update}\label{sec.compatible}
In the absence of source terms, the distortion field $\A$ and the specific thermal impulse $\J$ should remain curl-free for all times if they were initially curl-free. The use of our structure-preserving discrete nabla operators allows such a compatible update. The same argument carries over to linear source terms, since the curl of a source that depends linearly on the field is a linear function of the curl of that field, which vanishes whenever the field is curl-free. The compatible update for the distortion field reads as follows:
\begin{eqnarray}
   A_{ik}^{n+1,c} & = & A_{ik}^{n,c} - \Delta t \, \partial_k^{cp} \left( v_m^{n+1,p} A_{im}^{n,p} + \psi^{A,p}_i \right) 
   - \frac{\Delta t}{|\omega_c|} \sum \limits_{a \in \mathcal{N}_c} | \partial \omega_{ac} | F^{A,ac}_{ik}
   \nonumber    \\  
   && 
   - \frac{\Delta t}{|\omega_c|} \sum \limits_{p \in \mathcal{P}_c} |\omega_{pc}| v_m^{n+1,p} \left( \partial_m^{pa} A_{ik}^{n,a} -
   \partial_k^{pa} A_{im}^{n,a} \right) - \frac{\Delta t}{\theta_1^{n+1,c}} \alpha_{ik}^{n+1,c},
	\label{eqn.SP.A}
\end{eqnarray}
with 
\begin{equation}
	\psi^{A,p}_i = - \halb \ell^p \, s^p_{\max} \partial_m^{pc} A_{im}^{n,c}   
\end{equation}
and 
\begin{equation}
	F^{A,ac}_{ik} = - \halb s_{\max}^{ac} t^{ac}_k \left( A_{im}^{n,a} - A_{im}^{n,c} \right) t^{ac}_m 
	                - \halb s_{\max}^{ac} r^{ac}_k \left( A_{im}^{n,a} - A_{im}^{n,c} \right) r^{ac}_m  
\end{equation}
a compatible numerical viscosity that does not destroy the curl-free property if the field is curl-free, 
while it avoids checkerboard modes for a general field $\A$ that is not curl-free.
Here, $t^{ac}_m$ and $r^{ac}_m$ are two tangent vectors to the common edge $\partial \omega_{ac}$ between cells $\omega_a$ and $\omega_c$, which becomes a common face in three space dimensions. Since all fields are three-dimensional, also in the two-dimensional case the in-plane tangent vector $t^{ac}$ has to be properly extended to three dimensions by setting its third component to zero, while the second tangent vector is simply given by $r^{ac} = (0,0,1)$. 
In order to enforce the compatibility with the determinant constraint one may uniformly rescale all components of $A_{ik}^{n+1,c}$ so that
$|\A^{n+1,c}| = \rho^{n+1,c} / \rho_0$ holds, but in general this rescaling violates the curl-free property. For alternative compatible discretizations that verify the determinant constraint, see \cite{BoscheriGPRGCL}, which however is not curl-free, and \cite{HTCLagrangeGPR}, which satisfies both the algebraic constraint and the curl-free property by construction. 
As a consequence of this observation, we proceed as follows in the remainder of this paper: for simulations in the viscous relaxation regime, where the nonlinear relaxation source does not preserve the curl-free involution at the continuous level anyway, we apply the rescaling in order to enforce the algebraic constraint, since it is required for the correct Navier--Stokes asymptotics presented below. By contrast, in the homogeneous elastic solid regime, where the distortion field is curl-free, no rescaling is applied and the compatible update is retained unchanged.

Similarly, the final update for the thermal impulse is given by
\begin{eqnarray}
   J_{k}^{n+1,c} &=& J_{k}^{n,c} - \Delta t \, \partial_k^{cp} \left( v_m^{n+1,p} J_{m}^{n,p} + T^{**,p} + \psi^{J,p} \right) 
   - \frac{\Delta t}{|\omega_c|} \sum \limits_{a \in \mathcal{N}_c} | \partial \omega_{ac} | F^{J,ac}_{k}
   \nonumber    \\  
   && 
   - \frac{\Delta t}{|\omega_c|} \sum \limits_{p \in \mathcal{P}_c}  |\omega_{pc}| v_m^{n+1,p} \left( \partial_m^{pa} J_{k}^{n,a} -
\partial_k^{pa} J_{m}^{n,a} \right) - \frac{\Delta t}{\theta_2^{n+1,c}} \beta_{k}^{n+1,c},
	\label{eqn.SP.J}
\end{eqnarray}
with  
\begin{equation}
	\psi^{J,p} = - \halb \ell^p \, s^p_{\max} \partial_m^{pc} J_{m}^{n,c}   
\end{equation}
and 
\begin{equation}
	F^{J,ac}_{k} = - \halb s_{\max}^{ac} t^{ac}_k \left( J_{m}^{n,a} - J_{m}^{n,c} \right) t^{ac}_m 
	- \halb s_{\max}^{ac} r^{ac}_k \left( J_{m}^{n,a} - J_{m}^{n,c} \right) r^{ac}_m.  
\end{equation}
The starred values of $\A$ and $\J$ computed in the previous subsystems are intermediate auxiliary quantities, which serve to build the discrete stress tensors and the coefficients of the implicit systems; the fields at the new time level are those given by \eqref{eqn.SP.A} and \eqref{eqn.SP.J}, which start from the values at time $t^n$.
Both relaxation sources are discretized implicitly. For the specific thermal impulse this is simple, since $\beta_k$ is linear in $J_k$ and $\theta_2$ depends only on the density, which is already final after the convective step, so that \eqref{eqn.SP.J} can be solved directly for $J_k^{n+1,c}$. For the distortion field, on the other hand, both $\alpha_{ik}$ and $\theta_1$ depend nonlinearly on $\A^{n+1,c}$, so that \eqref{eqn.SP.A} constitutes a local nonlinear system in each cell, which is solved by Newton's method with a tolerance of $10^{-7}$.
The discrete total energy is finally updated as
\begin{equation}
	\mathcal{E}^{n+1,p} = \mathcal{E}^{**,p} - \Delta t \, \partial_k^{pc} \left( h^{**,c} \left( \rho v_k \right)^{n+1,c} + (\sigma_{ik}^{**,c} + \omega_{ik}^{***,c}) v^{n+1,c}_i\right)
	\label{eqn.SP.E}
\end{equation}
Here, $\sigma_{ik}^{**,c}$ and $\omega_{ik}^{***,c}$ are exactly the stress tensors defined in \eqref{eqn.si.sigma.omega}, i.e. the same ones that have already been used in \eqref{eqn.GJv.mom} for the update of the momentum. 
This concludes the numerical scheme.

\subsection{Asymptotic limits}\label{sec.asymptotic}
In this section, we analyze the asymptotic consistency of the numerical scheme with the viscous Navier-Stokes stress tensor as $\tau_1 \to 0$ and the Fourier law of heat conduction as $\tau_2 \to 0$. The behavior of semi-implicit staggered schemes of this type in the low Mach number limit has already been analyzed in \cite{SIGPR} and is therefore not repeated here. 

\subsubsection{The viscous Navier-Stokes stress tensor}

Since the focus lies on the viscous effects, without loss of generality, we omit the influence of the thermal impulses by setting $c_h = 0$. 
Note that averaging nodal values to cell-centered values and vice versa does not change the order in $\tau_1$ of the respective expansions under suitable boundary conditions. 
Following \cite{FourSplit}, as $\tau_1 \to 0$ the discrete relaxation factor 
$\theta_1^{n,c} = \rho^{n,c} \frac{\tau_1 c_s^2}{3} |\G^{n,c}|^{-5/6} \to 0$ and omitting higher order terms in $\tau_1$ the inverse matrix $\left(\Sigma^{n,c}\right)^{-1}$ can be approximated by 
\begin{equation}
	\ISigma{il} = \frac{1}{6} \frac{\tau_1}{\Delta t} |\G^{n,c}|^{-5/6} (G_{il}^{n,c})^{-1}.
\end{equation}
Applying the Chapman-Enskog expansion for small $\tau_1 \ll 1$ of the discrete metric tensor given by 
\begin{equation}
	\G^{**,c} = \G_0^{**,c} + \tau_1 \G_1^{**,c} + \mathcal{O}(\tau_1^2)
	\label{eqn.APG2}
\end{equation}
on the evolution equation \eqref{eqn.devG.1}, we obtain immediately at leading order 
\begin{equation}
	\mathring{\G}_0^{**,c} = 0, \quad \text{which implies} \quad \G^{**,c} = g^{**,c} \mathbf{I} + \mathcal{O}(\tau_1),
	\label{eqn.ChapmanDevG}
\end{equation}
where $g^{**,c} = |\G^{**,c}|^{1/3}$. 
Starting from well-prepared initial data as derived in \cite{GPRmodelMHD}, i.e. 
\begin{equation}\label{eqn.AP.Gn}
	\G^{n,c} = g^{n,c} \mathbf{I} + \tau_1 \G_1^{n,c} + \mathcal{O}(\tau_1^2), \quad \mathring{\G}_0^{n,c} = 0,
\end{equation}
in the discrete convective subsystem, the metric tensor is still well-prepared. 
In the viscous relaxation regime the rescaling described in Section \ref{sec.compatible} enforces the determinant constraint $|\A^{n,c}| = \rho^{n,c}/\rho_0$ already at the previous time level, hence $|\G^{n,c}| = \left(\rho^{n,c}/\rho_0\right)^2$, and since the well-prepared data \eqref{eqn.AP.Gn} give $|\G^{n,c}| = \left(g^{n,c}\right)^3 + \mathcal{O}(\tau_1)$, one obtains
\begin{equation}
	g^{n,c} = \left(\frac{\rho^{n,c}}{\rho_0}\right)^{\frac{2}{3}} + \mathcal{O}(\tau_1), \qquad
	\rho^{n,c} \left(g^{n,c}\right)^{-\frac{3}{2}} = \rho_0 + \mathcal{O}(\tau_1),
	\label{eqn.AP.rho0}
\end{equation}
which is what produces the constant reference density $\rho_0$ in the viscosity coefficient obtained below. 
Now, using \eqref{eqn.APG2} and \eqref{eqn.ChapmanDevG} in the update of the auxiliary velocity equation \eqref{eq.linsys.devGJv} we have 
\begin{eqnarray}
	\rho^{*,p}v_i^{**,p} \!
	%&=& \! b_i^{n,p}+\Delta t^2 \partial_k^{pc} \frac{\rho^{n,c} c_s^2 \tau_1}{6\Delta t} |\G^{n,c}|^{-\frac{5}{6}}  G_{il}^{n,c}\left(G_{al}^{n,c}\right)^{-1}\left(G_{am}^{n,c}\partial_k^{cp} v_m^{**,p}+G_{mk}^{n,c}\partial_a^{cp} v_m^{**,p}-\frac{2}{3}\delta_{ak}G_{nm}^{n,c}\partial_n^{cp}v_m^{**,p}\right) \nonumber \\
	&=& \rho^{*,p} v_i^{*,p}+\Delta t \partial_k^{pc} \frac{1}{6} \rho_0 c_s^2 \tau_1  \left(\partial_k^{cp} v_i^{**,p}+\partial_i^{cp}v_k^{**,p}-\frac{2}{3}\delta_{ik}\partial_m^{cp}v_m^{**,p}\right) + \mathcal{O}(\tau_1^2),
\end{eqnarray}
yielding the compatibility with the classical compressible Navier-Stokes stress tensor for the auxiliary velocity field $v_i^{**,p}$ as $\tau_1 \to 0$.
To verify the compatibility with the Navier-Stokes stress tensor for the centered momentum, we turn to the stress tensor $\sigma_{ik}^{**,c}$ in \eqref{eqn.si.sigma.omega} used in the update \eqref{eqn.GJv.mom}.
Using the expansions \eqref{eqn.AP.Gn} and \eqref{eqn.APG2} we have
\begin{equation}
	\label{eqn.APsigma}
\sigma_{ik}^{**,c} = \tau_1\rho^{n,c}c_s^2 g^{n,c}\delta_{il}\devG_{lk,1}^{**,c} + \mathcal{O}(\tau_1^2).
\end{equation}
Inserting the same expansions into equation \eqref{eqn.devG.1} and collecting the terms at leading order, one obtains
\begin{equation}
	g^{n,c} \left(\partial_k^{cp}v_i^{**,p} + \partial_i^{cp}v_k^{**,p} - \frac{2}{3} \delta_{ik} \partial_m^{cp}v_m^{**,p}\right) = - 6 (g^{n,c})^{7/2} \delta_{il}\devG_{lk,1},
\end{equation}
Inserting the expression into \eqref{eqn.APsigma}, we obtain 
\begin{equation}
	\sigma_{ik}^{**,c} = - \frac{1}{6} \tau_1 c_s^2 \rho_0 \left(\partial_k^{cp}v_i^{**,p} + \partial_i^{cp}v_k^{**,p} - \frac{2}{3} \delta_{ik} \partial_m^{cp}v_m^{**,p}\right) + \mathcal{O}(\tau_1^2),
\end{equation}
which yields the classical compressible Navier-Stokes stress tensor. 
Thus the update for the cell-centered momentum \eqref{eqn.GJv.mom} in the limit $\tau_1 \to 0$ is consistent with the compressible Navier-Stokes stress. 

\subsubsection{The Fourier heat conduction}
Following \cite{GPRmodelMHD}, well-prepared initial data for the thermal impulse is given by $J_k^{n,c} = \tau_2 J_{k,1}^{n,c} + \mathcal{O}(\tau_2^2)$. 
The convective step does not influence the expansion in $\tau_2$, thus we have $J_k^{*,c} = \tau_2 J_{k,1}^{*,c} + \mathcal{O}(\tau_2^2)$.
Assuming the expansion $J_k^{**,c} = J_{k,0}^{**,c} + \tau_2 J_{k,1}^{**,c} + \mathcal{O}(\tau_2^2)$, from \eqref{eqn.heatflux.disc.2} we obtain the conditions $J_{k,0}^{**,c} = 0$ and $J_{k,1}^{**,c} = -\partial_k^{cp} T^{**,p}$ which yield $J_k^{**,c} = - \tau_2 \, \partial_k^{cp} T^{**,p} + \mathcal{O}(\tau_2^2)$. 
Inserting this relation into the definition of the heat flux yields 
\begin{equation}
	q_k^{**,c} = \rho^{*,c} c_h^2 T^{**,c} J_k^{**,c} = - \rho^{*,c} c_h^2 \tau_2 T^{**,c} \partial_k^{cp} T^{**,p} + \mathcal{O}(\tau_2^2)
	\label{eqn.APheatFlux}
\end{equation}
which is a consistent discretization of the Fourier heat conduction as $\tau_2 \to 0$. 
The subsequent $\mathring{\mathbf{G}}$-$\J$-$\v$ subsystem, detailed in \cite{FourSplit}, preserves the expansion of $J_k^{***,c}$ in $\tau_2$, and so does the curl-free update \eqref{eqn.SP.J}, so that the thermal impulse is still well-prepared at the new time level.
Note that the final total energy update \eqref{eqn.SP.E} utilizes the consistent heat flux \eqref{eqn.APheatFlux} and therefore in the stiff relaxation limit $\tau_2 \to 0$ we recover a consistent discretization of the Fourier heat flux independently of $\Delta t$ at the new time level.

\subsection{Summary of the algorithm}\label{sec.summary}
Since the individual stages of the method have been presented one after the other over the preceding subsections, we collect here the complete sequence of operations that advances the discrete solution from time $t^n$ to time $t^{n+1}$. The equations are repeated in unnumbered form, together with the number of the original equation in each case, so that every formula of this summary can still be referred to unambiguously.

\subsubsection{Convective subsystem}
The transport of the conserved quantities with the material velocity of the medium is the only stage of the algorithm that is treated explicitly. The vertex-based mass density and total energy density are advanced by the special compatible discretization on the dual mesh 
\begin{equation*}
	\rho^{*,p} = \rho^{n,p} - \Delta t \, \partial_k^{pc} \, f^{\rho,c}_k, \qquad
	\mathcal{E}^{*,p} = \mathcal{E}^{n,p} - \Delta t \, \partial_k^{pc} \, f^{\mathcal{E},c}_k,
\end{equation*}
see \eqref{eqn.fv.mass} and \eqref{eqn.fv.rho}, with the multidimensional Rusanov-type fluxes \eqref{eqn.rusanov.mass} and \eqref{eqn.rusanov.md}, while the cell-centered momentum, the distortion field and the specific thermal impulse are transported by the edge-based finite volume scheme \eqref{eqn.mom.transport}, \eqref{eqn.A.transport} and \eqref{eqn.J.transport}. Since none of the subsequent stages is explicit, the admissible time step is governed by the material speed alone, see \eqref{eqn.cfl}.

\subsubsection{Temperature subsystem}
The nodal temperature $T^{*,p}$ from which this stage starts is recovered from $\mathcal{E}^{*,p}$ by means of the equation of state. Heat conduction is then taken into account implicitly, so that the propagation speed of the heat waves does not restrict the time step, which requires the solution of the scalar wave equation \eqref{eqn.wave.T},
\begin{equation*}
	\left( 1+\frac{\Delta t}{\tau_2} \right) \frac{c_v}{T^{n,p} c_h^2} \, T^{**,p} - \Delta t^2
	\partial^{pc}_k \left( \partial^{cq}_k T^{**,q} \right)
	= \left( 1+\frac{\Delta t}{\tau_2} \right) \frac{c_v}{T^{n,p} c_h^2} \, T^{*,p} - \Delta t \,
	\partial^{pc}_k \left( J_k^{*,c} \right),
\end{equation*}
whose coefficient matrix is symmetric and positive definite, so that a matrix-free conjugate gradient method can be employed. The thermal impulse and the total energy are subsequently updated according to \eqref{eqn.heat.disc.final} and \eqref{eqn.heat.disc.final.E}.

\subsubsection{$\mathring{\mathbf{G}}$-$\J$-$\v$ subsystem}
This stage evolves the trace-free part $\mathring{\mathbf{G}}$ of the metric tensor and not $\G$ itself, which is what permits the asymptotic-preserving property with respect to the Navier-Stokes stress tensor to be established and which is therefore an essential ingredient of the method. After the elimination of $\devG_{ik}^{**,c}$ and $J_k^{***,c}$ in favour of the velocity, the discrete vector wave equation \eqref{eqn.velocity.wave},
\begin{equation*}
	\rho^{*,p}v_i^{**,p} - \Delta t^2 \partial_k^{pc} H_{ik\,nm}^{n,c} \partial_n^{cp}v_m^{**,p} =  b_i^{n,p},
\end{equation*}
remains to be solved, with the rank 4 tensor $H_{ik\,nm}^{n,c}$ given in \eqref{eqn.H.r4.1} and the right hand side in \eqref{eqn.b.rhs}. Once the nodal velocity field is available, $\devG_{ik}^{**,c}$ and $J_k^{***,c}$ are recovered from \eqref{eqn.devG.4} and \eqref{eqn.J.disc.5}, the two stress tensors are evaluated directly according to \eqref{eqn.si.sigma.omega}, which is essential for the asymptotic-preserving property of the scheme, and the cell-centered momentum is updated by \eqref{eqn.GJv.mom}.

\subsubsection{Pressure subsystem}
The last implicit stage couples the momentum with the pressure gradient and with the enthalpy flux in the energy equation. Eliminating the momentum leads to the discrete scalar wave equation \eqref{eqn.wave.p} for the new pressure,
\begin{equation*}
	\frac{1}{\gamma - 1} p^{n+1,p}  - \Delta t^2
	\partial^{pc}_i \left( h^{**,c} \, \partial^{cq}_i p^{n+1,q} \right)
	= \frac{1}{\gamma - 1} p^{**,p} -   \Delta t \,
	\partial^{pc}_i \left( h^{**,c} \, (\rho v_i)^{**,c} \right),
\end{equation*}
which degenerates into the pressure Poisson equation of the incompressible limit as $M_a \to 0$ and which therefore yields a momentum update that is compatible with the low Mach number limit, see \eqref{eqn.mom.disc.4}.

\subsubsection{Compatible update of $\A$ and $\J$ and final total energy update}
The distortion field and the specific thermal impulse are not advanced within the four subsystems, but in one single compatible step, see \eqref{eqn.SP.A} 
\begin{eqnarray*}
		A_{ik}^{n+1,c} & = & A_{ik}^{n,c} - \Delta t \, \partial_k^{cp} \left( v_m^{n+1,p} A_{im}^{n,p} + \psi^{A,p}_i \right) 
		- \frac{\Delta t}{|\omega_c|} \sum \limits_{a \in \mathcal{N}_c} | \partial \omega_{ac} | F^{A,ac}_{ik}
		\nonumber    \\  
		&& 
		- \frac{\Delta t}{|\omega_c|} \sum \limits_{p \in \mathcal{P}_c} |\omega_{pc}| v_m^{n+1,p} \left( \partial_m^{pa} A_{ik}^{n,a} -
		\partial_k^{pa} A_{im}^{n,a} \right) - \frac{\Delta t}{\theta_1^{n+1,c}} \alpha_{ik}^{n+1,c},
\end{eqnarray*}
and \eqref{eqn.SP.J}, 
\begin{eqnarray*}
	J_{k}^{n+1,c} &=& J_{k}^{n,c} - \Delta t \, \partial_k^{cp} \left( v_m^{n+1,p} J_{m}^{n,p} + T^{**,p} + \psi^{J,p} \right) 
	- \frac{\Delta t}{|\omega_c|} \sum \limits_{a \in \mathcal{N}_c} | \partial \omega_{ac} | F^{J,ac}_{k}
	\nonumber    \\  
	&& 
	- \frac{\Delta t}{|\omega_c|} \sum \limits_{p \in \mathcal{P}_c}  |\omega_{pc}| v_m^{n+1,p} \left( \partial_m^{pa} J_{k}^{n,a} -
	\partial_k^{pa} J_{m}^{n,a} \right) - \frac{\Delta t}{\theta_2^{n+1,c}} \beta_{k}^{n+1,c},
\end{eqnarray*}
which preserves their curl-free character at the discrete level in case the source terms are absent or linear. All components of $A_{ik}^{n+1,c}$ may be uniformly rescaled so that the determinant constraint $|\A^{n+1,c}| = \rho^{n+1,c}/\rho_0$ is met, but this in general then violates the curl-free property.  The time step is completed by the total energy update \eqref{eqn.SP.E},
\begin{equation*}
	\mathcal{E}^{n+1,p} = \mathcal{E}^{**,p} - \Delta t \, \partial_k^{pc} \left( h^{**,c} \left( \rho v_k \right)^{n+1,c} + (\sigma_{ik}^{**,c} + \omega_{ik}^{***,c}) v^{n+1,c}_i\right).
\end{equation*}

% % % % % % % % % % % % % % % % % % % % % % % % % % % % % %
%            Numerical experiments                        %
% % % % % % % % % % % % % % % % % % % % % % % % % % % % % %

\section{Numerical results}\label{sec.results}
In the following, the new four-split structure-preserving semi-implicit finite volume scheme on unstructured meshes presented in this paper is abbreviated as UnSPSIFV in the legends of the figures.

\subsection{Taylor-Green vortex at low Mach number}

\begin{figure}
	\begin{center}
		\begin{tabular}{cc}
	    \includegraphics[trim=10 10 10 10,clip,width=0.49\textwidth]{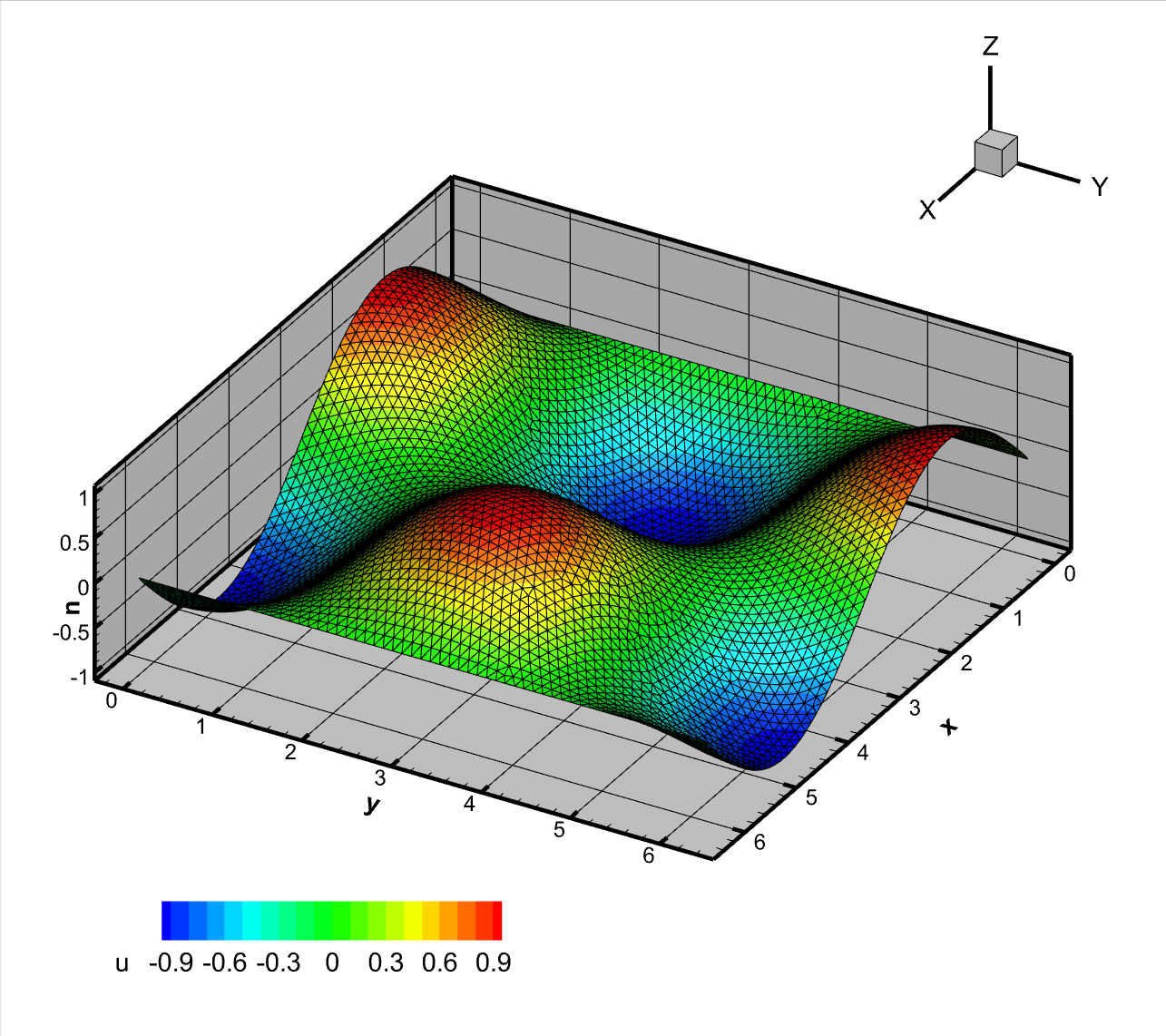} &
	    \includegraphics[trim=10 10 10 10,clip,width=0.49\textwidth]{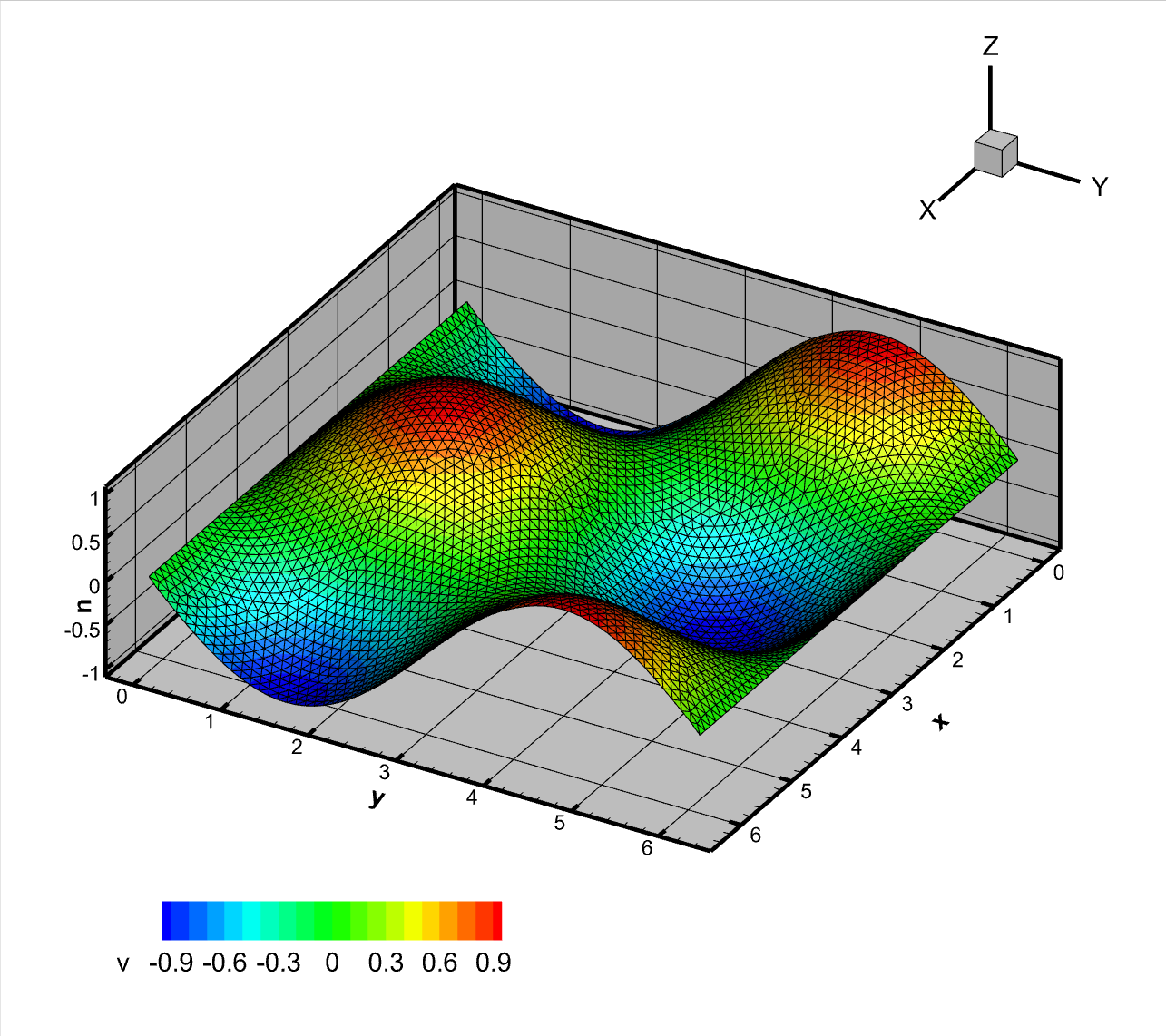}
 	    \end{tabular}
    \end{center}
    \caption{Computed velocity component $u$ (left) and $v$ (right) at the final time $t=1.0$ with the adopted unstructured mesh.}
    \label{fig.tgv.grid}
\end{figure}

\begin{figure}[htbp]
    	    \begin{center}
        		\begin{tabular}{cc}
            			\includegraphics[width=0.45\textwidth]{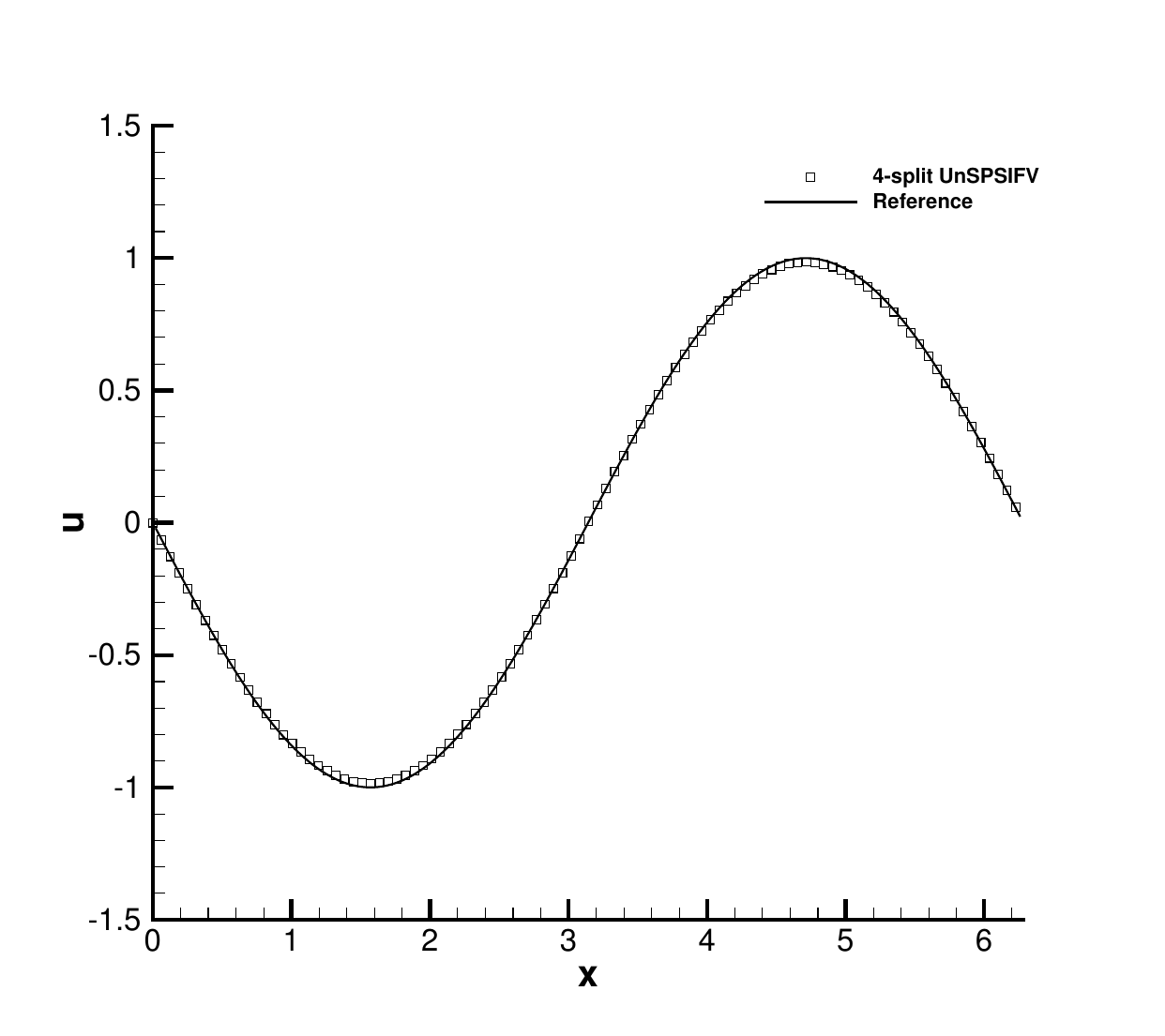}   &
            			\includegraphics[width=0.45\textwidth]{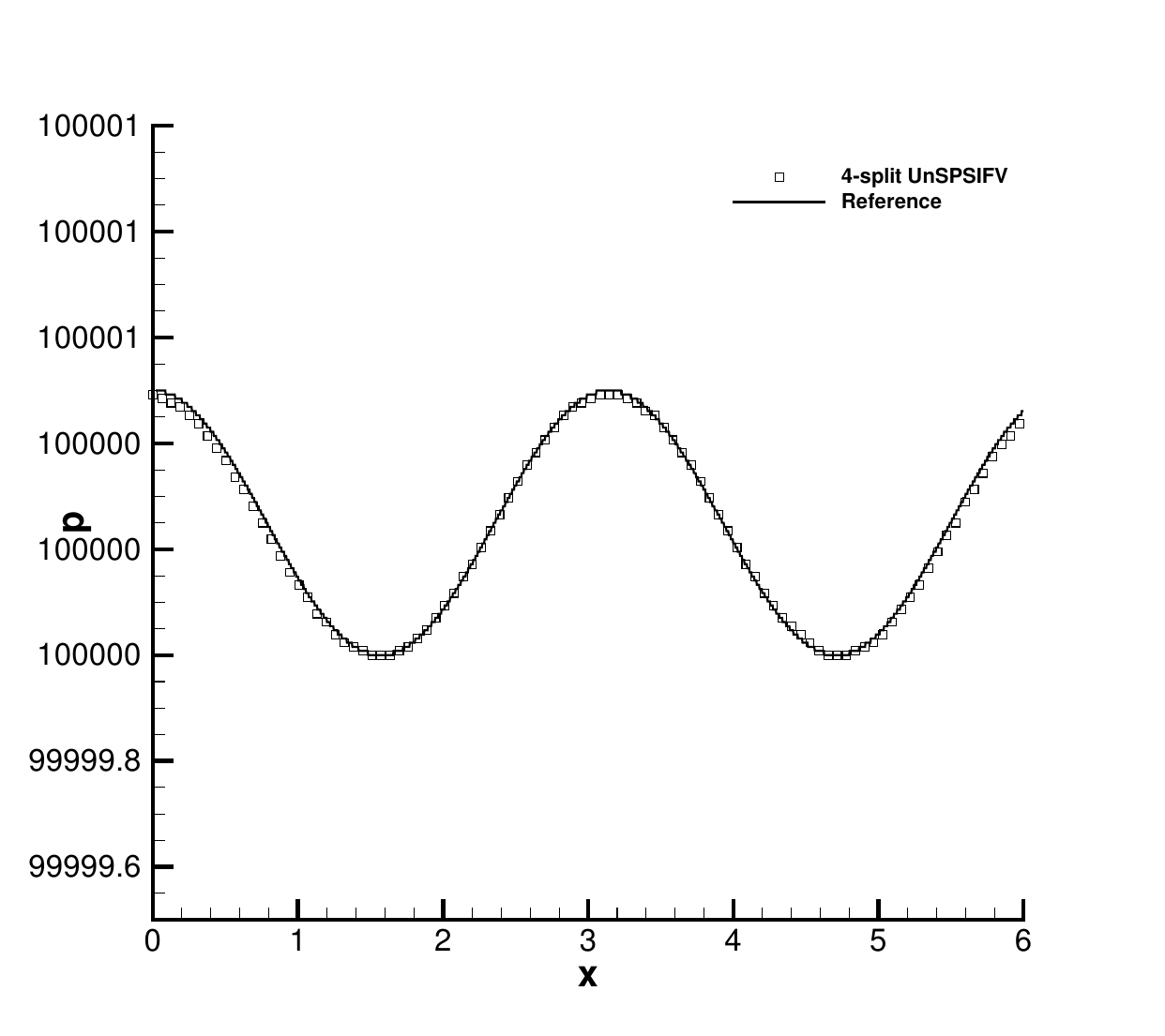}
            		\end{tabular}
        	\end{center}
    \caption{Numerical solution of the GPR model for the Taylor-Green vortex problem in the low Mach number limit and in the stiff relaxation limit with $\tau_1 = 10^{-8}$ and $\tau_2 = 10^{-10}$ at time $t=1.0$ using the new four-split scheme. 1D cuts along the $x$-axis and comparison with the exact solution of the incompressible Navier-Stokes equations for the velocity components $u$ (left) and the pressure $p$ (right).  }
    \label{fig.tgv}
\end{figure}

To assess the low Mach number properties of the novel four-split scheme on unstructured grids we consider the two-dimensional Taylor-Green vortex
\begin{eqnarray}
    \rho(x,y,t) &=& \rho_0, \label{eq:TG_rho} \\
    u(x,y,t)&=&\sin(x)\cos(y)e^{-2\nu t}, \label{eq:TG_0} \\
    v(x,y,t)&=&-\cos(x)\sin(y)e^{-2\nu t}, \label{eq:TG_1} \\
    p(x,y,t)&=& p_0 + \frac{1}{4}(\cos(2x)+\cos(2y))e^{-4\nu t},
    \label{eq:TG_2}
\end{eqnarray}
which is an exact solution of the incompressible Navier-Stokes equations and as such satisfies the compressible system as the acoustic Mach number $M_a$ tends to zero.
The computational domain $\Omega = [0,2\pi]^2$ is equipped with periodic boundary conditions and is discretized with $8092$ triangular elements. The final time of the simulation is set to $t = 1.0$.
The material and relaxation parameters used are
$\gamma=1.4$, $\rho_0=1$, $\tau_1=10^{-8}$, $\tau_2=10^{-10}$, so that $\nu=1.667\cdot10^{-3}$. Finally $c_v=c_p/\gamma$ with $c_p=1004$, $c_s=1000$, $c_h=100$. The initial condition is given by \eqref{eq:TG_rho}--\eqref{eq:TG_2} for density, velocity and pressure, where we set $p_0=10^5$, leading to a low
acoustic Mach number regime of $M_a=0.0027$ and a low shear Mach number regime of $M_s=10^{-3}$.
The initial distortion field is given by a material at rest $\mathbf{A}=\mathbf{I}$ and the specific thermal impulse at the beginning of the simulation is $\mathbf{J}=0$.
The resulting velocity field, as well as the adopted triangular mesh, is reported in Figure \ref{fig.tgv.grid}.
As can be seen from Figure \ref{fig.tgv}, the numerical solution obtained with the new scheme is in good agreement with the reference solution given by the exact solution of the incompressible Navier-Stokes equations.

To assess the asymptotic-preserving property of the new scheme, i.e. that the divergence of the velocity field scales with $M_a^2$, see e.g.  \cite{KlaMaj,KlaMaj82,Klein2001,KleinMach,MunzPark,MunzDumbserRoller}, we solve the Taylor-Green vortex again up to a final time of $t=0.1$, but with increasing background pressure $p_0$ and thus with decreasing acoustic Mach number, while keeping the mesh fixed to $8092$ cells and the time step constant.
Moreover, for this test case $\rho \to \rho_0$ as $M_a \to 0$.
In Table \ref{tab.Machconv} the density $L^2$ and $L^\infty$ errors and the $L^\infty$ error of the divergence of the velocity field as a function of $p_0$ and the acoustic Mach number $M_a$ are displayed.
As expected, second order convergence of the density and divergence errors with respect to the acoustic Mach number is achieved with the new scheme.
\begin{table} [htpb]
    \caption{Convergence results with respect to the acoustic Mach number, obtained with the new four-split scheme for the GPR model in the low Mach number limit. The error norms refer to the density $\rho$ in $L^2$ and $L^\infty$ norm as well as to the divergence of the velocity in $L^\infty$ norm at a final time of $t=0.1$. The observed convergence orders refer to the acoustic Mach number $M_a$ and to the density and divergence errors, and clearly show second order convergence in the Mach number. }
    		\begin{center}
        			\renewcommand{\arraystretch}{1.1}
        			\begin{tabular}{cccccccc}
            				\hline
            				$M_a$ & $p_0$ &  ${L^2}(\rho)$ & ${L^\infty}(\rho)$ &  ${L^\infty}(\nabla \cdot \mathbf{v})$ & $\mathcal{O}^2_M(\rho)$ & $\mathcal{O}^\infty_M(\rho)$
            				& $\mathcal{O}^\infty_M(\nabla \cdot \mathbf{v})$ \\
            				\hline
								1.00E-01 & 7.14E+01	&	5.4853E-04 & 2.2043E-04	&  6.5297E-03 & 	&		&			\\
								5.00E-02 & 2.86E+02	&	2.1350E-04 & 8.2583E-05	&  1.1148E-03 & 1.4	&	1.4	&	2.6		\\
								1.00E-02 & 7.14E+03	&	7.7547E-06 & 2.9929E-06	&  3.4839E-05 & 2.1	&	2.1	&	2.2		\\
								5.00E-03 & 2.86E+04	&	1.9388E-06 & 7.4854E-07	&  6.0876E-06 & 2.0	&	2.0	&	2.5		\\
								1.00E-03 & 7.14E+05	&	7.7552E-08 & 2.9942E-08	&  2.7801E-07 & 2.0	&	2.0	&	1.9		\\
								5.00E-04 & 2.86E+06	&	1.9388E-08 & 7.4860E-09	&  7.7642E-08 & 2.0	&	2.0	&	1.8		\\
            				\hline
            			\end{tabular}
        		\end{center}
    \label{tab.Machconv}
\end{table}
\subsection{Riemann problems in the fluid and solid limit}
\begin{figure}[htbp]
    \begin{center}
        		\includegraphics[width=0.32\textwidth]{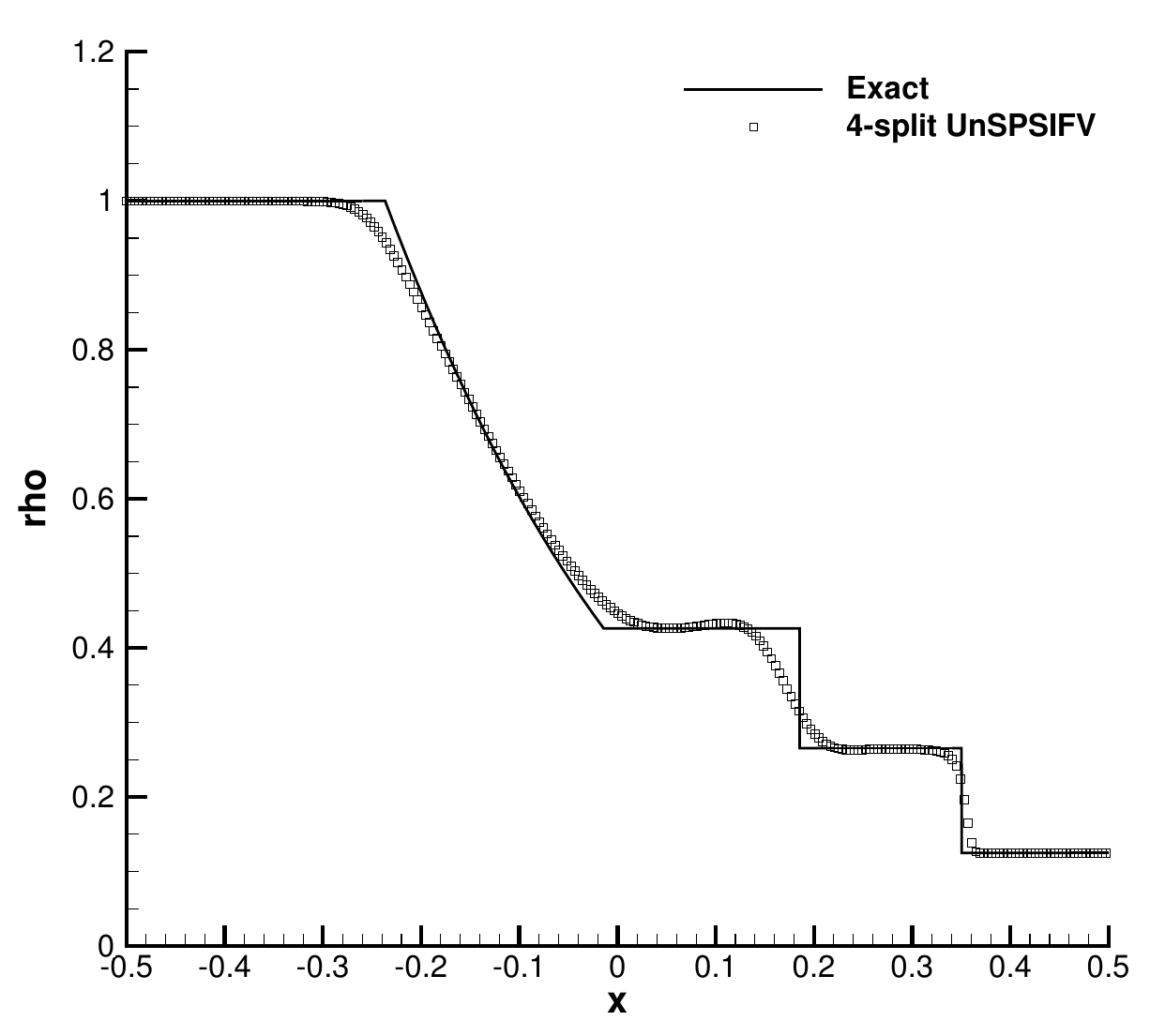}
        		\includegraphics[width=0.32\textwidth]{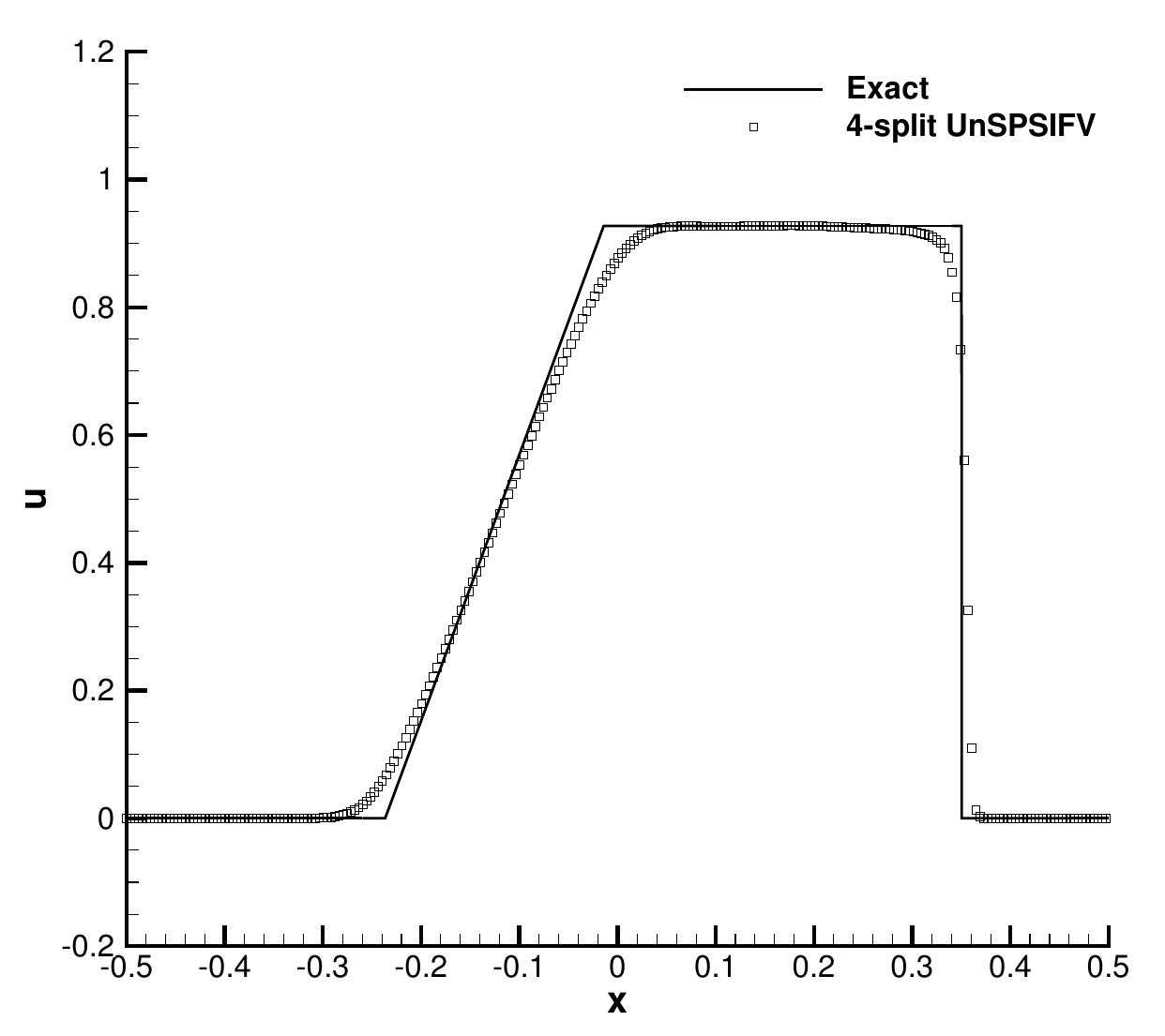}
        		\includegraphics[width=0.32\textwidth]{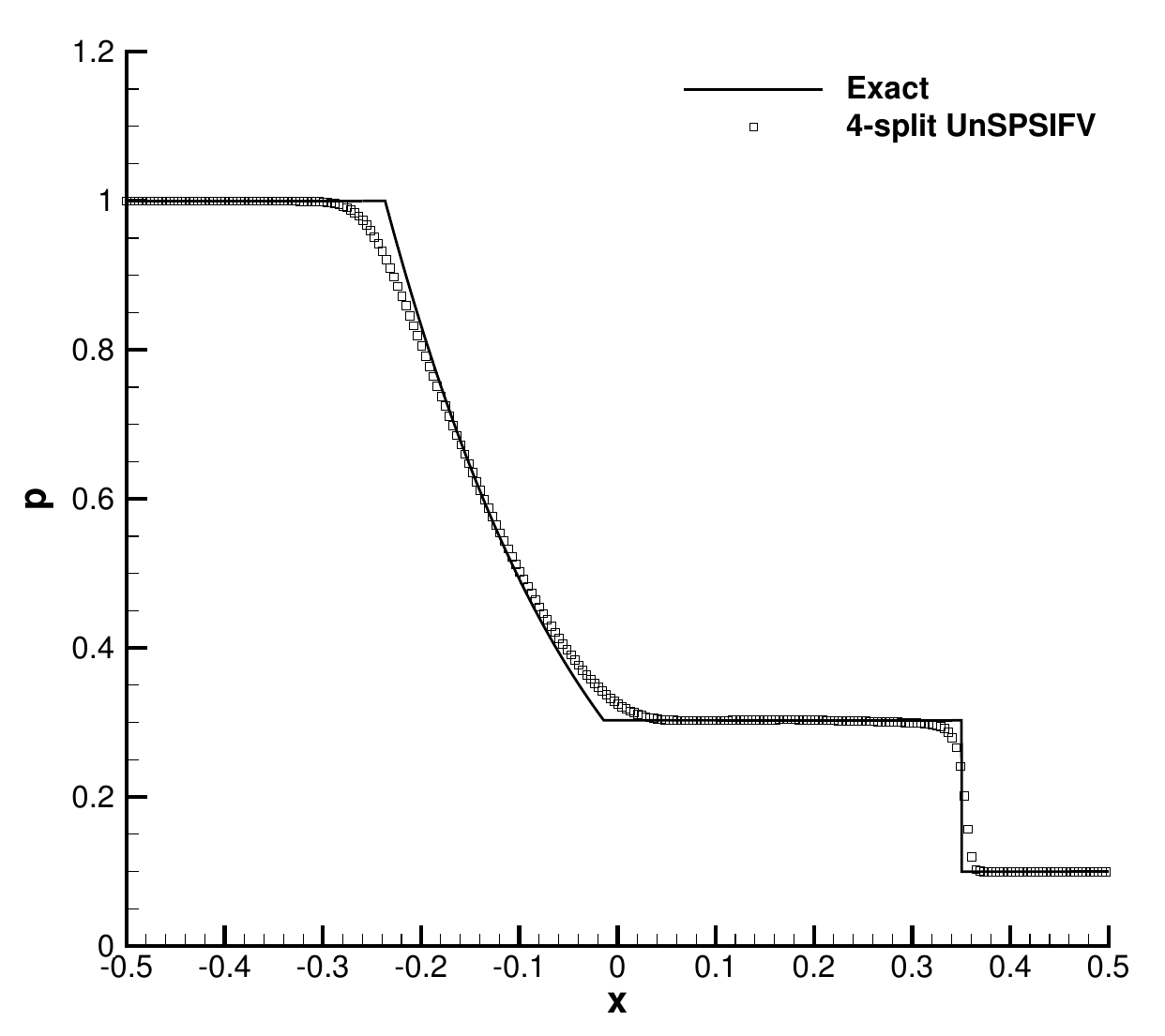}
        \caption{Exact solution of the compressible Euler equations and numerical solution of the GPR model in the stiff relaxation limit ($\tau_1 = 10^{-10}, \tau_2 = 10^{-12}$) for Riemann problem RP1 (Sod shock tube) obtained with the new four-split scheme. The density $\rho$, the velocity component $u$ and the pressure $p$ are shown at a final time of $t=0.2$.}
        \label{fig.rp1}
    \end{center}
\end{figure}

\begin{figure}[htbp]
    \begin{center}
        		\includegraphics[width=0.32\textwidth]{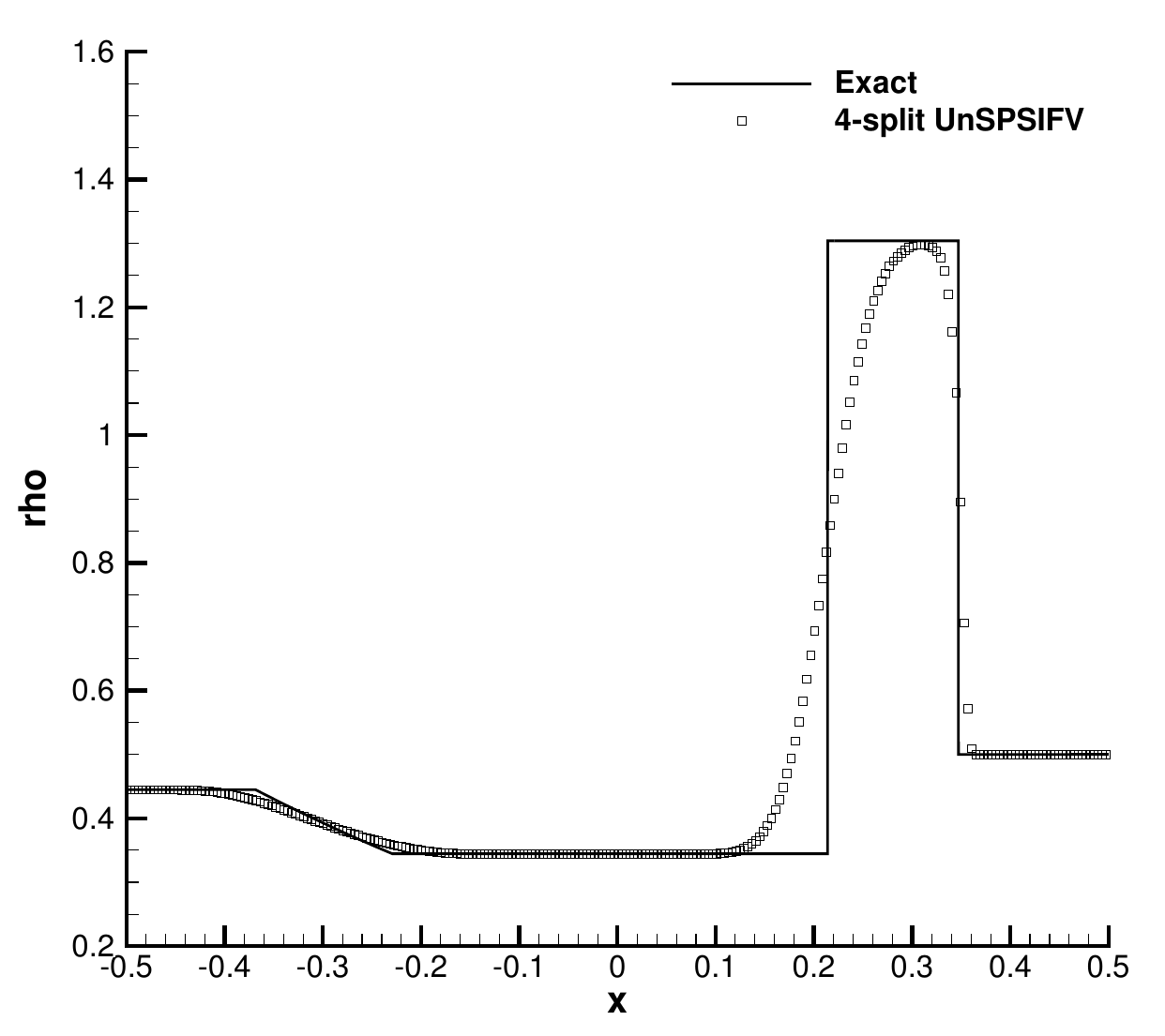}
        		\includegraphics[width=0.32\textwidth]{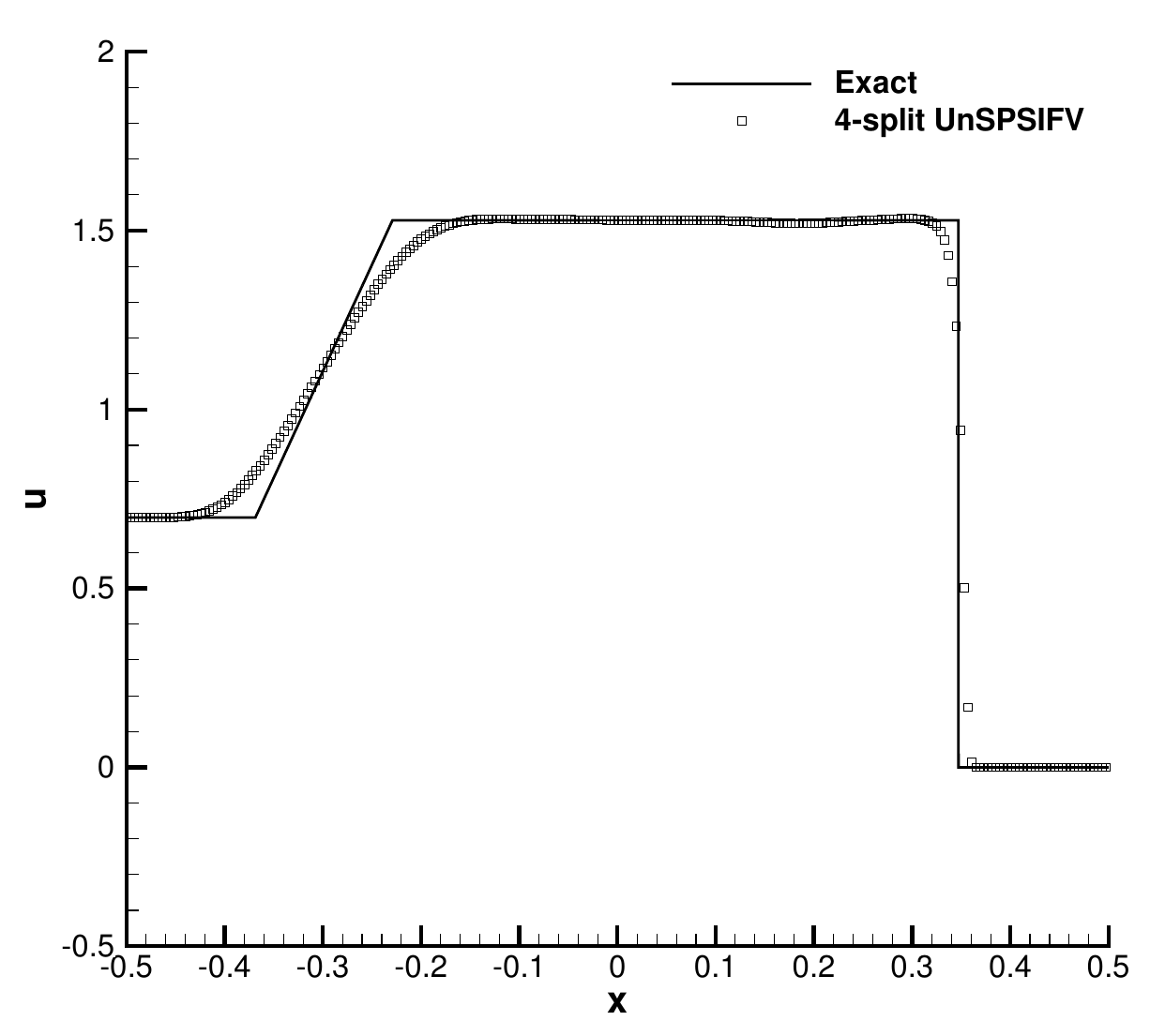}
        		\includegraphics[width=0.32\textwidth]{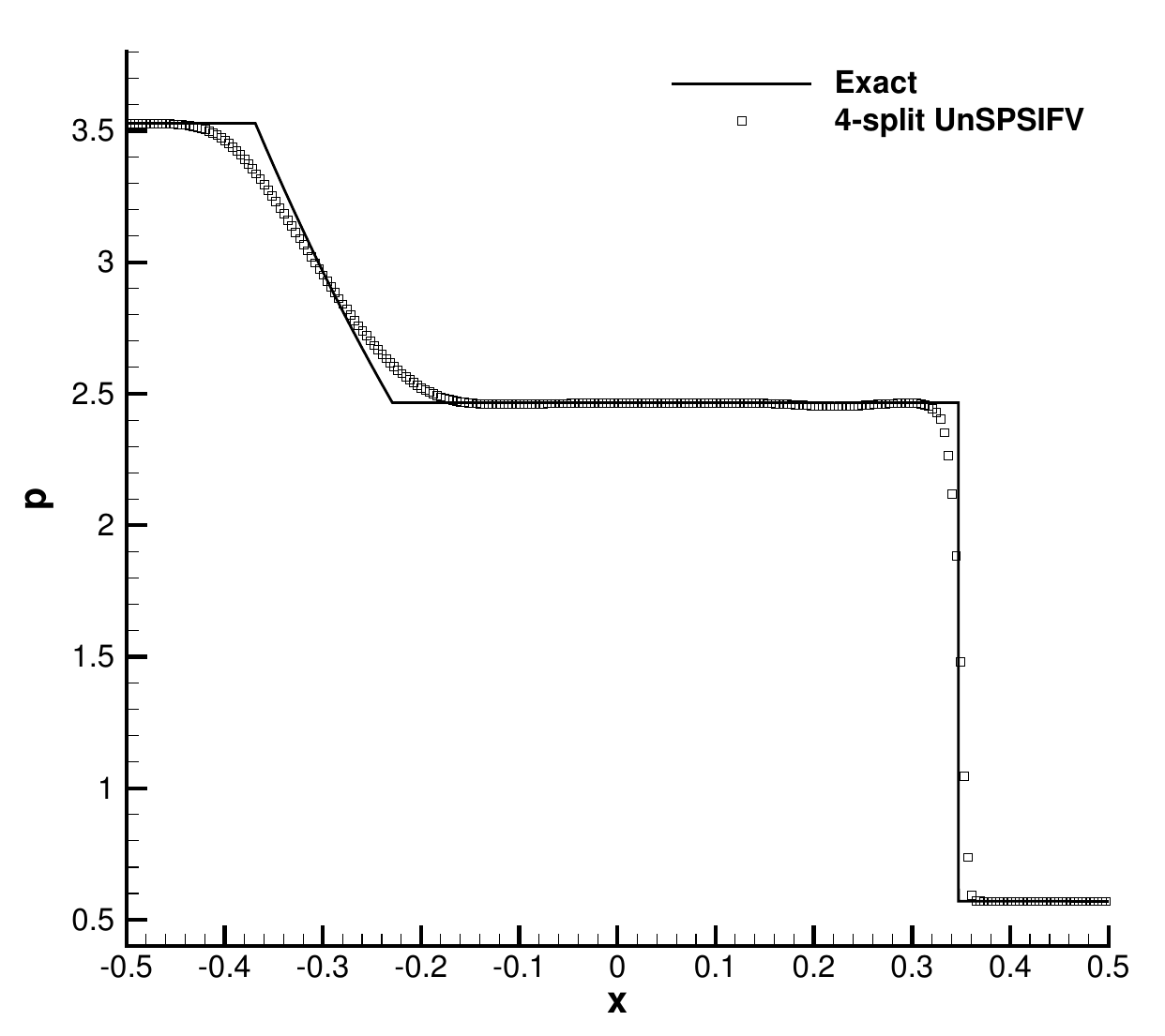}
        \caption{Exact solution of the Euler equations and numerical solution of the GPR model in the stiff relaxation limit ($\tau_1 = 10^{-10}, \tau_2 = 10^{-12}$) for Riemann problem RP2 (Lax shock tube). The density $\rho$, the velocity component $u$ and the pressure $p$ are shown at a final time of $t=0.14$.}
        \label{fig.rp2}
    \end{center}
\end{figure}

\begin{figure}[htbp]
    \begin{center}
        \begin{tabular}{ccc}
            		\includegraphics[width=0.32\textwidth]{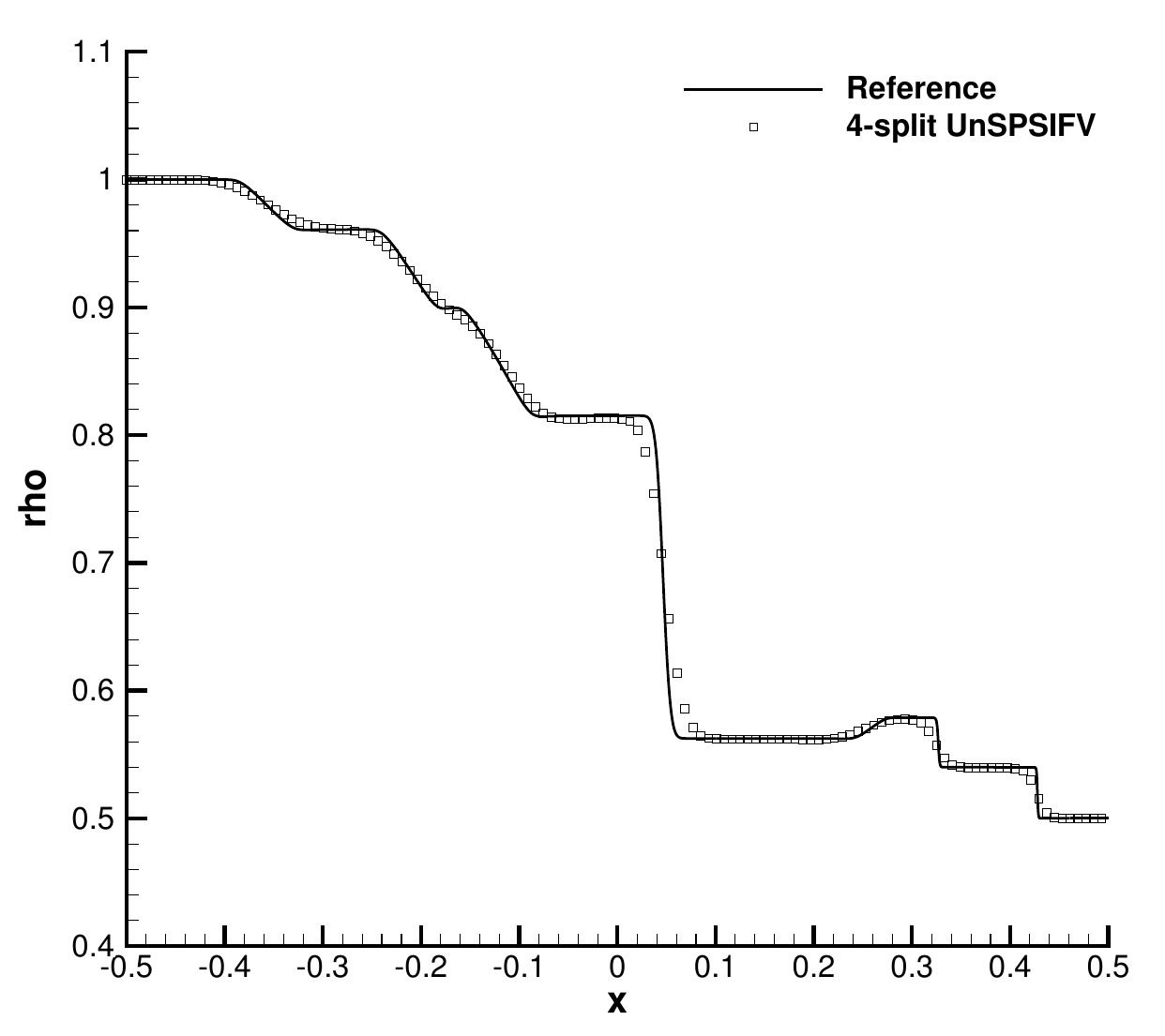}   &
            		\includegraphics[width=0.32\textwidth]{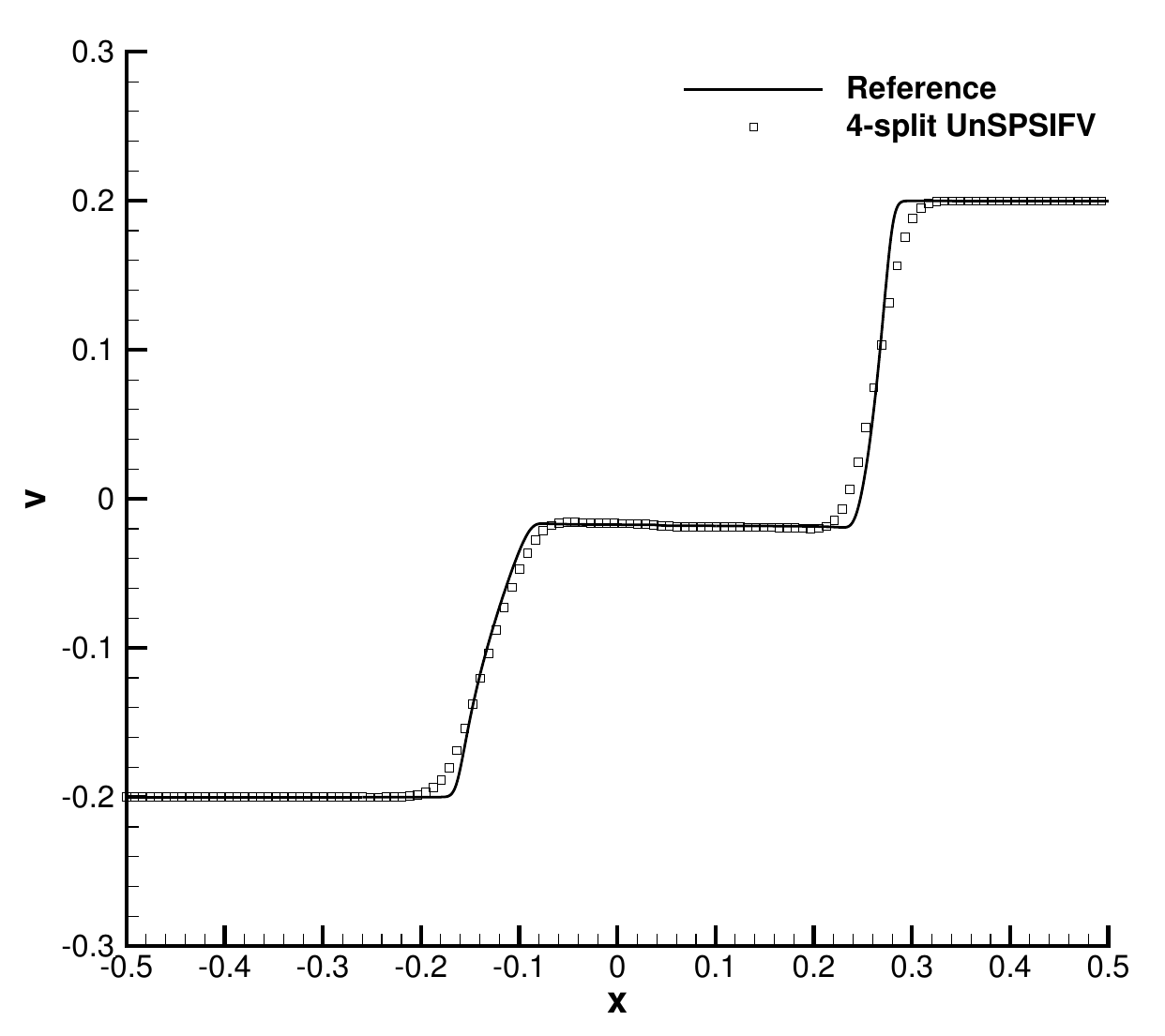}     &
            		\includegraphics[width=0.32\textwidth]{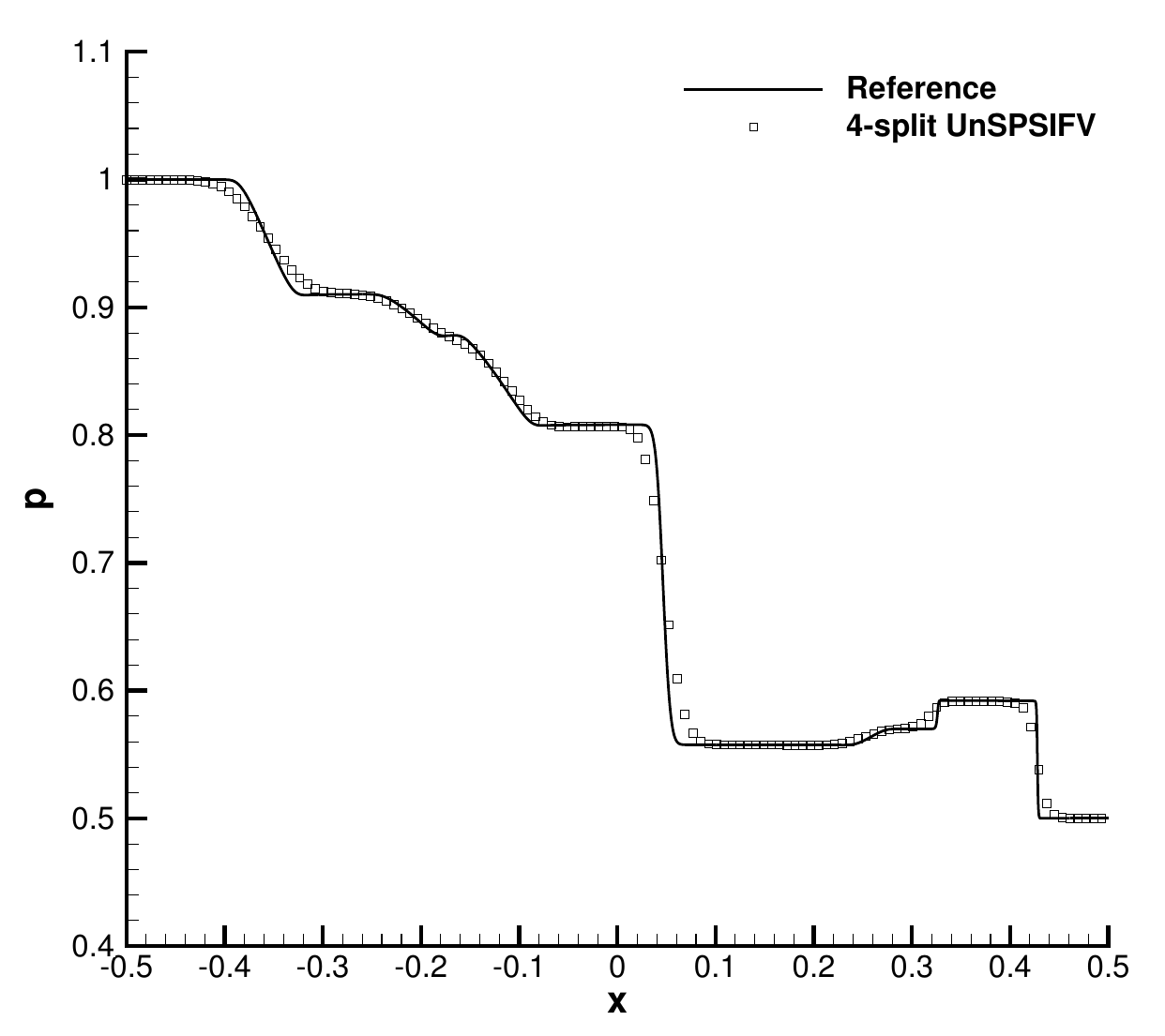}     \\
            		\includegraphics[width=0.32\textwidth]{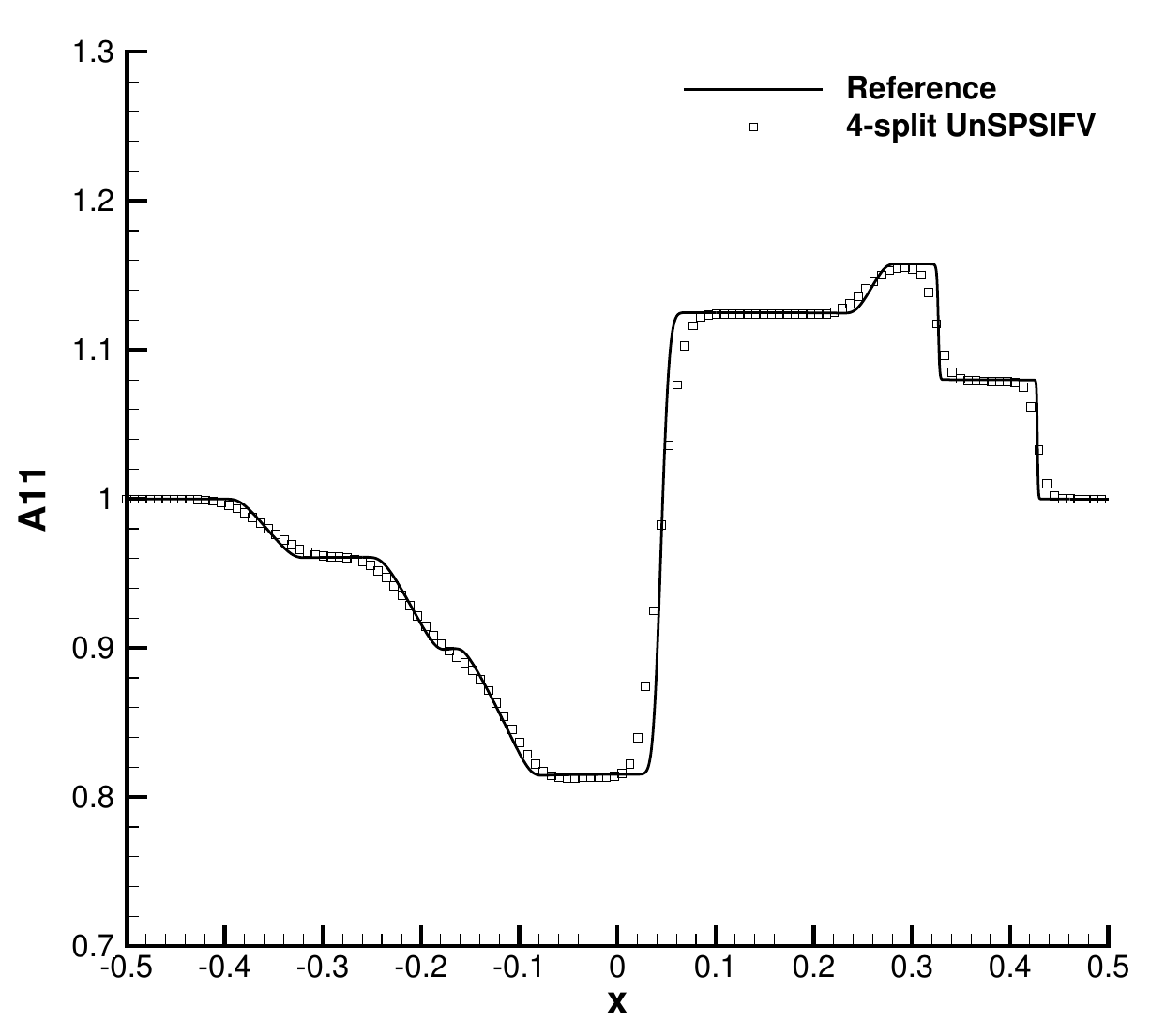}   &
            		\includegraphics[width=0.32\textwidth]{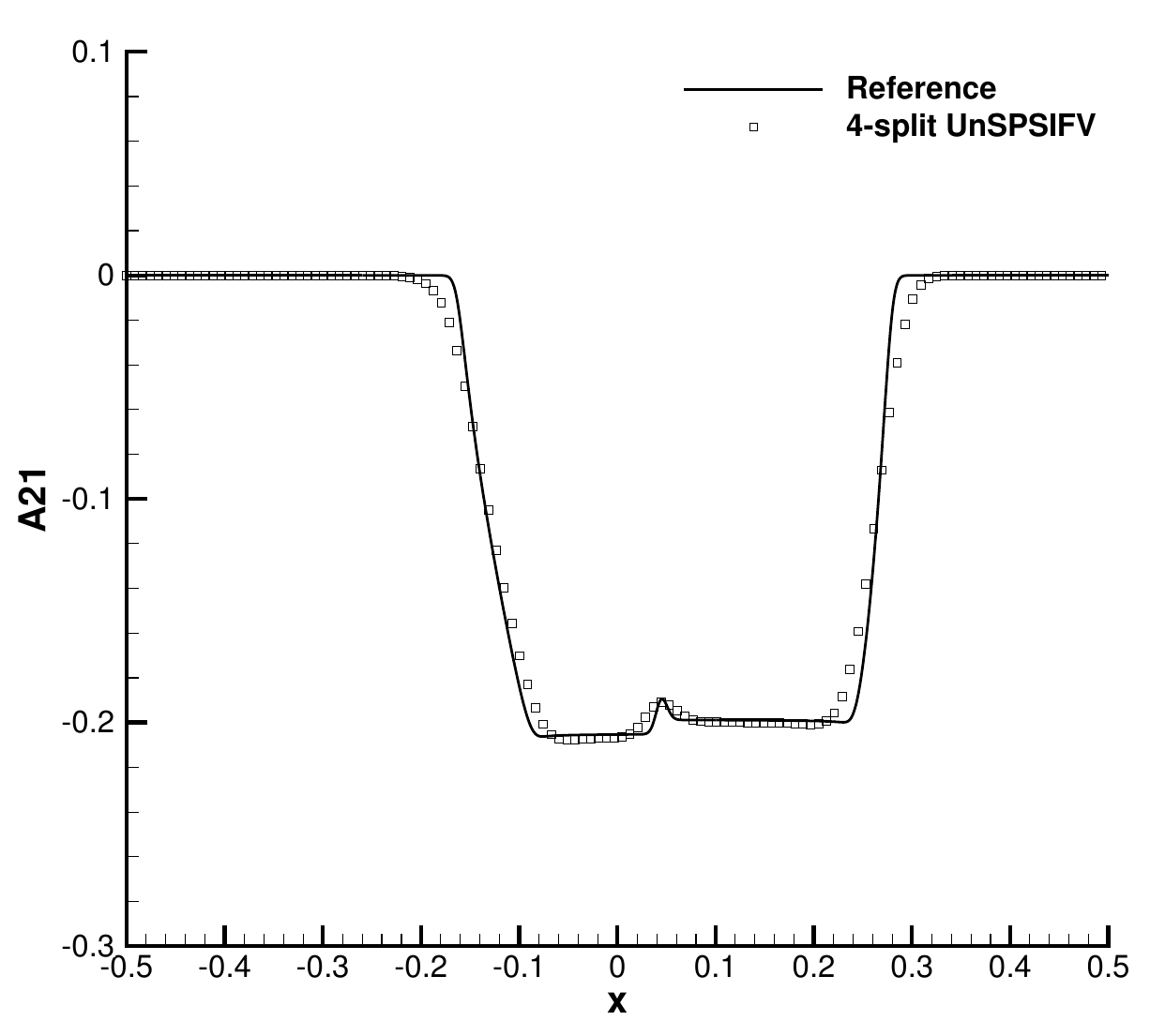}   &
            		\includegraphics[width=0.32\textwidth]{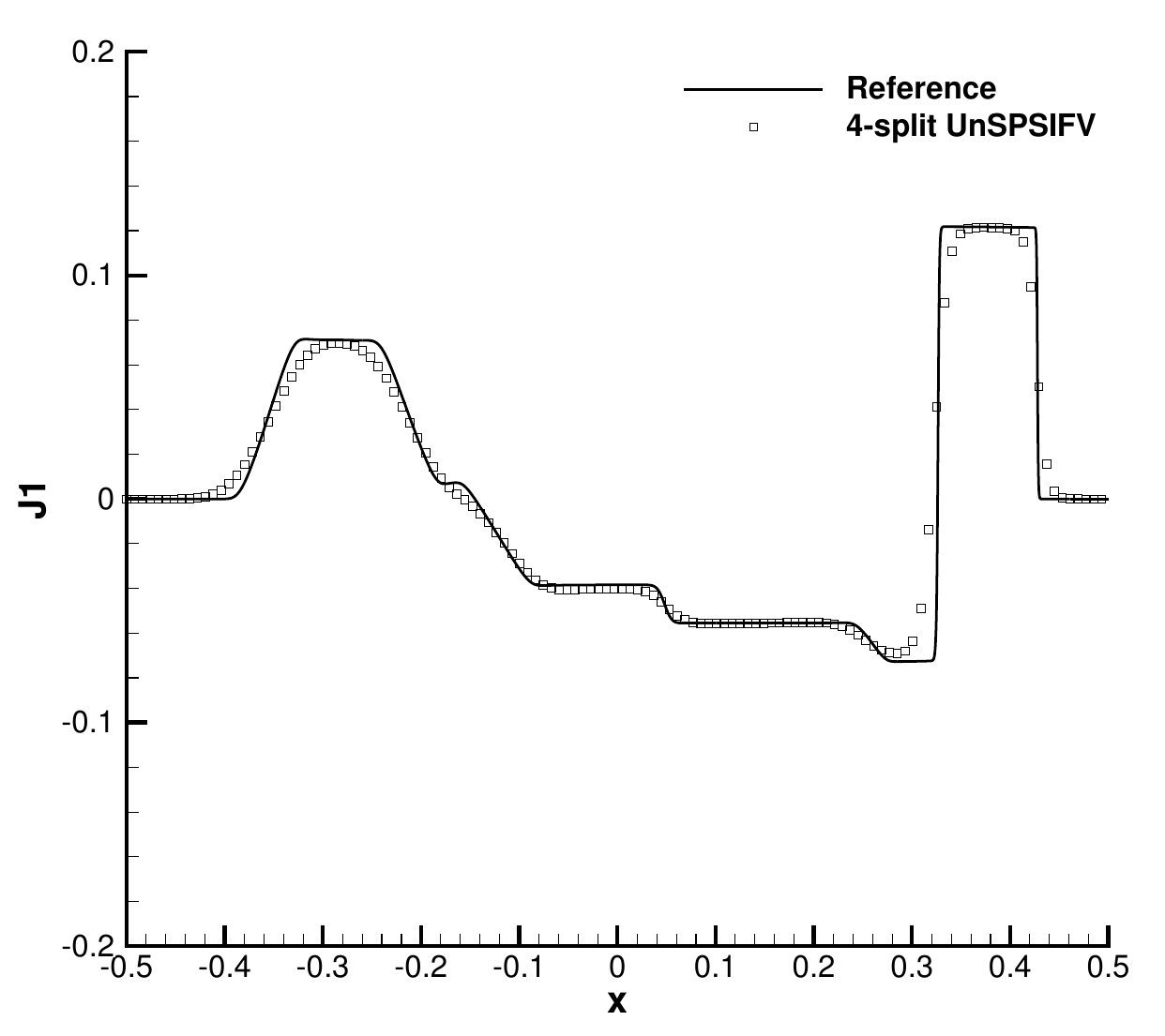}
        \end{tabular}
        \caption{Reference solution and numerical solution of the homogeneous GPR model without source terms ($\tau_1 = \tau_2 = 10^{20}$) for Riemann problem RP3 at a final time of $t=0.2$. Top row: density $\rho$, the velocity component $v$ and the pressure $p$. Bottom row: distortion field components $A_{11}$, $A_{21}$ and thermal impulse component $J_1$.
            One can note seven waves that are contained in the homogeneous part of the GPR model.}
        \label{fig.rp3}
    \end{center}
\end{figure}

\begin{figure}[htbp]
    \begin{center}
        		\includegraphics[width=0.32\textwidth]{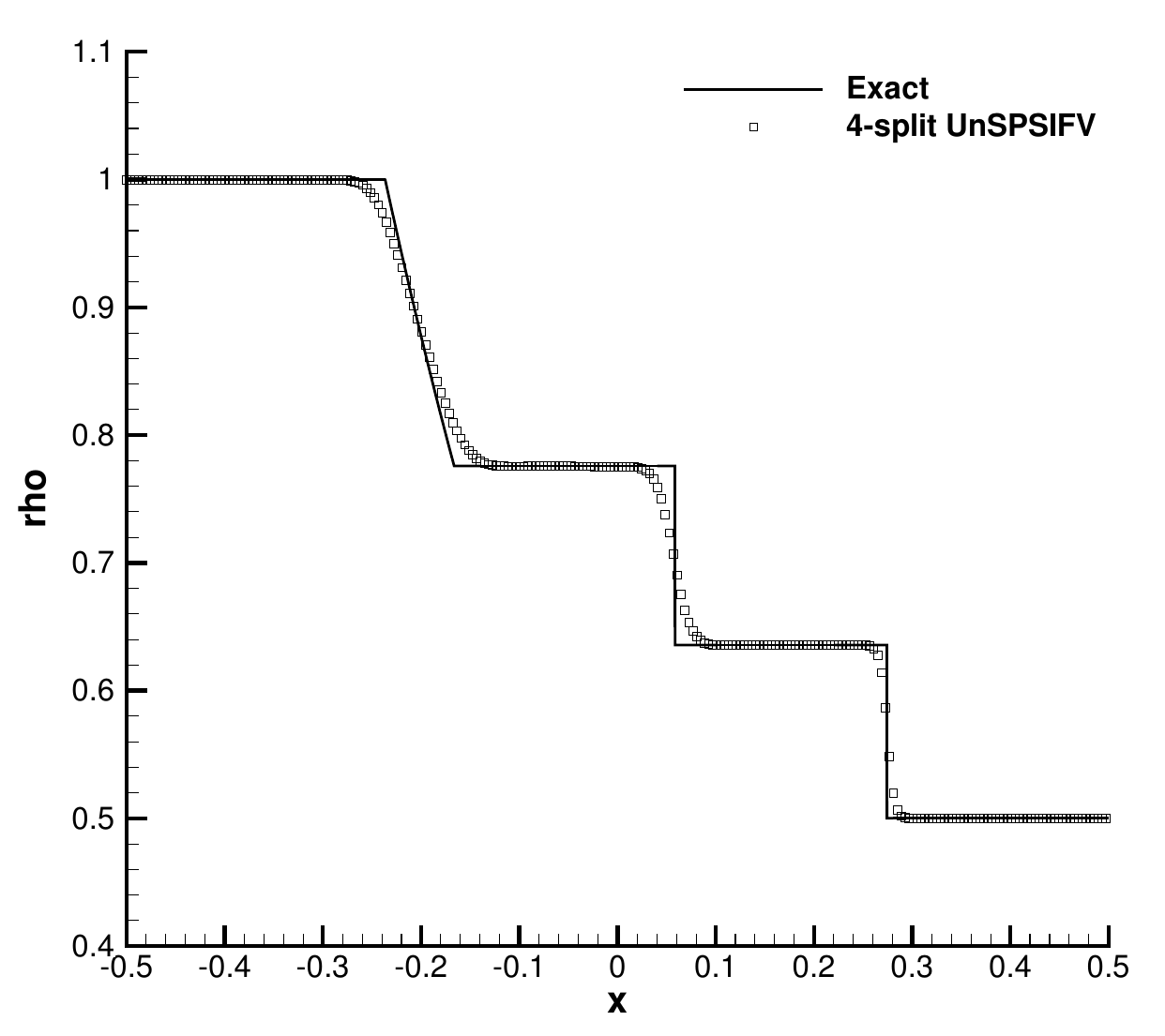}
        		\includegraphics[width=0.32\textwidth]{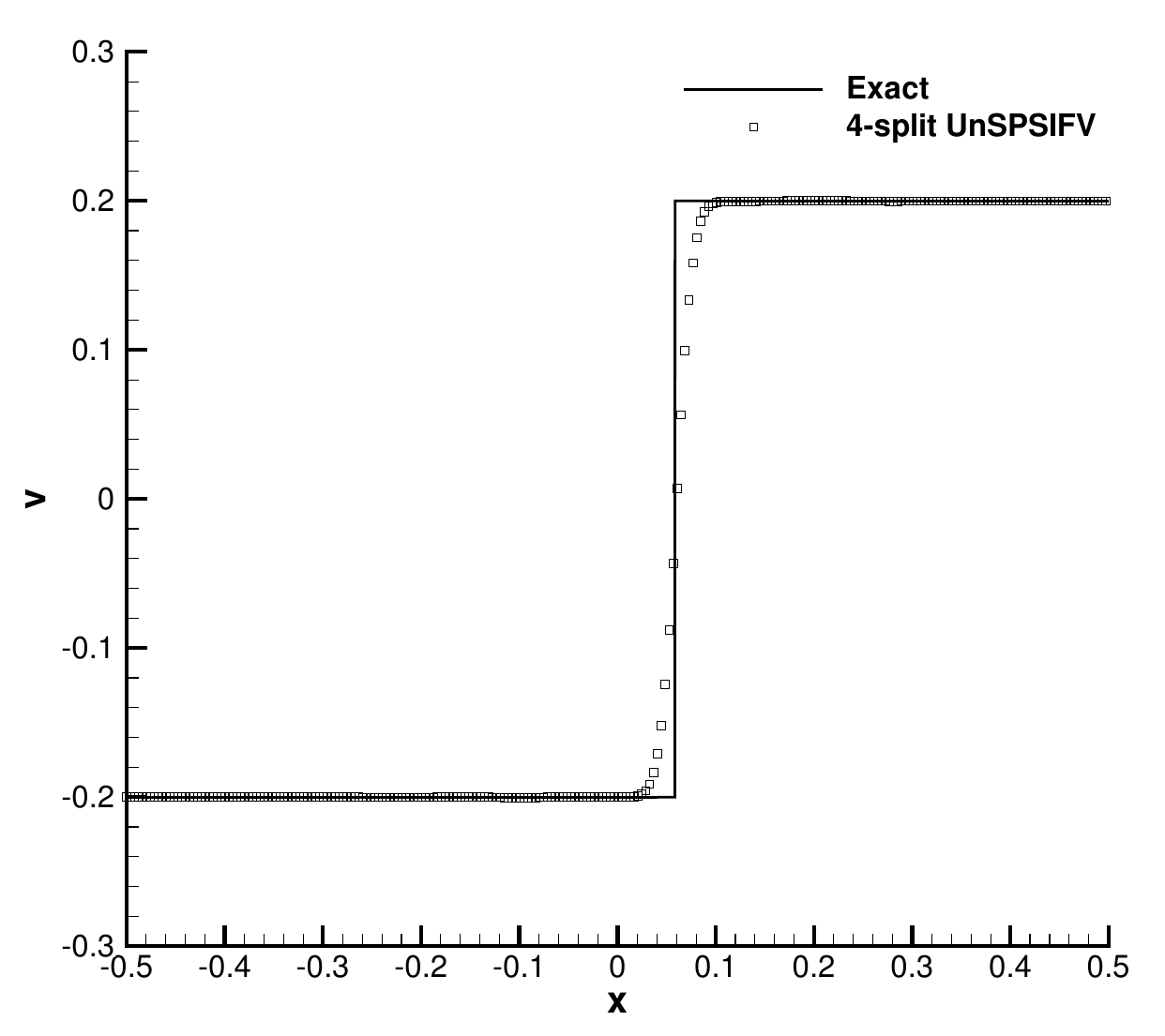}
        		\includegraphics[width=0.32\textwidth]{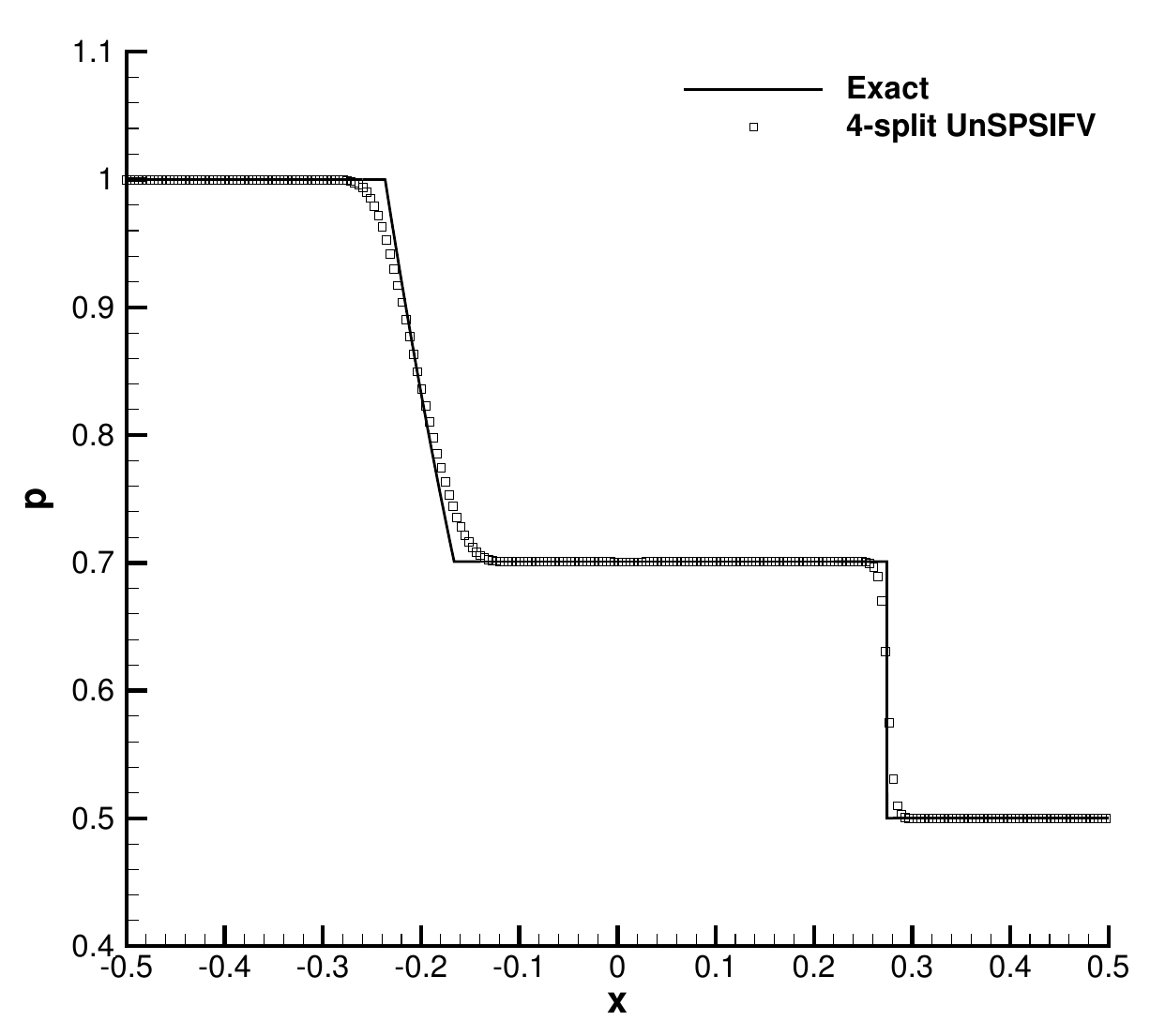}
        \caption{Exact solution of the Euler equations and numerical solution of the GPR model in the stiff relaxation limit ($\tau_1 = 10^{-10}, \tau_2 = 10^{-12}$) for Riemann problem RP4. The density $\rho$, the velocity component $v$ and the pressure $p$ are shown at a final time of $t=0.2$.}
        \label{fig.rp4}
    \end{center}
\end{figure}

Next, we consider a series of Riemann problems (RPs) in 1D regarding the fluid and solid limit of the GPR model.
Therein, RP1 is the classical Sod shock tube, RP2 is the well-known Lax problem, RP3 is a Riemann problem placed in the solid limit and RP4 is RP3 in the fluid limit.
As reference solution in the fluid limit, given by $\tau_1, \tau_2 \to 0$, the exact solution of the Riemann problem of the compressible Euler equations \cite{toro-book} is used while in the solid limit the reference is given by numerical results published in \cite{SIGPR,HTCGPR}.
Due to the non-linear hyperbolic nature of the model, the numerical scheme has to be applicable to flows containing discontinuities such as shock and contact waves typically appearing in high Mach number regimes.
For all RPs, the computational domain $\Omega = [-0.5, 0.5]\times [-0.05,0.05]$ is discretized with $8862$ triangles, except for RP3 where $35998$ elements are used since the solution is composed of a large number of different waves. In this case periodic boundary conditions are imposed only in the $y$-direction while the exact solution is imposed at the left and the right boundaries.
The initial data for density, velocity field and pressure are given in Table \ref{tab.ic.ideal} which includes also the applied relaxation times $\tau_1$ and $\tau_2$. The initial conditions for the $z$-velocity component, distortion field and thermal impulse are for all RPs set to $v_3=0$, $\A = \mathbf{I}$ and $\J = 0$, respectively.
The parameters of the GPR model for RP1 and RP2 are given by $\rho_0 = 1$, $\gamma = 1.4$, $c_s = 100$, $c_h = 10$ and $c_v=1$.
For RP3 and RP4 the parameters are set to $\rho_0 = 1$, $\gamma = 1.4$ and $c_s = 1$, $c_h = 1$, $c_v=1$.
In Figures  \ref{fig.rp1} - \ref{fig.rp4} the numerical results obtained with the new scheme are depicted showing a good agreement between numerical and reference solutions.
\begin{table}[!htbp]
    \renewcommand{\arraystretch}{1.25}
    \caption{Initial states left (L) and right (R) for density $\rho$, velocity $\mathbf{v}$ and pressure $p$
        for a set of Riemann problems solved on the domain $\Omega=[-\frac{1}{2},+\frac{1}{2}] \times [-\frac{1}{20},+\frac{1}{20}]$ using the new four-split scheme.
        The Riemann problems include the fluid limit (RP1, RP2 and RP4) as well as the solid limit (RP3). The parameters $c_h$ and $c_s$ as well as the relaxation times $\tau_1$ and $\tau_2$ are also specified. In all cases we set $\gamma=1.4$. }
    \begin{center}
        \begin{tabular}{ccccccccccccc}
            \hline
            RP & $\rho_L$ & $u_L$ & $v_L$ & $p_L$ & $\rho_R$ & $u_R$ & $v_R$ & $p_R$ & $c_s$ & $c_h$ & $\tau_1$ & $\tau_2$  \\
            \hline
            RP1 &  1.0      &  0.0       &  0.0 & 1.0     & 0.125      &  0.0        &  0.0 & 0.1      & 100 & 10 & $10^{-10}$ & $10^{-12}$ \\
            RP2 &  0.445    &  0.698     &  0.0 & 3.528   & 0.5        &  0.0        &  0.0 & 0.571    & 100 & 10 & $10^{-10}$ & $10^{-12}$ \\
            RP3 &  1.0      &  0.0       & -0.2 & 1.0     & 0.5        &  0.0        & +0.2 & 0.5      & 1 & 1 & $10^{20}$ & $10^{20}$ \\
            RP4 &  1.0      &  0.0       & -0.2 & 1.0     & 0.5        &  0.0        & +0.2 & 0.5      & 1 & 1 & $10^{-10}$ & $10^{-12}$ \\
            \hline
        \end{tabular}
    \end{center}
    \label{tab.ic.ideal}
\end{table}
\subsection{Simple shear layer}
\begin{figure}[htbp]
    \begin{center}
        \begin{tabular}{cc}
            			\includegraphics[width=0.38\textwidth]{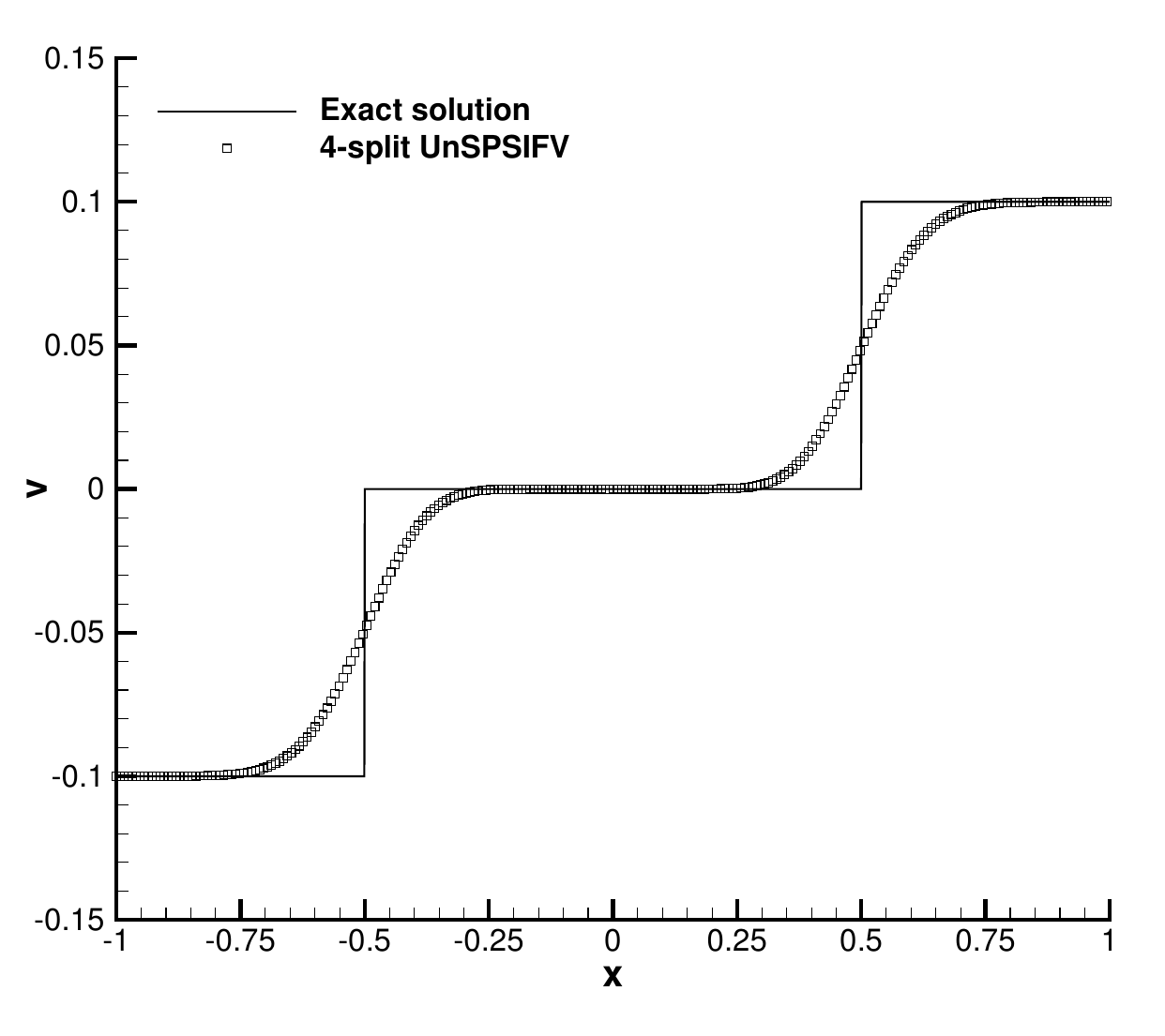}   &
            			\includegraphics[width=0.38\textwidth]{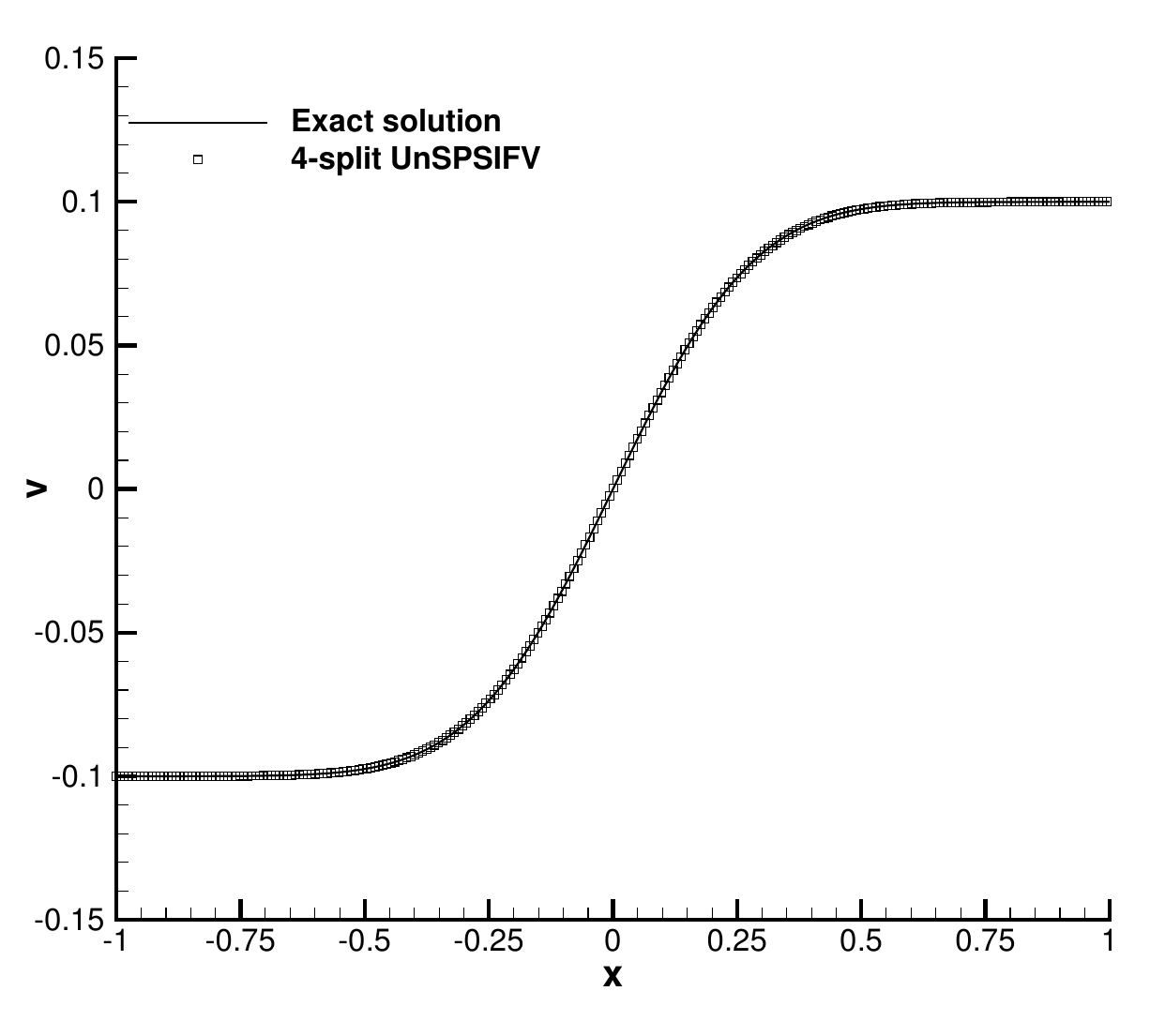}  \\
            			\includegraphics[width=0.38\textwidth]{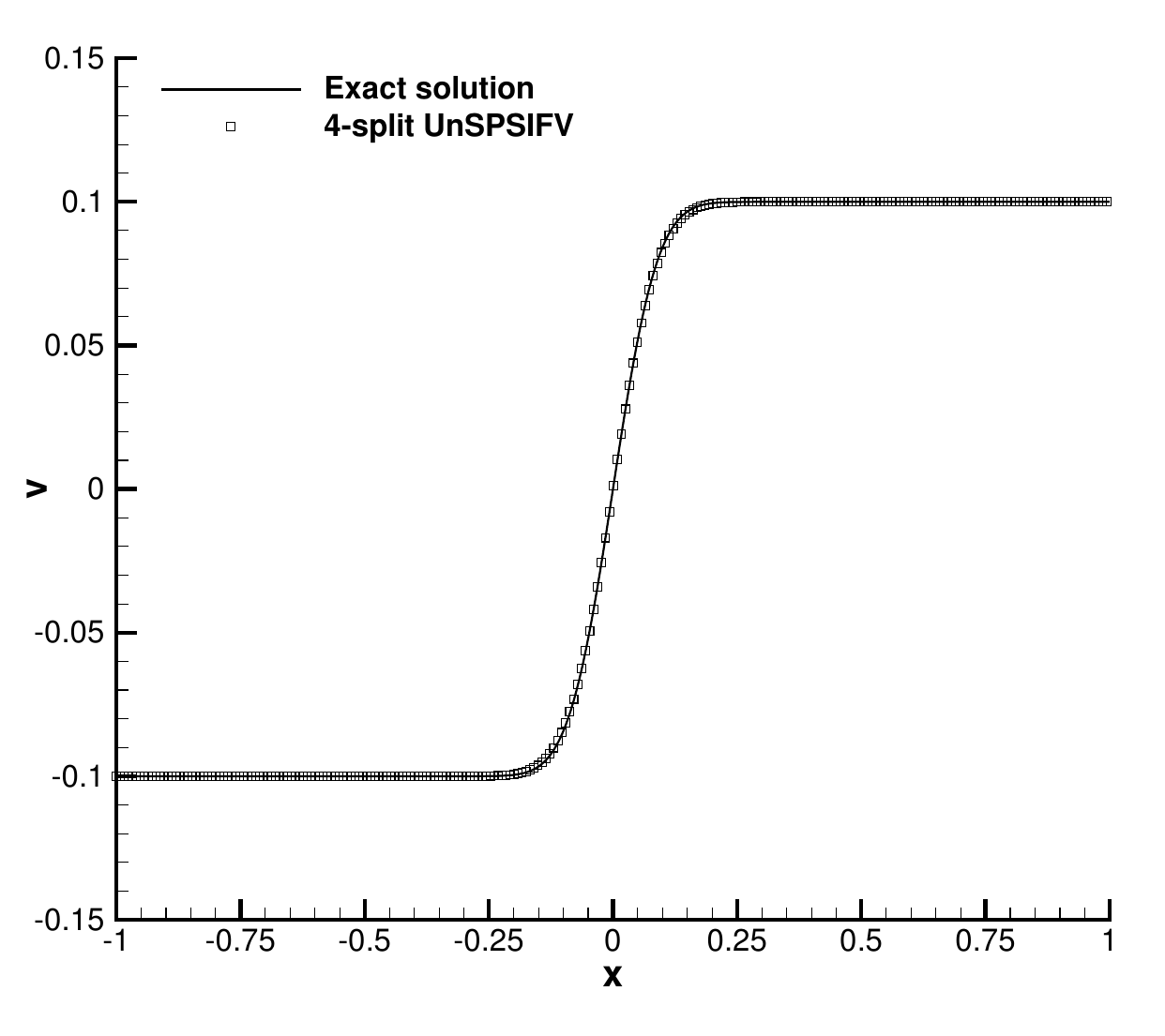}   &
            			\includegraphics[width=0.38\textwidth]{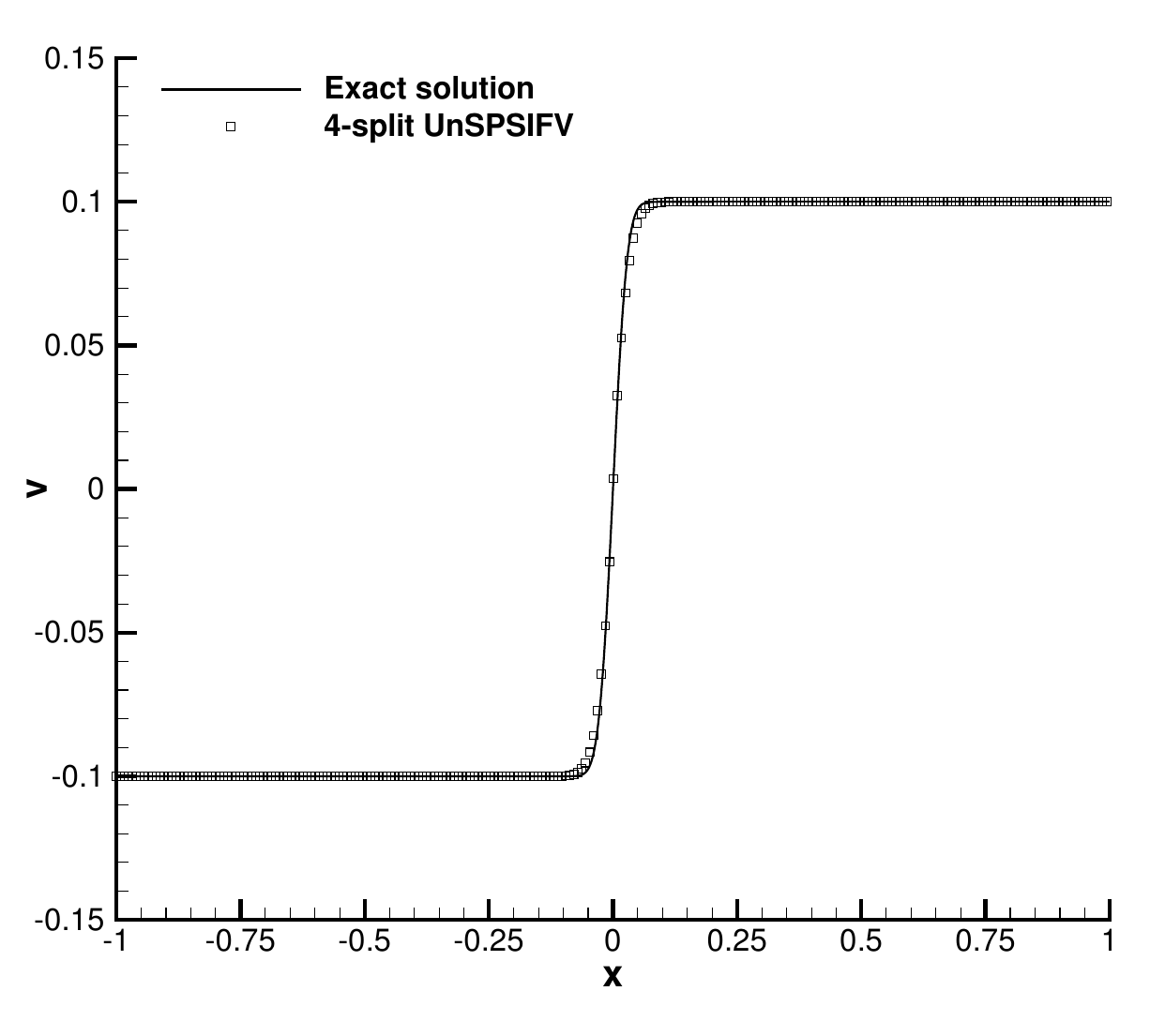}
        \end{tabular}
        \vspace*{-4mm}
        \caption{Numerical solution obtained with the new four-split scheme at the corresponding final times for the GPR model applied to a simple shear flow problem. Results for the solid (top left) and for fluids with different viscosities: $\mu=10^{-1}$ (top right), $\mu=10^{-2}$ (bottom left) and $\mu=10^{-3}$ (bottom right).  }
        \label{fig.shear}
    \end{center}
\end{figure}

The next test cases concern simple shear flow in fluids and solids on a two-dimensional computational domain $\Omega = [-1,+1] \times [-0.25,+0.25]$.
The initial conditions are given by
$\rho=1$, $v_1=v_3=0$, $p=10^5$, $\A=\mathbf{I}$, $\mathbf{J}=0$, while the velocity
component $v_2$ is $v_2 = -v_0$ for $x<0$
and $v_2=+v_0$ for $x\geq 0$, with $v_0=0.1$.
Further parameters are set as $\gamma=1.4$, $c_v=1$, $\rho_0=1$, $c_s=10^3$, $c_h=10^2$ and $\tau_2=10^{-12}$.
The simulation is carried out on a grid composed of $22502$ cells where we impose Neumann boundary conditions in the $x$-direction and periodic boundaries in the $y$-direction up to a final time of $t=0.25$ for the fluid case and $t=5 \cdot 10^{-4}$ for the solid case.
The related acoustic Mach number regime is given by $M_a \approx 2.7 \cdot 10^{-4}$, while the shear Mach number is $M_s = 10^{-4}$.
We compare the numerical results of the new four-split structure-preserving scheme with a reference solution obtained for the fluid limit by the exact solution of the incompressible Navier-Stokes equations for the first problem of Stokes \cite{GPRmodel,SIGPR}.
For the solid limit ($\tau_1 \to \infty$), this initial condition leads to a left-travelling and a right-travelling shear wave, both polarized in the $y$-direction and propagating with speed $c_s$. This makes it possible to obtain the reference solution by solving the corresponding Riemann problem for the linearized equations exactly.
The used time step in the fluid case is given by $\Delta t = 5 \cdot 10^{-3}$ and by $\Delta t = 10^{-5}$ in the solid case.
In Figure \ref{fig.shear} the numerical solution for different values of $\mu$ is displayed showing an excellent agreement with the reference solutions although the Mach numbers are very low.

\subsection{Lid-driven cavity}

\begin{figure}[!htbp]
    \begin{center}
        \begin{tabular}{cc}
         			\includegraphics[trim=10 10 10 10,clip,width=0.47\textwidth]{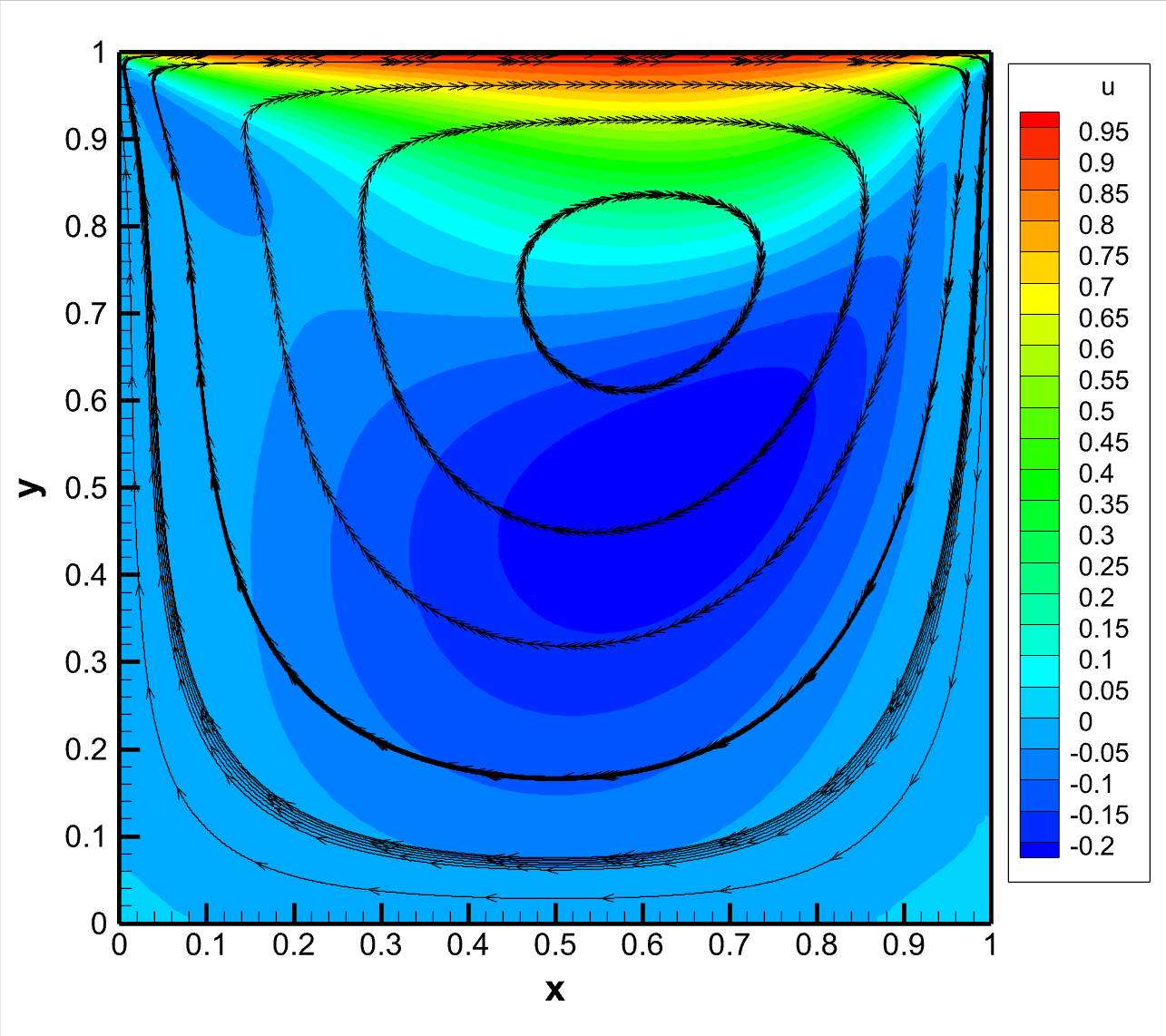}  &
         			\includegraphics[width=0.47\textwidth]{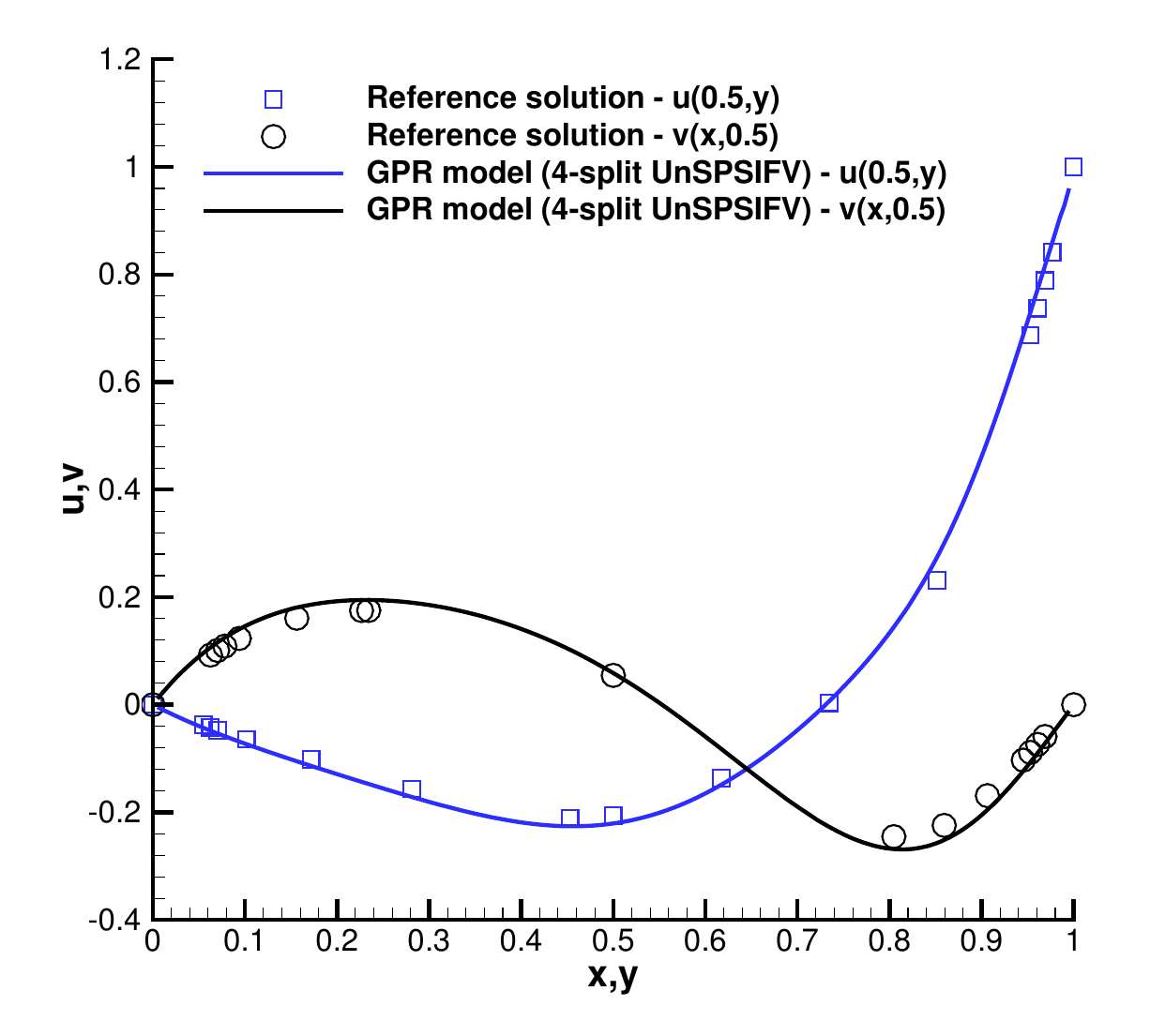}
        \end{tabular}
        \caption{Lid-driven cavity problem at $t=10$ for acoustic Mach number $M_a=2.67 \cdot 10^{-3}$, Reynolds number $Re=100$ and $M_s=10^{-3}$. Color contours of the velocity component $u$ (left) and
            comparison of the velocity components $u$ and $v$ on the centerlines $x=0.5$ and $y=0.5$ with the reference solution of Ghia \textit{et al.} \cite{Ghia1982} (right). }
        \label{fig.cavity}
    \end{center}
\end{figure}

A classical benchmark for the incompressible Navier-Stokes equations is the lid-driven cavity test case  \cite{Ghia1982} which can also be used to validate the performance of compressible flow solvers, such as the new four-split structure-preserving scheme, in the low Mach number regime, see \cite{TavelliDumbser2017,DumbserCasulli2016}.
In the context of the GPR model it has also been successfully used as benchmark problem in \cite{GPRmodel,SIGPR,HTCGPR,HTCAbgrall}.
Here, the problem is formulated on the computational domain $\Omega = [0,1] \times [0,1]$ with initial conditions given by $\rho=1$, $\mathbf{v}=0$, $p=10^5$, $\A=\mathbf{I}$ and $\mathbf{J}=0$.
Furthermore we have $\gamma=1.4$, $c_v = 10^5$, $c_s = 10^3$, $\rho_0=1$, $\tau_2 = 10^{-14}$ and $c_h=0$. To obtain a Reynolds number of $Re=100$ the viscosity coefficient is set to $\mu=10^{-2}$.
Moving the lid at the upper boundary generates the flow inside the cavity; hence the velocity is set there to $\mathbf{v}=(1,0,0)$ while on all other boundaries a no-slip wall boundary condition with $\mathbf{v}=0$ is imposed.
The overall acoustic Mach number is given by $M_a=2.67 \cdot 10^{-3}$ while the shear Mach number is  $M_s = 10^{-3}$ with respect to the lid velocity.
The simulation is performed with the new scheme up to a final time of $t=10$ on a computational grid made of $30906$ cells.
The numerical results are compared with the reference solution of  Ghia \textit{et al.} \cite{Ghia1982} obtained for the incompressible Navier-Stokes equations.
From Figure \ref{fig.cavity} good agreement between the numerical solution obtained with the new scheme for the GPR model and the incompressible Navier-Stokes reference solution can be observed.
We emphasize that compared to the GPR solvers \cite{GPRmodel} and \cite{SIGPR} the time step restriction of the new four-split structure-preserving scheme is only with respect to the fluid velocity and is thus independent of the adiabatic sound speed, shear sound speed $c_s$ and heat wave speed $c_T$.

\subsection{Solid rotor}

In order to numerically verify the discrete curl-free property of the new scheme we consider the solid rotor problem following \cite{SIGPR} setting $\tau_1  = 10^{20}$, i.e. no strain relaxation is present.
Further, the heat conduction is considered in the stiff Fourier limit by setting $\tau_2=10^{-14}$ in model \eqref{eqn.GPR}.
On the computational domain $\Omega = [-1,+1]^2$, the initial condition is given uniformly by $\rho = 1$, $p = 10^5$, $\mathbf{A} = \mathbf{I}$ and $\mathbf{J}=0$ while the initial velocity field is $u = -y/R$, $v = +x/R$ and $w=0$ within the circular region $r \leq R$, where $r = \left\| \mathbf{x} \right\|_2$ and $R=0.2$, while $\mathbf{v}=0$ for $r > R$. The remaining parameters are $\gamma = 1.4$, $c_v=1004/\gamma$, $c_s = 1.0$ and $c_h = 100$.
Hence the flow is placed in the low acoustic and low heat Mach number regime.
The simulation with the new scheme is carried out on a mesh made of $89844$ cells up to the final time $t=0.3$.
The numerical reference solution is obtained by solving the same problem in the absence of heat conduction with the same computational parameters employing the second-order two-split structure-preserving scheme introduced in \cite{SIGPR}.
The reason for setting $c_h=0$ when computing the reference solution is that the two-split scheme is not designed for low heat Mach numbers $M_h$.

\begin{figure}[!htbp]
    \begin{center}
        \begin{tabular}{ccc}
            			\includegraphics[trim=10 10 10 10,clip,width=0.3\textwidth]{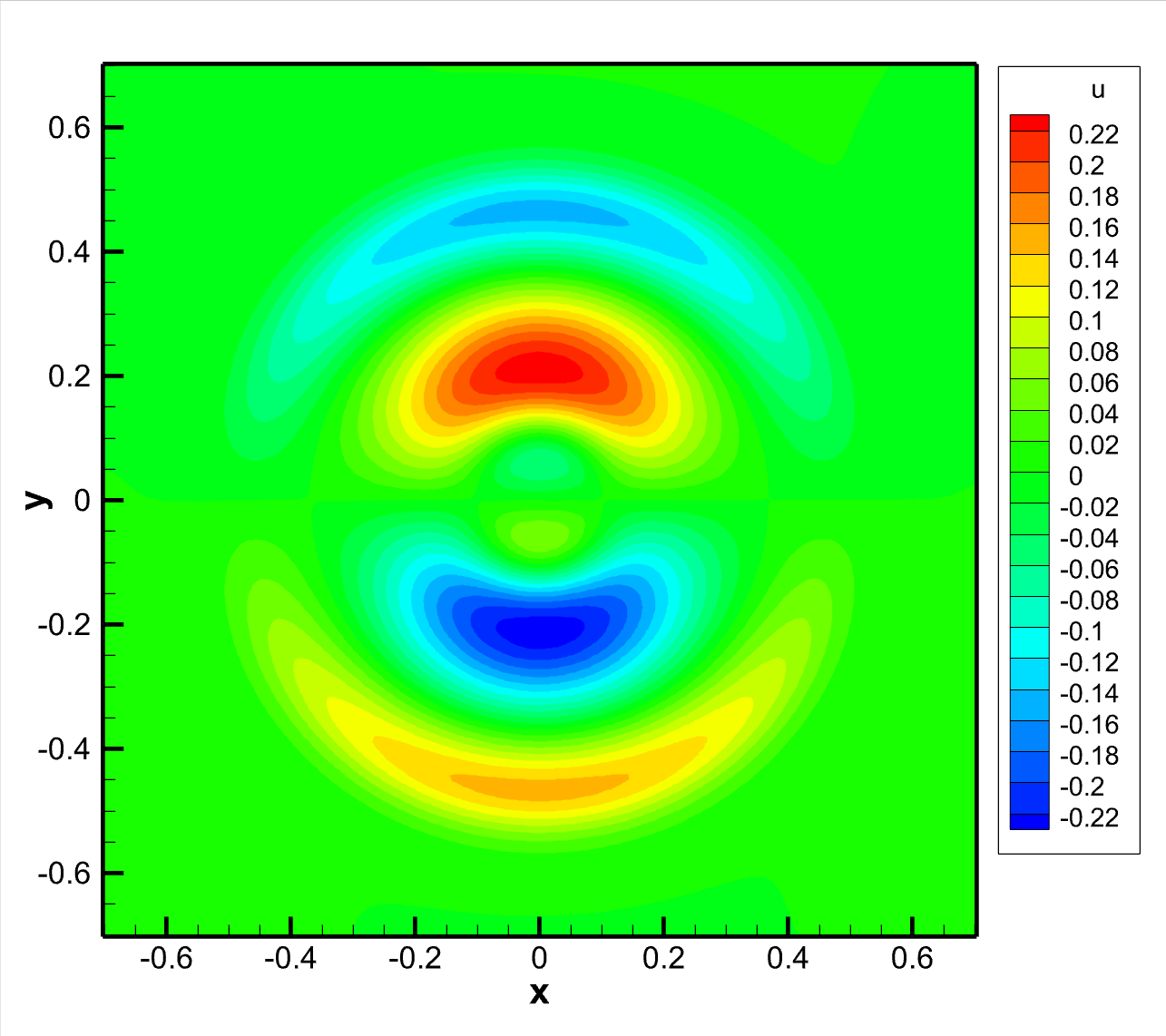}  &
            			\includegraphics[trim=10 10 10 10,clip,width=0.3\textwidth]{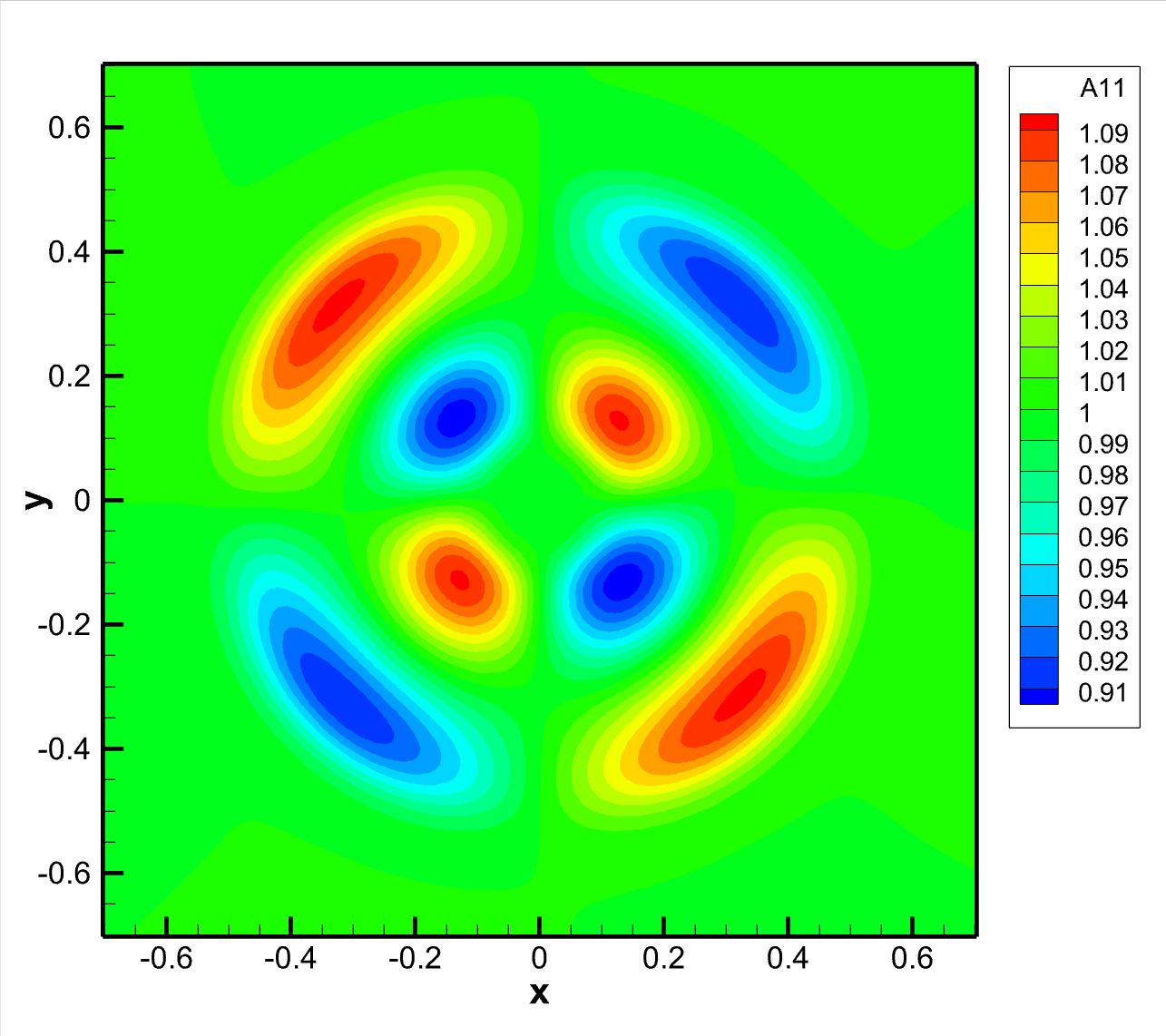}  &
            			\includegraphics[trim=10 10 10 10,clip,width=0.3\textwidth]{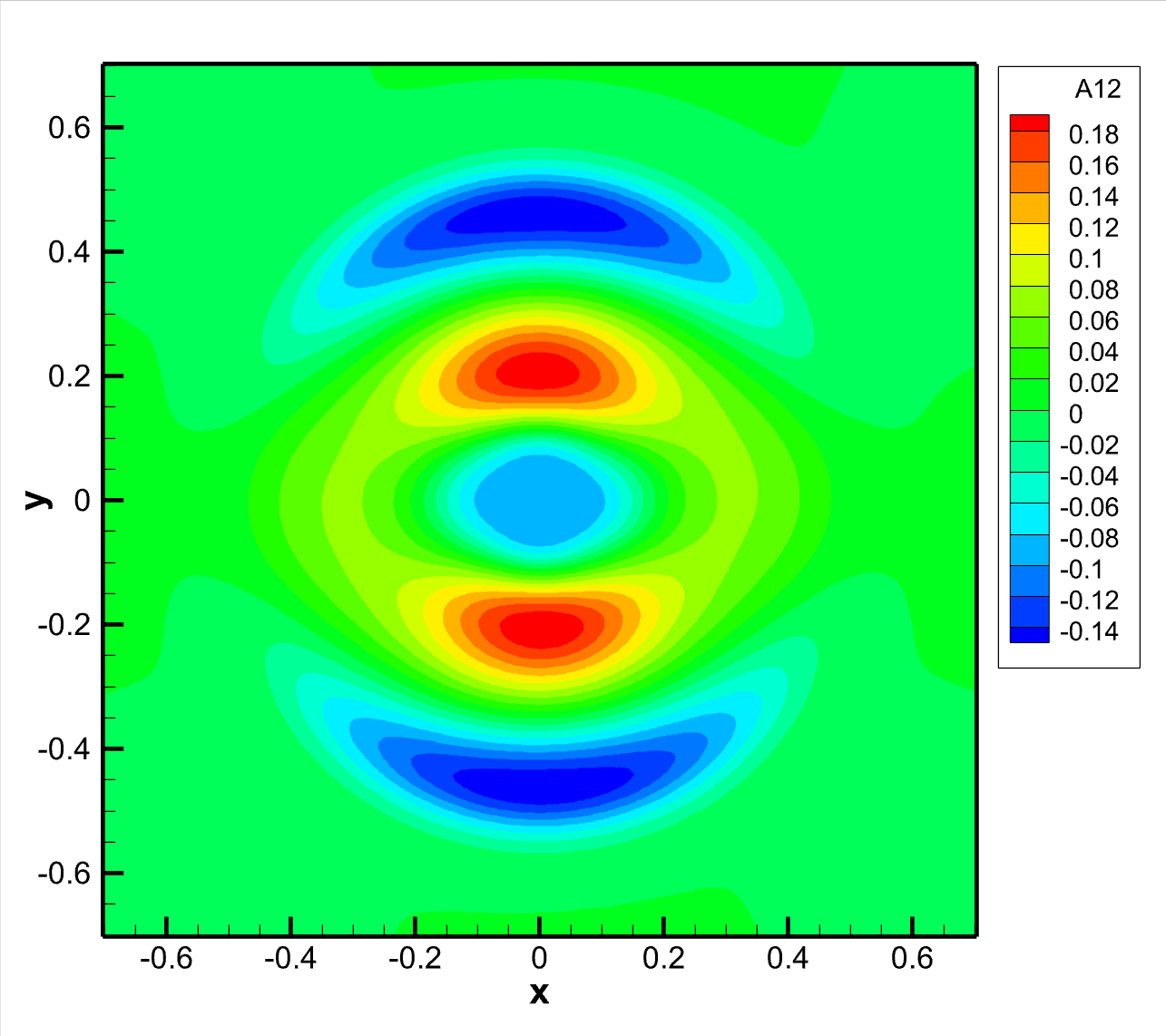}
        \end{tabular}
        \caption{Solid rotor problem at time $t=0.3$. Contour colors of the horizontal velocity $u$ (left), component $A_{11}$ (center) and component $A_{12}$ (right) of the distortion field. }
        \label{fig.solidrotor}
    \end{center}
\end{figure}

\begin{figure}[!htbp]
    \begin{center}
        \begin{tabular}{cc}
            			\includegraphics[width=0.45\textwidth]{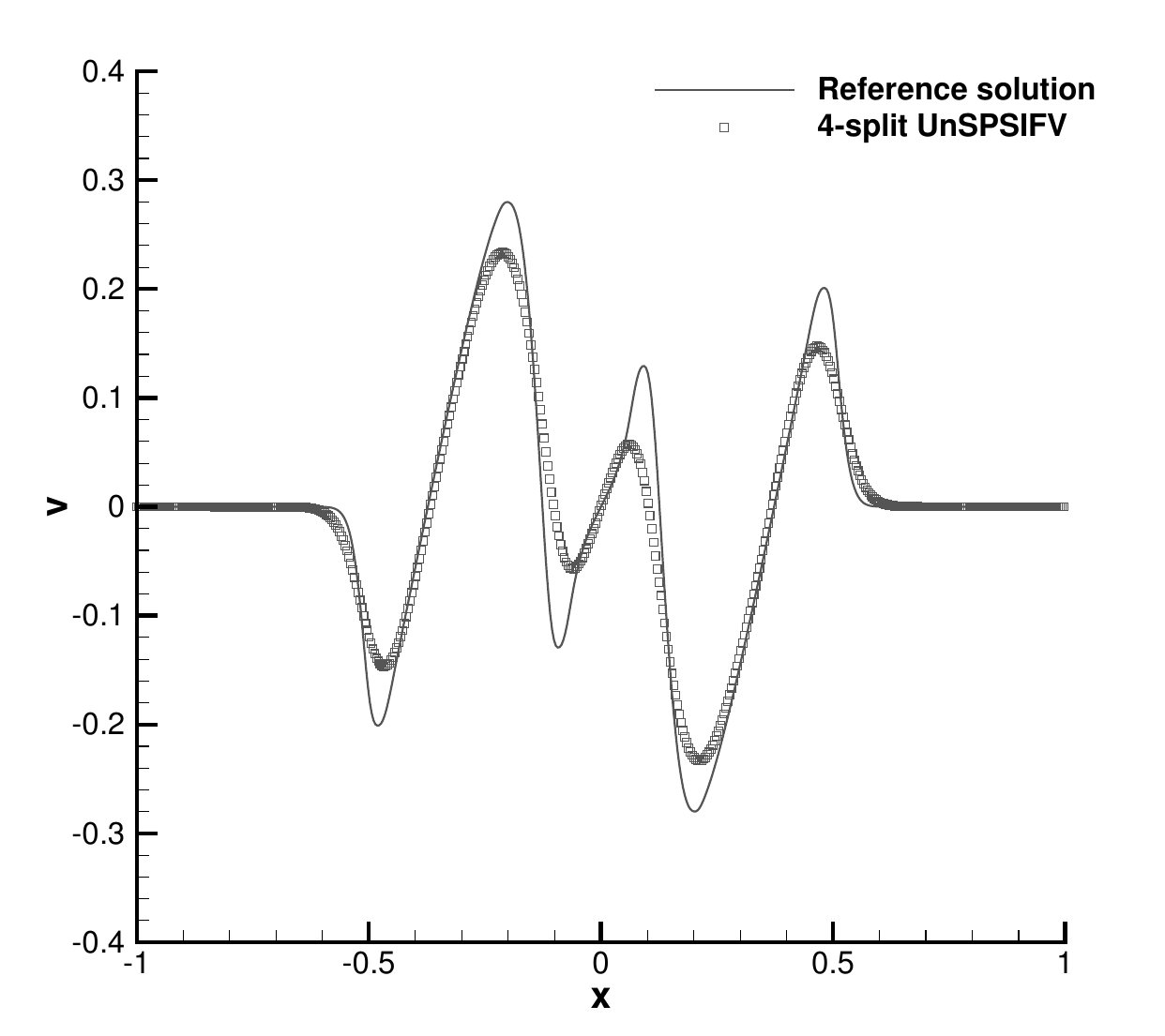}  &
            			\includegraphics[width=0.45\textwidth]{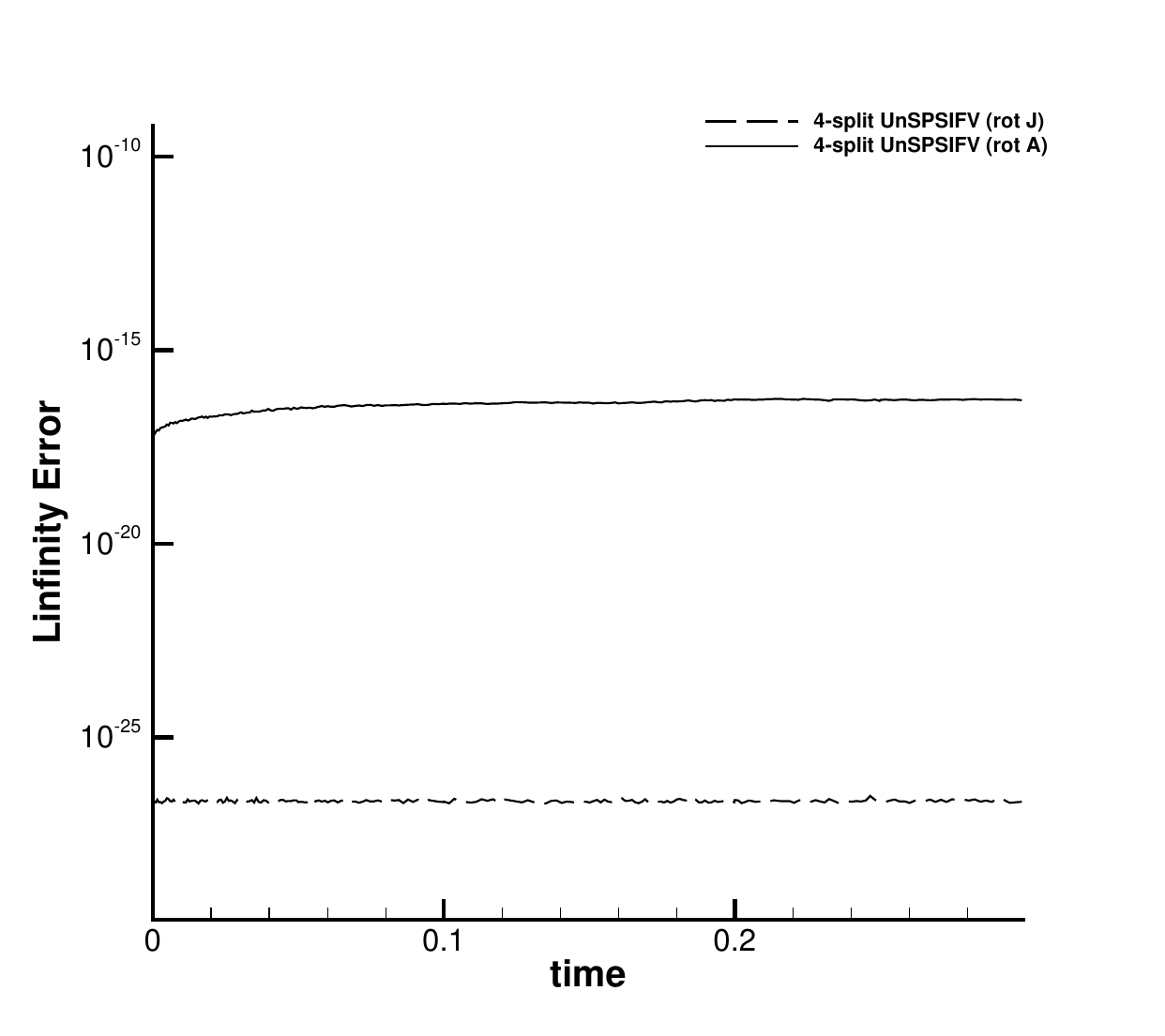}
        \end{tabular}
        \caption{Solid rotor problem. 1D cut through the solution for velocity component $v$ at $y=0$ at time $t=0.3$ and comparison with the numerical reference solution (left). Time series of the curl errors of $\mathbf{A}$ and $\mathbf{J}$ in $L^{\infty}$ norm (right). }
        \label{fig.solidrotor.cutcurl}
    \end{center}
\end{figure}

In Figure \ref{fig.solidrotor} the contour colors of the velocity component $u$ and the components $A_{11}$ and $A_{12}$ of the distortion field obtained with the new scheme are shown.
The comparison with the reference solution is depicted in the left panel of Figure \ref{fig.solidrotor.cutcurl} in terms of a 1D cut at $y=0$.
In the right panel of Figure \ref{fig.solidrotor.cutcurl} the development of the curl errors of $\mathbf{A}$ and $\mathbf{J}$ in the $L^\infty$ norm is given.
We observe a good agreement with the reference solution and, as expected, the curl errors are of the order of machine precision.

\subsection{Circular explosion problem}

\begin{figure}[!htbp]
    \begin{center}
        \begin{tabular}{cc}
            			\includegraphics[trim= 10 10 10 10,clip,width=0.47\textwidth]{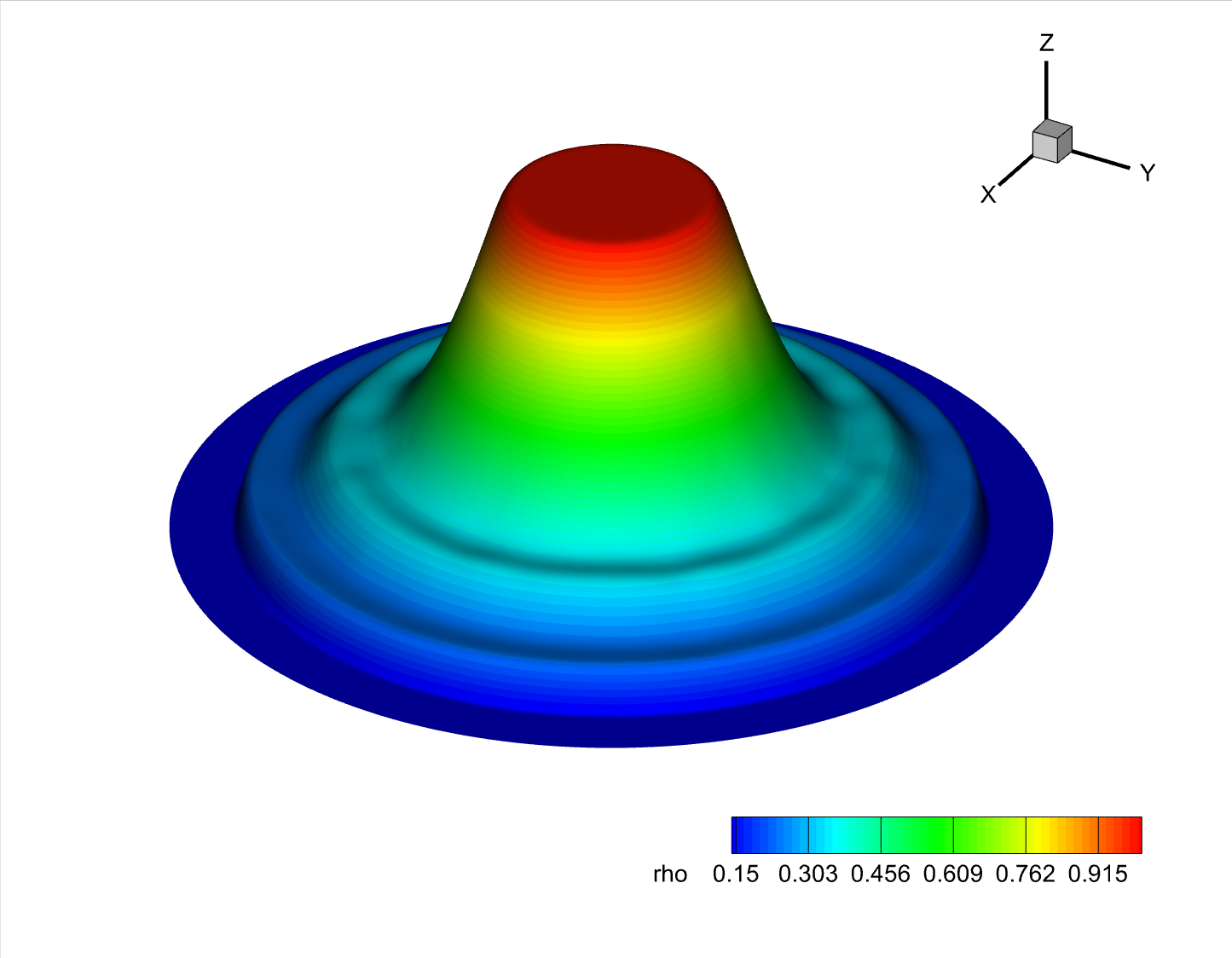}  &
            			\includegraphics[trim= 10 10 10 10,clip,width=0.47\textwidth]{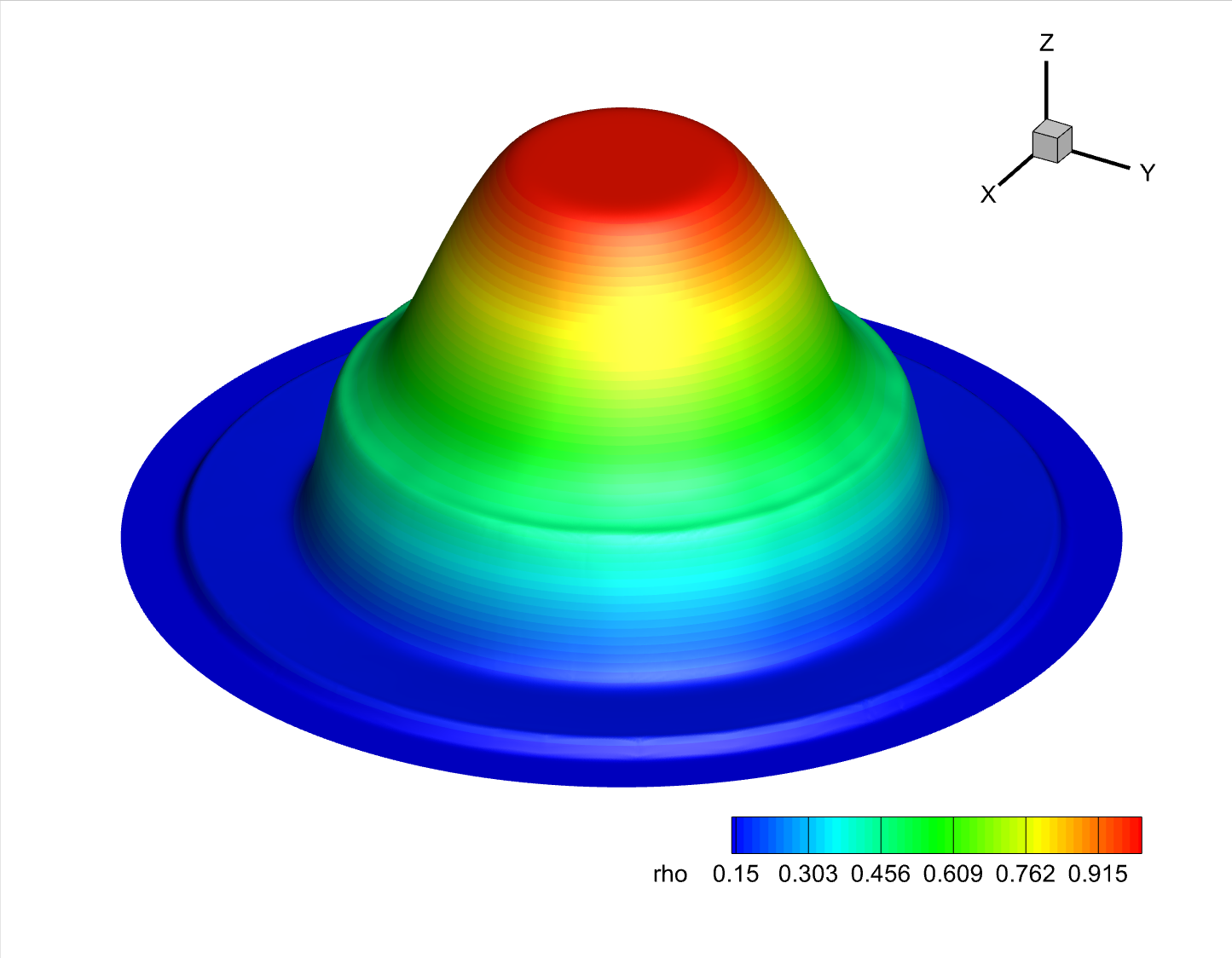}
        \end{tabular}
        \caption{2D explosion problems EP1 and EP2 at their final times $t=0.2$ and $t=0.15$, respectively. 3D density contour color plot for EP1 (left) and EP2 (right).  }
        \label{fig.ep2dfluidsolid}
    \end{center}
\end{figure}

\begin{figure}[!htbp]
    \begin{center}
        \begin{tabular}{ccc}
            			\includegraphics[width=0.3\textwidth]{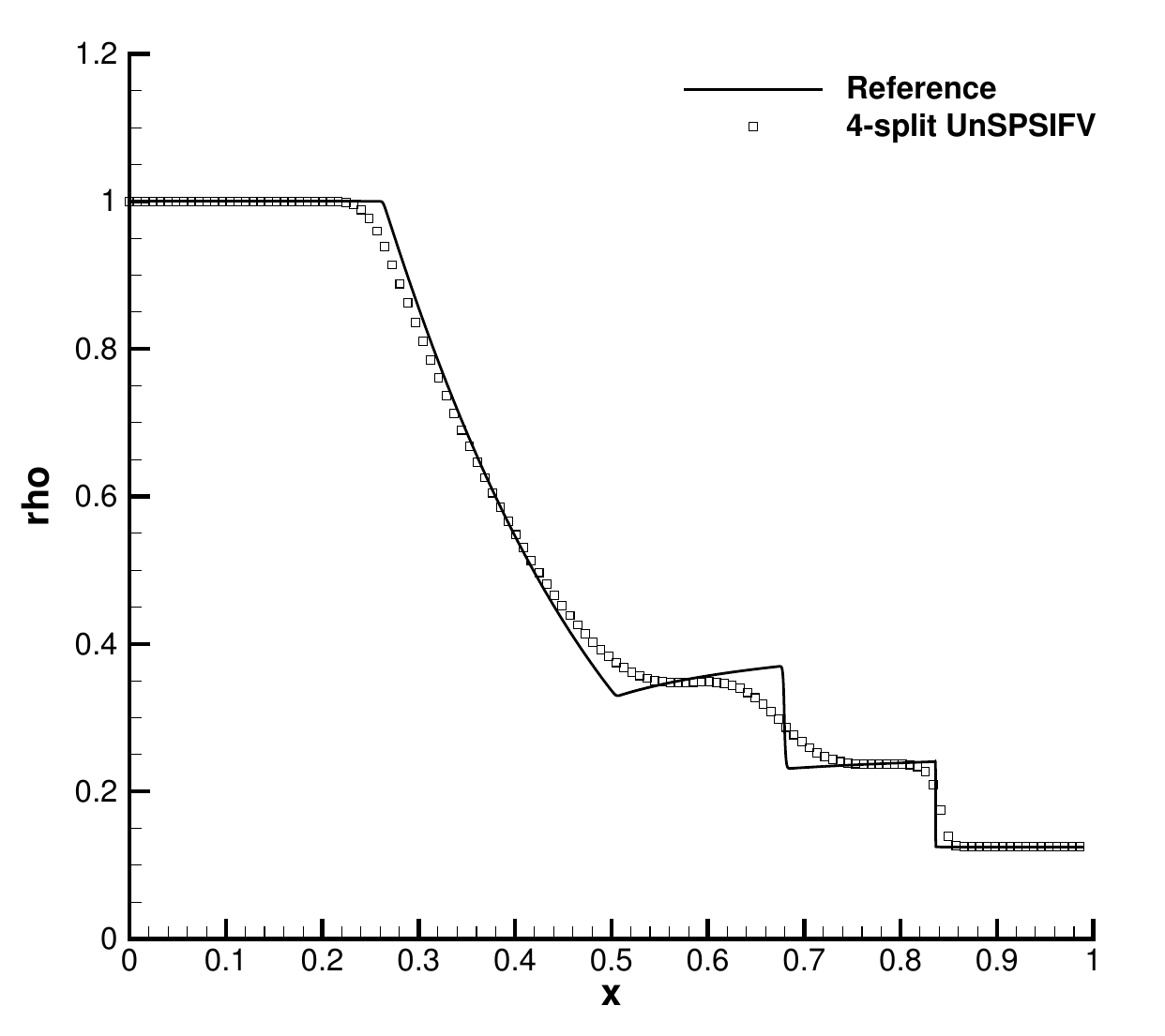}  &
            			\includegraphics[width=0.3\textwidth]{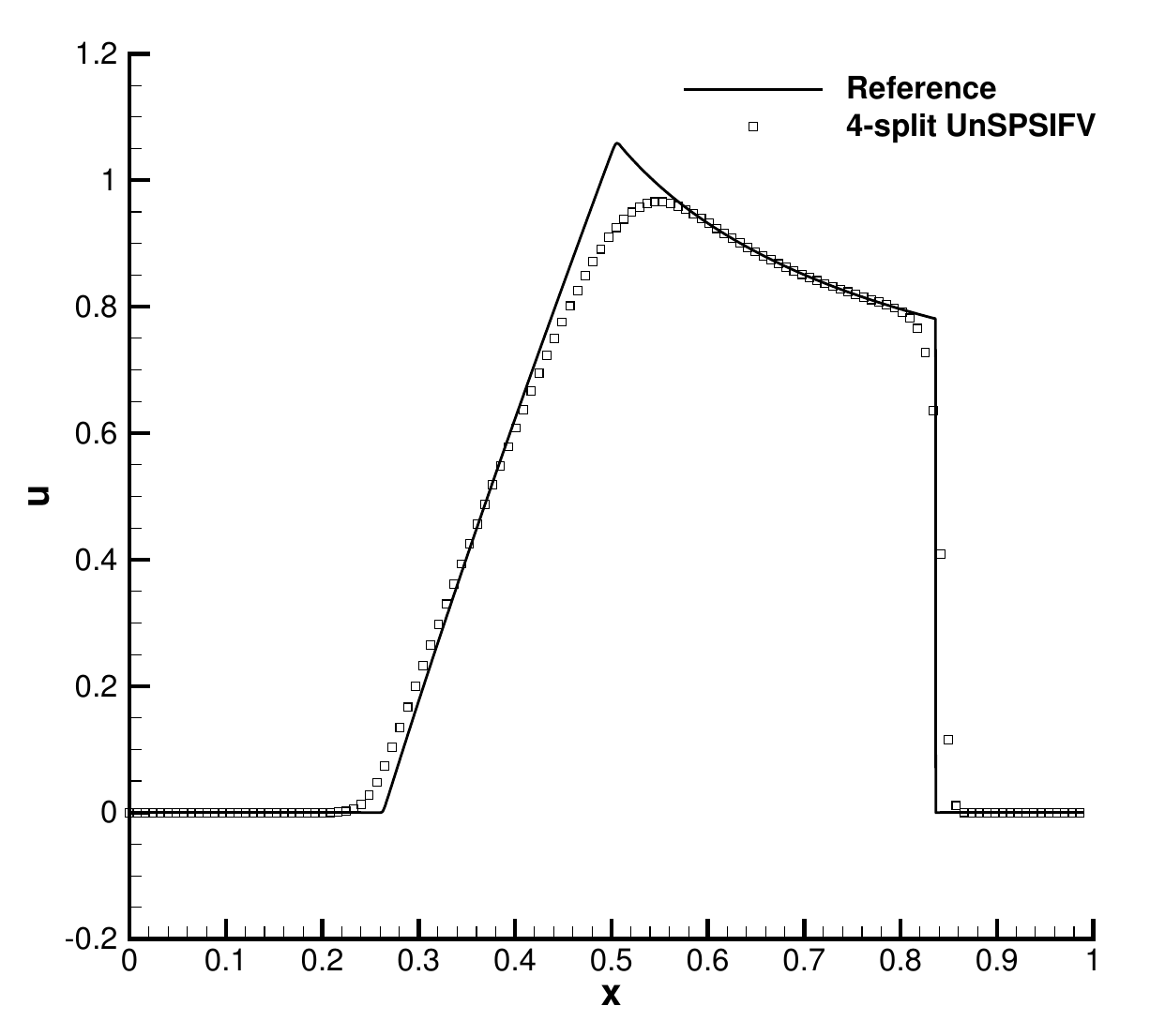}  &
            			\includegraphics[width=0.3\textwidth]{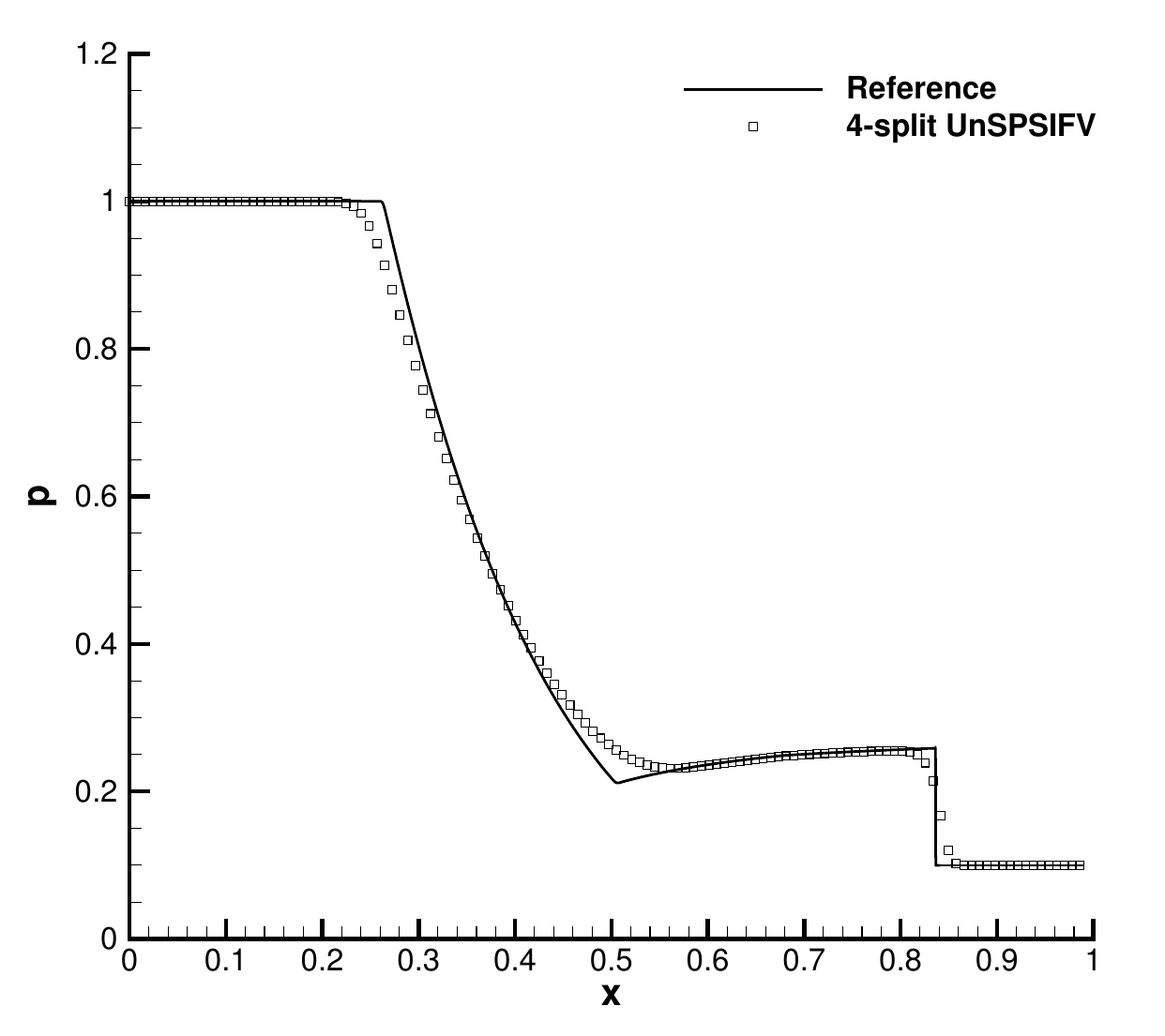}  \\
        \end{tabular}
        \caption{2D explosion problem EP1 at time $t=0.2$. 1D cuts along the $x$-axis providing a comparison of the Euler reference solution (solid line) against the numerical solution of the GPR model obtained
            with the new four-split scheme (square symbols) in the stiff relaxation limit (fluid limit, $\tau_1=10^{-8}$, $\tau_2=10^{-10}$) for density (left), radial velocity (center) and pressure (right).   }
        \label{fig.ep2dfluid}
    \end{center}
\end{figure}

\begin{figure}[!htbp]
    \begin{center}
        \begin{tabular}{cc}
            			\includegraphics[width=0.45\textwidth]{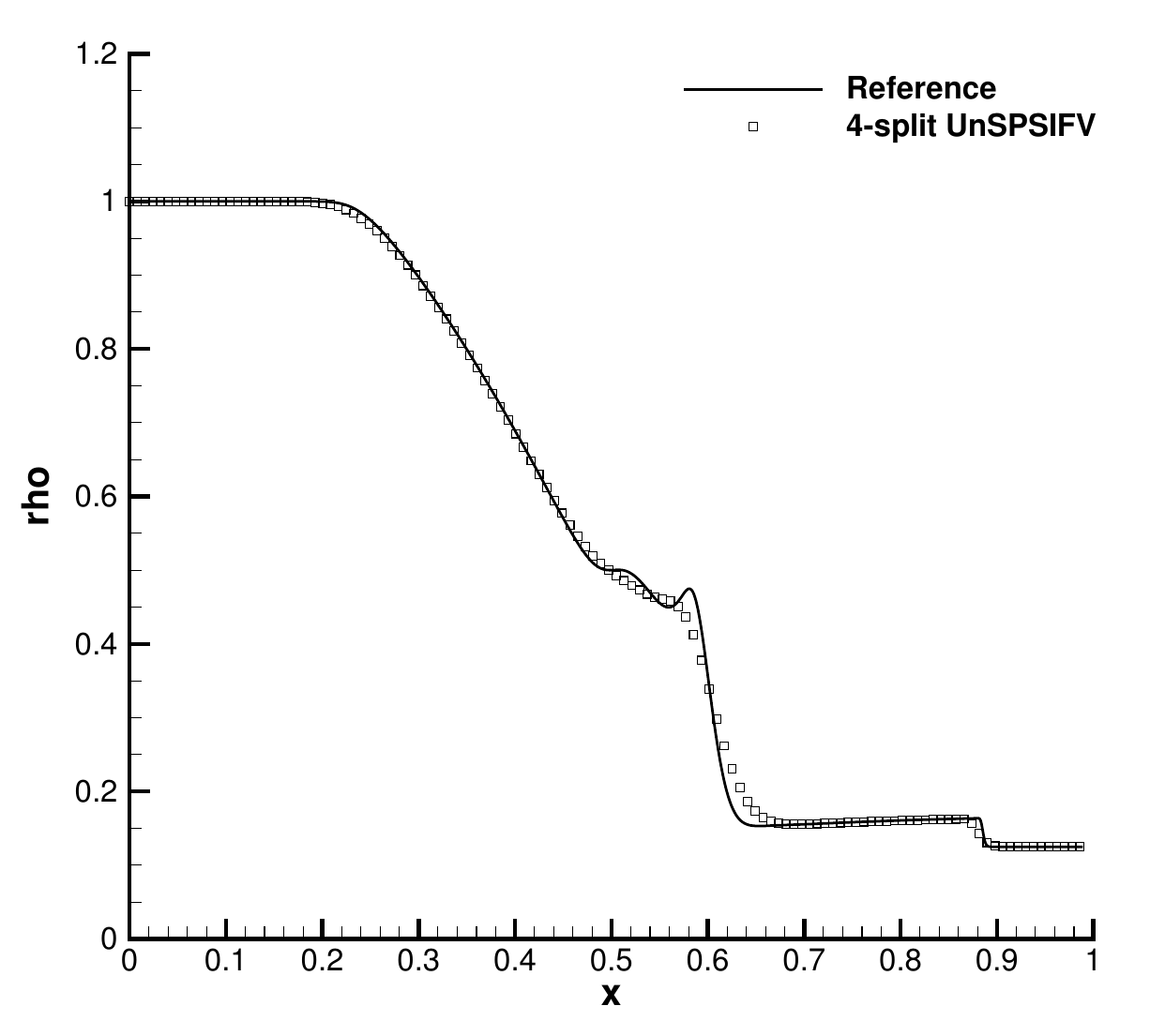}  &
            			\includegraphics[width=0.45\textwidth]{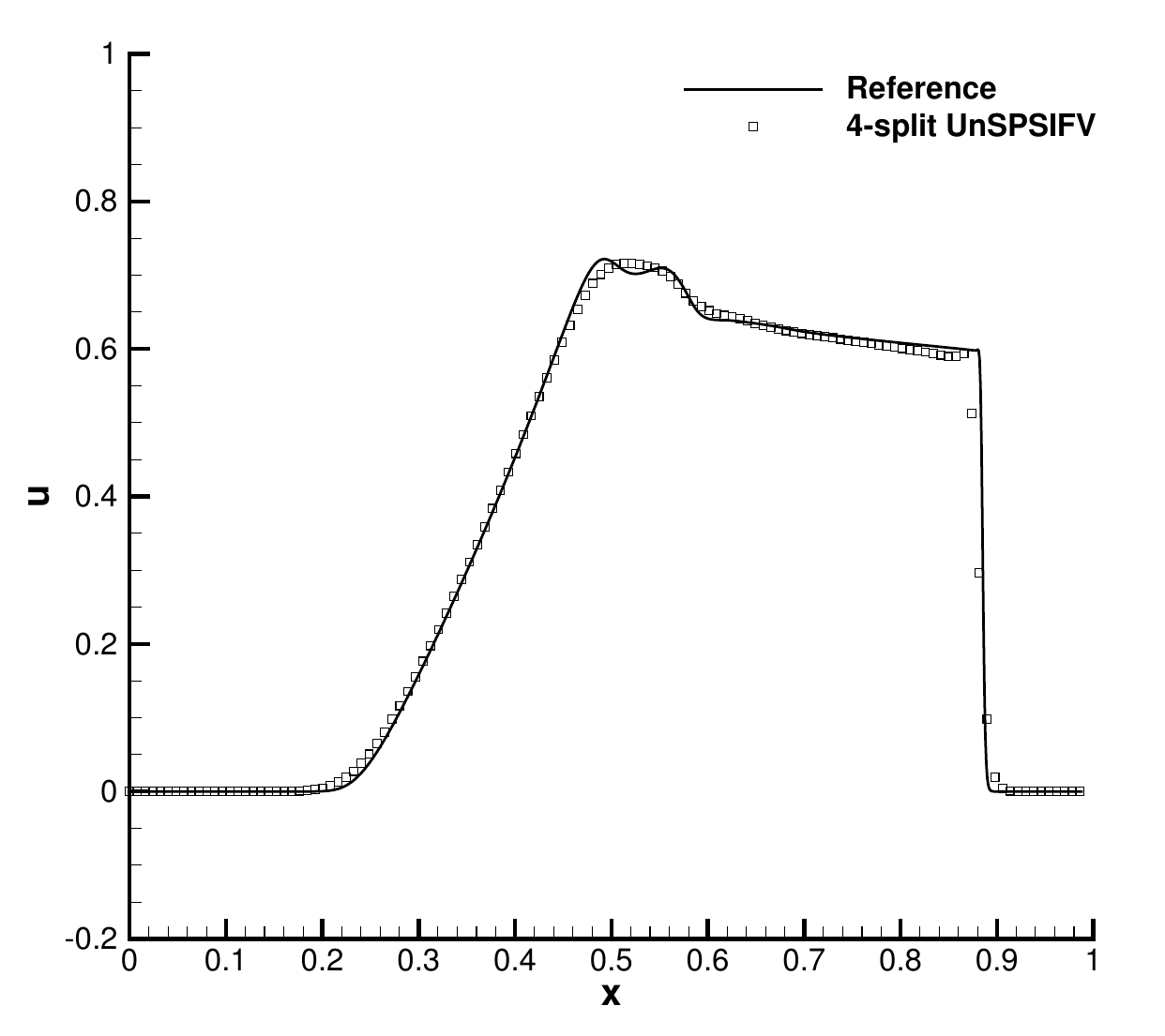}   \\
            			\includegraphics[width=0.45\textwidth]{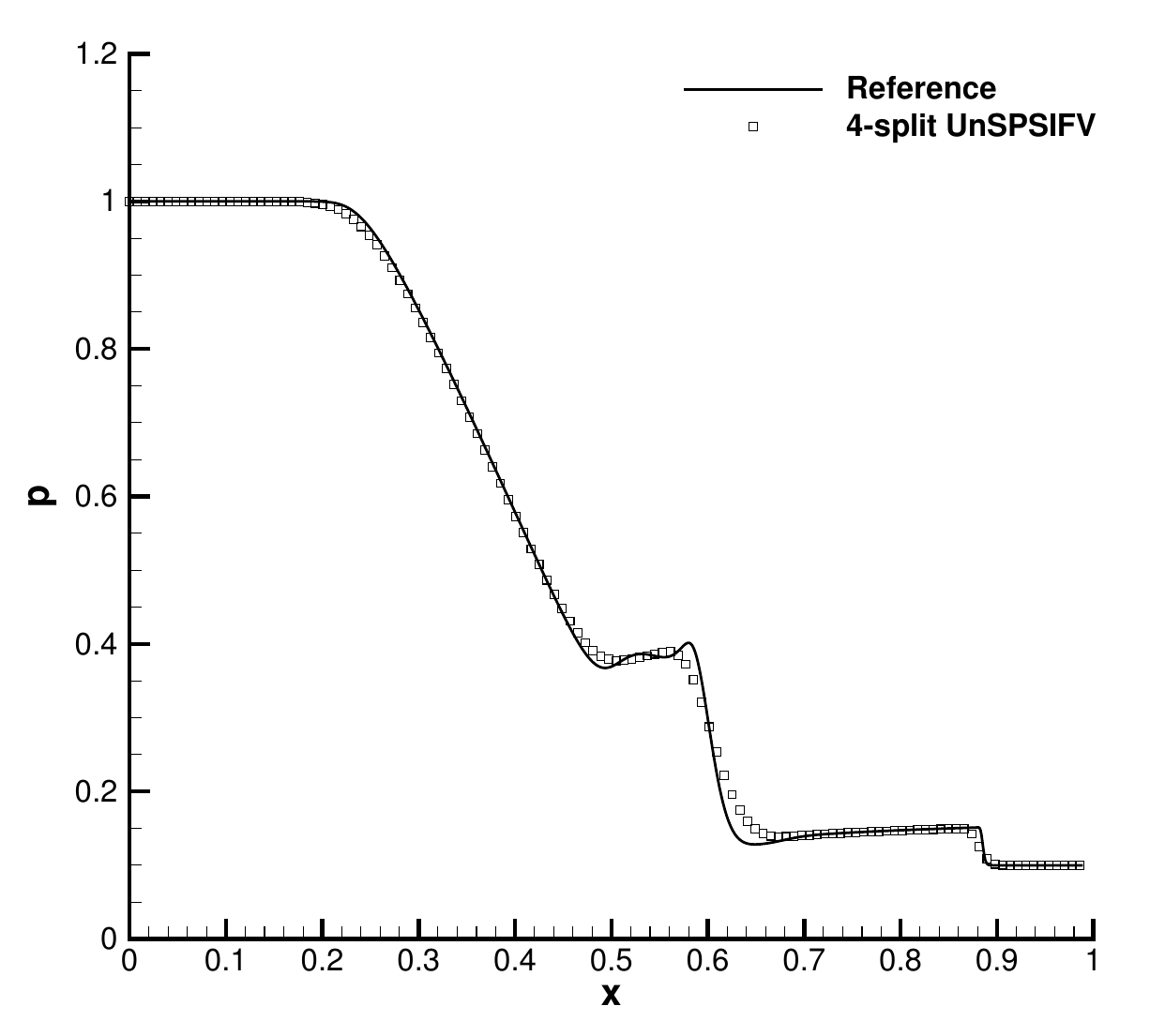}  &
            			\includegraphics[width=0.45\textwidth]{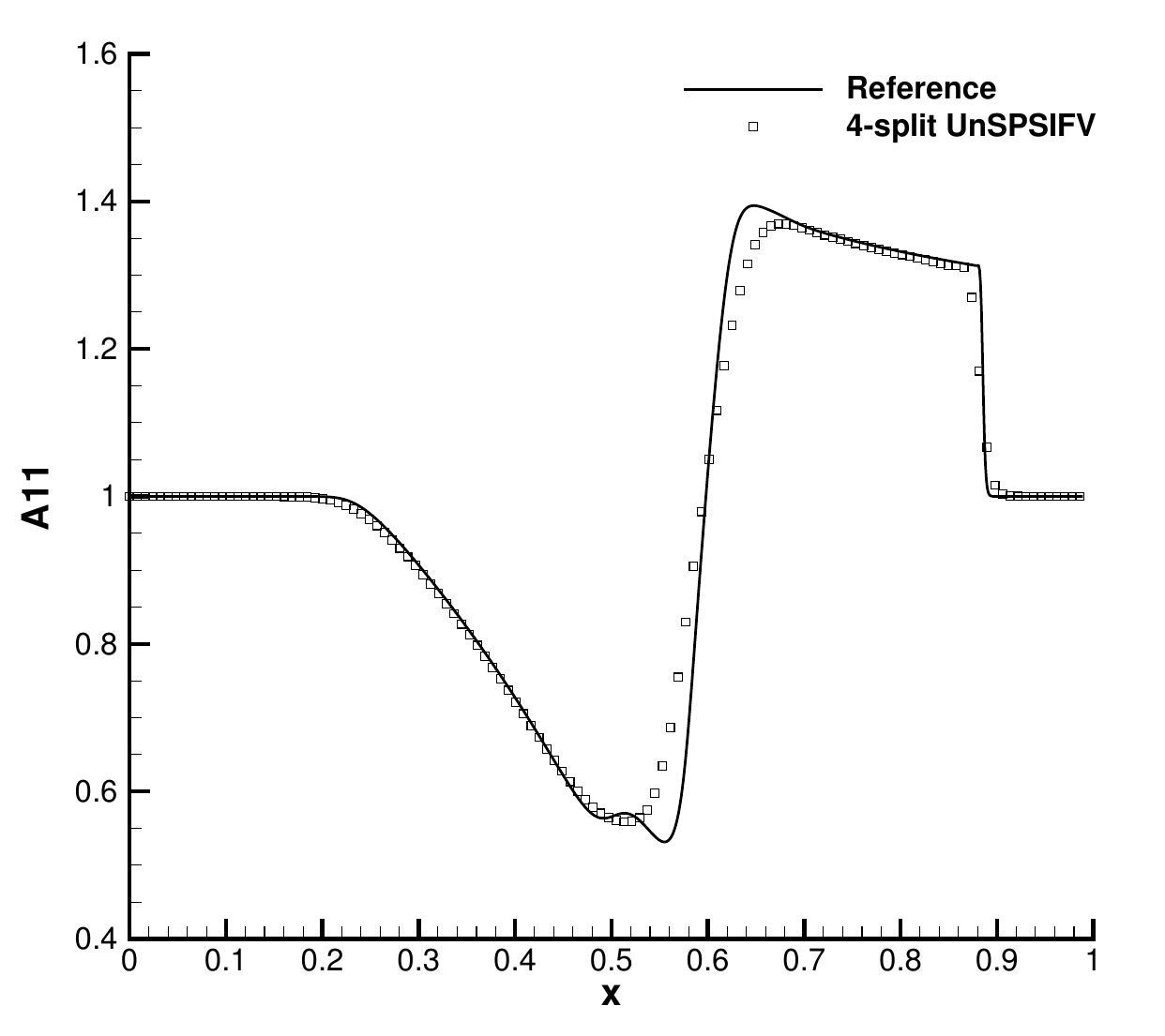}
        \end{tabular}
        \caption{Cut along the $x$-axis for the 2D explosion problem EP2 at time $t=0.15$. Comparison of a fine-grid reference solution computed with the thermodynamically compatible HTC scheme of \cite{HTCGPR,HTCAbgrall} (solid line) against the numerical solution of the homogeneous GPR model (solid limit, $\tau_1=\tau_2=10^{20}$)  obtained with the new four-split scheme (square symbols). Density $\rho$ (top left),
            velocity component $u$ (top right),  pressure $p$ (bottom left) and distortion field component $A_{11}$ (bottom right).   }
        \label{fig.ep2dsolid}
    \end{center}
\end{figure}

As a final test case we consider two radial explosion problems in the quasi-inviscid fluid limit of the model and the solid limit, denoted by EP1 and EP2, respectively.
The initial condition is given on the computational domain $\Omega=[-1,1]^2$ and reads in dependence of an inner and outer state $\mathbf{Q}_{in}$ and $\mathbf{Q}_{out}$ as follows
\begin{equation}
    \mathbf{Q}(x,y,0)=\begin{cases}
        \mathbf{Q}_{in} \quad \text{if} \quad r\leq R \\
        \mathbf{Q}_{out} \quad \text{if} \quad r>R. \\
    \end{cases}
\end{equation}
Therein, $r=\sqrt{x^2+y^2}$ is the radial coordinate and the initial discontinuity is placed at a radius of $R=0.5$.
\paragraph{EP1: Fluid limit of the model}
The quasi-inviscid fluid limit of the model is achieved by setting $\tau_1 \ll 1$ and $\tau_2 \ll 1$.
The inner and outer states are given by the left and right states of the 1D Riemann problem RP1, i.e.
$\rho_{in}=1$, $p_{in}=1$, $\rho_{out}=0.125$ and $p_{out}=0.1$ with a global zero velocity field $\mathbf{v}=0$. Accordingly, the thermal impulse is initially set to zero, i.e. $\J=0$, and the distortion field is initialized as $\A=\mathbf{I}$.
The remaining parameters are set to  $\gamma=1.4$, $c_v=2.5$, $c_s=1$, $c_h=1$, $\rho_0=1$, $\tau_1 = 10^{-8}$ and $\tau_2 = 10^{-10}$.
These settings, together with the choice of the initial states, lead to all Mach numbers being of the order of unity.
The final time of the simulation is $t=0.2$ and the grid is made of $146664$ cells.
As reference solution we solve the associated 1D radial Euler equations with source terms using a classical second-order MUSCL-Hancock scheme. For details on its computation see \cite{toro-book,SIGPR}.
The numerical results obtained with the new scheme are given in Figure \ref{fig.ep2dfluidsolid}. The 3D density contours are depicted in the left panel, and the comparison with the 1D reference solution is shown in Figure \ref{fig.ep2dfluid}.
The new scheme captures all waves accurately and shows a good agreement with the reference solution.
\paragraph{EP2: Solid limit of the model}
Analogously, we now solve the model with the same initial condition in the solid limit which is achieved by setting
$\tau_1 = \tau_2 = 10^{20}$ and using the same parameters as above for EP1 otherwise.
The simulation is carried out on the same computational grid up to a final time of $t=0.15$.
The reference solution is now obtained with the thermodynamically compatible HTC scheme presented in \cite{HTCGPR,HTCAbgrall} to which we compare the numerical solution of the new scheme.
In the right panel of Figure  \ref{fig.ep2dfluidsolid} the 3D density contour plot is shown while a comparison along 1D cuts in the $x$-direction is given in Figure \ref{fig.ep2dsolid}, showing good agreement between the two numerical solutions.

% % % % % % % % % % % % % % % % % % % % % % % % % % % % % %
%                  Conclusions                            %
% % % % % % % % % % % % % % % % % % % % % % % % % % % % % %

\section{Conclusions}\label{sec.conclusions}
This work introduced a novel structure-preserving semi-implicit four-split finite volume scheme on unstructured meshes for the numerical approximation of the unified first-order continuum mechanics model proposed in \cite{PeshRom2014}. The governing equations are decomposed into four physically meaningful subsystems associated with transport, thermal effects, material deformation, and pressure dynamics. While the transport step is advanced explicitly, the remaining contributions are integrated implicitly through a sequential procedure. As a result, the time-step restriction depends solely on the material velocity, eliminating the severe limitations imposed by acoustic, shear, and thermal wave propagation. This feature makes the method particularly attractive for all-Mach-number applications and for solid mechanics problems characterized by large elastic wave speeds. Furthermore, the scheme naturally recovers the asymptotic behavior of the model in the stiff relaxation limit, and thus reproduces, at the discrete level, the compressible Navier-Stokes stress and the Fourier heat flux. At the same time, the pressure update is designed to remain consistent with the incompressible regime, yielding the correct low Mach number behavior. A key ingredient of the method is the use of mimetic discrete operators on vertex-staggered unstructured triangular grids, which preserve the curl involution constraints associated with the distortion field and thermal impulse whenever the source terms are linear or absent, as confirmed numerically by the curl errors of $\A$ and $\J$, which remain at the level of machine precision throughout the simulation, see Figure \ref{fig.solidrotor.cutcurl}. The nonlinear transport part is instead handled by a robust finite-volume shock-capturing discretization, allowing the method to accurately resolve smooth solutions and to robustly capture discontinuities. The performance of the new algorithm has been assessed through a large suite of numerical experiments involving fluids and solids across different physical regimes. The reported results demonstrate that the method remains accurate and robust over a wide spectrum of Mach numbers, ranging from nearly incompressible flows to highly compressible dynamics, and from inviscid fluids to elastic materials.

Future developments will focus on enriching the structure-preserving properties of the proposed method, in particular by incorporating thermodynamic compatibility in the spirit of \cite{HTCGPR,HTCAbgrall,BoscheriGPRGCL,HTCLagrange,HTCLagrangeGPR}. Another important direction concerns the extension of the framework to coupled electro-magneto-mechanical phenomena, following the unified approach outlined in \cite{GPRmodelMHD}. Furthermore, the construction of high-order accurate discretizations will be investigated, building upon the compatible discretization techniques recently proposed in \cite{CompatibleDG1} as well as on advanced semi-implicit IMEX time-integration strategies developed in \cite{Spiteri3Split,BosFil2016}.

% % % % % % % % % % % % % % % % % % % % % % % % % % % % % %
%                Acknowledgment                          %
% % % % % % % % % % % % % % % % % % % % % % % % % % % % % %

\section*{Acknowledgments}

	M.D. was funded by the Fondazione Caritro via the project SOPHOS and by the European Research Council (ERC) under the European Union's Horizon Europe research and innovation programme via the project SOPHOS, grant agreement no ERC-ADG-2025-101265878-SOPHOS. Views and opinions expressed are however those of the authors only and do not necessarily reflect those of the European Union or the European Research Council Executive Agency. Neither the European Union nor the granting authority can be held responsible for them.
	W.B. received financial support from the "Institut des Math\'ematiques pour la Plan\`ete Terre" (France).
	A.T. acknowledges the financial support of the Agence Nationale de la Recherche (ANR) via the project DELFIN, project no. ANR-25-CE46-7729.
	M.D. and A.T. were also financially supported by the Italian Ministry of University
	and Research (MUR) via the  Departments of Excellence  Initiative 2018--2027 attributed to DICAM of the University of Trento (grant L. 232/2016).

    \vspace{3mm}

	All authors are members of the Gruppo Nazionale per il Calcolo Scientifico dell'Istituto Nazionale di Alta Matematica (GNCS-INdAM).

\section*{In memoriam}

This paper is dedicated to the memory of Professor Philip Lawrence Roe (4 May 1938 -- 26 April 2026), whom we all remember as a brilliant scientist and one of the fathers of our research field of numerical analysis for hyperbolic conservation laws, but, most importantly, also as an exceptional human being and dear friend, whom we would like to thank infinitely for all the deep and very inspiring discussions we had with him.  

%
%
%
%
%% % % % % % % % % % % % % % % % % % % % % % % % % % % % % %
%% % % % % % % % % % % % % % % % % % % % % % % % % % % % % %
%%              Bibliography
%% % % % % % % % % % % % % % % % % % % % % % % % % % % % % %
%% % % % % % % % % % % % % % % % % % % % % % % % % % % % % %
\bibliographystyle{elsarticle-num}
\bibliography{biblio}

% % % % % % % % % % % % % % % % % % % % % % % % % % % % % %
% % % % % % % % % % % % % % % % % % % % % % % % % % % % % %
%                   Appendix                              %
% % % % % % % % % % % % % % % % % % % % % % % % % % % % % %
% % % % % % % % % % % % % % % % % % % % % % % % % % % % % %
%\appendix

\end{document}